\documentclass[11pt]{amsart}
\usepackage[all]{xy}
\usepackage{yfonts}
\usepackage{amsmath}
\usepackage{amsfonts}
\usepackage{amssymb}
\usepackage{amsthm}
\usepackage{mathrsfs}
\usepackage{graphicx}
\usepackage{longtable}
\usepackage{tikz}
\usetikzlibrary{positioning}
\usepackage[
bookmarks=false,
colorlinks=true,
debug=true,
naturalnames=true,
pdfnewwindow=true,
citecolor=blue,
linkcolor=blue,
urlcolor = blue]{hyperref}
\usepackage{mathtools}
\numberwithin{equation}{section}

\newtheorem{Thm}{Theorem}[section]
\newtheorem{Cor}[Thm]{Corollary}
\newtheorem{Prop}[Thm]{Proposition}
\newtheorem{Lem}[Thm]{Lemma}
\newtheorem{Conj}[Thm]{Conjecture}

\theoremstyle{definition}
\newtheorem{Def}[Thm]{Definition}
\newtheorem{Rem}[Thm]{Remark}
\newtheorem{Ex}[Thm]{Example}

\newcommand{\Irr}{\mathop{\mathrm{Irr}}\nolimits}

\newcommand{\Z}{\mathbb{Z}}
\newcommand{\Q}{\mathbb{Q}}
\newcommand{\C}{\mathbb{C}}
\newcommand{\kk}{\Bbbk}

\newcommand{\id}{\mathrm{id}}

\title[Monoidal analog of cluster-tilting theory]
{Monoidal analog of cluster-tilting theory and monoidal reachability of acyclic cluster algebras}

\date{}

\author[H.~Sakamoto]{Heizo Sakamoto}
\address{Graduate School of Mathematical Sciences, the University of Tokyo, 
3-8-1, Komaba, Meguro-ku, Tokyo, 153-8914, Japan}
\email{sakaheisw0327@g.ecc.u-tokyo.ac.jp}

\begin{document}
\begin{abstract}
In this paper, we prove that many analogous facts of cluster-tilting theory holds in 
monoidal categories from quantum affine algebras or quiver Hecke algebras. 
We show that a maximal commuting family of real prime modules, 
which is an analogous object of a cluster-tilting object, 
is always a laurent family 
introduced by Kashiwara-Kim-Oh-Park. 
By using $\widetilde{\Lambda}$-invariants, we construct an exchange matrix of each laurent family 
and show that a laurent family is always completely $\Lambda$-admissible monoidal seed. 
We obtain that a pair of real prime modules is an exchange pair
if and only if $\mathfrak{d}$-invariant is $1$. 

When the category is a monoidal categorification of a cluster algebra of acyclic type, 
we prove that $\mathfrak{d}$-invariant is not less than generic Ext invariant by using 
tropical $F$-invariant introduced by Cao and Fei. 
Moreover, we show that, the conjecture asserting that real simple modules correspond 
exactly to cluster monomials holds in these cases. 
\end{abstract} 

\maketitle

\tableofcontents

\section{Introduction}\label{sec : intro}
Quantum groups, introduced by Drinfeld ~\cite{Dr85} and Jimbo~\cite{Ji85}, 
are $q$-deformations of universal enveloping algebras. 
In particular, when the Lie algebra is taken to be an affine Lie algebra, 
the corresponding quantum group is called a \emph{quantum affine algebra}. 
The study of finite-dimensional representations of quantum affine algebras 
has been an active area of research since the 1990s, 
owing to its significant applications not only in mathematics but also in physics, 
particularly in the theory of integrable systems and solvable lattice models
(see, for example, \cite{Baxter, FR92, DFJMN}). 
For a finite-dimensional simple Lie algebra $\mathfrak{g}$, 
let $\mathscr{C}_{\mathfrak{g}}$ be the category of finite-dimensional representations of 
the quantum affine algebra $U_q(\hat{\mathfrak{g}})$. 
In \cite{KKOP24}, Kashiwara-Kim-Oh-Park introduced a full subcategory 
$\mathscr{C}_{\mathfrak{g}}^{[a, b], \mathfrak{s}}$ 
of $\mathscr{C}_{\mathfrak{g}}$ and proved that it is a monoidal categorification of a cluster algebras, 
which is introduced in \cite{HL10}. 

Quiver Hecke algebras, introduced by Khovanov-Lauda and Rouquier~\cite{KL09, Ro08}, 
is infinite dimensional graded algebras associated with root data. 
Let $R$ be a symmetric quiver Hecke algebra over a field of characteristic $0$. 
Let $R$-gmod be a category of finite-dimensional graded module. 
It categorifies a quantum unipotent coordinate ring. 
For $w$ in the Weyl group $W$, 
the full subcategory $\mathscr{C}_{w}$ is introduced in \cite{KKKO18} and 
it is proved that $\mathscr{C}_{w}$ categorify a cluster algebra. 

Let $\mathscr{C}$ be $\mathscr{C}_{\mathfrak{g}}^{[a, b], \mathfrak{s}}$ or $\mathscr{C}_{w}$. 
Let $N = b-a+1$ for quantum affine case and $N = (\text{the length of }w)$ for quver Hecke case. 
In what follows, for $L, M \in \mathscr{C}$, even if $\mathscr{C} = \mathscr{C}_w$, 
we write the tensor product or the convolution product of  $L$ and $M$ as $L\otimes M$. 
Also, we ignore the grading shifts when $\mathscr{C} = \mathscr{C}_w$. 
For $L \in \Irr{\mathscr{C}}$, we call $L$ \emph{real} if $L\otimes L \in \Irr{\mathscr{C}}$ and 
\emph{imaginary} otherwise. 
The real and imaginary modules is really important issue to study the structure of 
the tensor product of two simple modules. 

To classify real and imaginary modules, the following conjecture is essential : 
\begin{Conj}
\label{Conj : real vs monomial intro}
    Let $\mathcal{C}$ be a monoidal categorification of a cluster algebra $\mathscr{A}$ 
    and let $\iota : K(\mathcal{C})_{\Z} \overset{\sim}{\to} \mathscr{A}$ be the corresponding 
    isomorphism. 
    Then, there is a bijection between the set of real simple modules 
    and the set of cluster monomials, i.e., 
    \[
    \{ \iota([L]) \mid L  \text{ is a real simple module in $\mathcal{C}$}\} 
    = \{ \text{cluster monomial in $\mathscr{A}$}\}. 
    \]
\end{Conj}
In this paper, we consider a monoidal analog of cluster-tilting theory and 
prove this conjecture when $\mathscr{C}$ is a monoidal categorification of a cluster algebra of acyclic type. 

Let $Q$ be a finite quiver and $\mathcal{C}_Q$ be the (generalized) cluster category. 
In the cluster-tilting theory, 
the following facts are known by \cite{BMRRT, IY08, DK08, FG} : 
\begin{enumerate}
  \item[(CT1)] For a basic cluster-tilting object $T$, $|T| = |Q_0|$, 
  \item[(CT2)] for an almost complete basic cluster-tilting object $\overline{T}$, 
  there exists exactly two indecomposable rigid objects which are complement of $\overline{T}$, 
  \item[(CT3)] $(M, M')$ is an exchange pair if and only if $\operatorname{Ext}^1_{\mathcal{C}_Q}(M, M') =1$, 
  \item[(CT4)] if $Q$ is acyclic, then any basic cluster-tilting object is reachable. 
\end{enumerate}
It is natural to ask that these fact also hold in the monoidal settings.

First, we introduce \emph{maximal real set} in $\mathscr{C}$(Definition~\ref{Def : maximal real set}), 
which is a monoidal analog of a maximal rigid set. 
We call a family $\mathcal{M}$ of simple modules in $\mathscr{C}$ a maximal real set 
if it is a commuting family of real prime modules and it is maximal in such families. 
A maximal real set $\mathcal{M}$ is also a monoidal analog of cluster-tilting object 
by the following property(Proposition~\ref{Prop : simple modules commuting with a maximal real set}) : 
 if $X \in \Irr \mathscr{C}$ is commute with any element of $\mathcal{M}$, 
then $X$ is a monomial in $\mathcal{M}$. 
We prove that $\mathcal{M}$ is a maximal real set if and only if it is a Laurent family introduced in 
\cite{Laurentfamily} and prove the following corresponding to (CT1) : 
\begin{Thm}[Theorem~\ref{Thm : cardinality of maximal real sets}]
  Any maximal real set has exactly $N$ elements. 
\end{Thm}
By using properties of Laurent families in \cite{KK19, Laurentfamily}, 
we prove the following corresponding to (CT3) : 
\begin{Thm}[Theorem~\ref{Thm : characterization of exchange pairs}]
Let $L, L' \in \Irr\mathscr{C}$ be non-isomorphic real prime modules.
Then $(L, L')$ is an exchange pair if and only if
$\mathfrak{d}(L, L') = 1$.
\end{Thm}
For a maximal real set $\mathcal{M}$, 
we define the matrix $\widetilde B_{\mathcal M}$ using $\widetilde{\Lambda}$-invariant and 
prove the following corresponding to (CT2) : 
\begin{Thm}[Corollary~\ref{Cor : complete Lambda admissibility of maximal real sets}]
  A pair $(\mathcal{M}, \widetilde{B}_{\mathcal{M}})$ is a completely
$\Lambda$-admissible monoidal seed.
\end{Thm}
Also, for a real prime module $L \in \Irr\mathscr{C}$, 
we prove that $[L] \in K(\mathscr{C})_{\Z}$ is a prime element
(Proposition~\ref{Prop : real prime modules give prime elements}). 
By the property of $\mathscr{C}_{\mathfrak{g}}^{[a, b], \mathfrak{s}}$, 
we have the following : 
\begin{Thm}[Corollary~\ref{Cor : unique prime tensor factorization in Cg}]
For a real module $L \in \mathscr{C}_{\mathfrak{g}}$,
tensor factorization of $L$ into prime modules
is unique up to permutation.
\end{Thm}

For the cluster category $\mathcal{C}_{Q}$ and $\mathbf{g} \in \Z^Q$, 
there is a \emph{canonical decomposition} $\mathbf{g} = \mathbf{g}_1 \oplus \cdots\oplus \mathbf{g}_s$, 
which express the decomposition of generic presenintation in $\operatorname{PHom}(\mathbf{g})$ 
into indecomposable presentations. 
Viewing decomposition into indecomposable objects in a cluster category 
as the counterpart of prime factorization on the monoidal side, 
we want to transfer the structure of canonical decompositions to the monoidal setting. 
In the cluster category, $e$-invariant is introduced in \cite{DF15}. 
Since $\mathbf{g}=\mathbf{g}_1\oplus\cdots\oplus\mathbf{g}_s$ is the canonical decomposition 
if and only if $e(\mathbf{g}_i,\mathbf{g}_j)=0$ for all $i\neq j$, 
it is natural to compare the $\mathfrak d$-invariant and $e$-invariant. 
In \cite{Cao, Fei23}, 
the \emph{tropical $F$-invariant} is introduced and they compared the tropical $F$-invariant and 
$\mathfrak{d}, e$-invariants. 
In \cite{Cao}, it is proved that 
$\mathfrak{d}(L_1, L_2) = F_{\iota([L_1])}^t[(I\ 0)\deg^t\iota([L_2])]
+ F_{\iota([L_2])}^t[(I\ 0)\deg^t\iota([L_1])]$ if at least one of $L_1, L_2$ is reachable. 
We prove, without any assumptions, 
$\mathfrak{d}(L_1, L_2) \geq F_{\iota([L_1])}^t[(I\ 0)\deg^t\iota([L_2])]
+ F_{\iota([L_2])}^t[(I\ 0)\deg^t\iota([L_1])]$
 for $L_1, L_2 \in \Irr{\mathscr{C}}$
 (Theorem~\ref{Thm : tropical F invariants and Lambda inequalities}). 

Let consider the case that the cluster algebra which is categorified by $\mathscr{C}$ is acyclic type. 
For the case that  $Q$ is acyclic, 
it is proved in \cite{Fei23} that 
$F_{\mathbf{g}_1}^{\mathrm{gen}}[\mathbf{g}_2] = e(\mathbf{g}_1, \mathbf{g}_2)$ 
for $\mathbf{g}_1, \mathbf{g}_2\in \Z^{Q_0}$, 
where $F_{\mathbf{g}}^{\mathrm{gen}}$ is the generic $F$-polynomial. 
By comparing the Newton polytope of $F$-polynomial from simple module and generic $F$-polynomial, 
we obtain the following inequality : 
\begin{Thm}[Corollary~\ref{Cor : d dominates generic extensions}]
  Let $L_1, L_2 \in \Irr\mathscr{C}$.
Let $\widehat{\mathbf{g}}_i
\coloneq \deg^{t_1}\iota([L_i])$ and
$\mathbf{g}_i
\coloneq \operatorname{pr}_{\mathrm{uf}}\widehat{\mathbf{g}}_i$.
Then we have
\[
\mathfrak{d}(L_1, L_2)
\geq e(\mathbf{g}_1, \mathbf{g}_2)
+ e(\mathbf{g}_2, \mathbf{g}_1).
\]
\end{Thm}

As a result, if $\mathfrak{d} = 0$, then $e = 0$. 
By using this property, we prove the following corresponding to (CT4) : 
\begin{Thm}[Theorem~\ref{Thm : monoidal reachability for acyclic cluster algebras}]
Assume that the cluster algebra categorified by $\mathscr{C}$ is of acyclic type. 
For $L \in \Irr\mathscr{C}$,
$L$ is real if and only if $L$ is reachable.
Moreover, $L$ is real prime if and only if
$L$ corresponds to a cluster variable.
\end{Thm}

\subsection*{Acknowledgement}
I would like to express my sincere gratitude to my supervisor, 
Noriyuki Abe, for his continuous support and helpful discussions. 
Also I am deeply grateful to Hironori Oya for taking the time to discuss this work. 
This work was supported by JSPS KAKENHI Grant Number 26KJ0951. 

\subsection*{On the use of AI}
The author used Chat GPT to search related paper, check the draft and correctness of the argument, 
and improve the expression. 
All the ideas in this paper are the author's own. 
The author takes full responsibility for its content.

\section{Quantum loop algebras and their finite-dimensional representations}
\label{Sec : qloop}

\subsection{Quantum loop algebras}
\label{ssec : quantum loop algebras}

Let $\mathfrak{g}$ be a complex finite-dimensional simple Lie algebra
with the set of Dynkin indices $I = \{1, \ldots, n\}$
and Cartan matrix $C = (c_{ij})_{i, j \in I}$.
Let $r$ be the maximum number of edges in the Dynkin diagram
and $D = \mathrm{diag}(d_i \mid i \in I)$ be the unique diagonal matrix
such that $d_i \in \{1, r\}$ for each $i \in I$
and $DC = (d_i c_{ij})_{i, j \in I}$ is symmetric.

Let $q$ be an indeterminate and $\kk$
the algebraic closure of the rational function field $\Q(q)$.
Let $q_i \coloneq q^{d_i}$.
For $k \in \Z$, we set 
\[
[k]_{q_i} \coloneq \frac{q_i^k - q_i^{-k}}{q_i - q_i^{-1}}. 
\]
For integers $k \geq l \geq 0$, we set
\[
[k]_{q_i}! \coloneq [k]_{q_i}[k - 1]_{q_i}\cdots [2]_{q_i}[1]_{q_i},
\quad
\left[\begin{matrix} k \\ l \end{matrix}\right]_{q_i}
\coloneq \frac{[k]_{q_i}!}{[k - l]_{q_i}! [l]_{q_i}!}.
\]

\begin{Def}
\label{Def : quantum loop algebra}
\emph{The quantum loop algebra} $U_q(L\mathfrak{g})$
is the $\kk$-algebra given by the following generators and relations :

generators : \quad
$x_{i, k}^{\pm}\ (i \in I, k \in \Z), \quad
k_i^{\pm 1}\ (i \in I), \quad
h_{i, l}\ (i \in I, l \in \Z \setminus \{0\})$

relations :
\begin{itemize}
    \item $k_i k_i^{-1} = k_i^{-1} k_i = 1,
    k_i k_j = k_j k_i$ for $i, j \in I$,
    \item $k_i x_{j, k}^{\pm} k_i^{-1}
    = q_i^{\pm c_{ij}} x_{j, k}^{\pm}$
    for $i, j \in I$ and $k \in \Z$,
    \item $[k_i, h_{j, l}] = [h_{i, l}, h_{j, m}] = 0$
    for $i, j \in I$ and $l, m \in \Z \setminus \{0\}$,
    \item $[h_{i, l}, x_{j, m}^{\pm}]
    = \pm \frac{1}{l}[l c_{ij}]_{q_i}x_{j, l + m}^{\pm}$
    for $i, j \in I, l \in \Z \setminus \{0\}$ and $m \in \Z$,
    \item $x_{i, k + 1}^{\pm} x_{j, l}^{\pm}
    - q_i^{\pm c_{ij}} x_{j, l}^{\pm} x_{i, k + 1}^{\pm}
    = q_i^{\pm c_{ij}} x_{i, k}^{\pm} x_{j, l + 1}^{\pm}
    - x_{j, l + 1}^{\pm} x_{i, k}^{\pm}$
    for $i, j \in I$ and $k, l \in \Z$,
    \item $[x_{i, k}^{+}, x_{j, l}^{-}]
    = \displaystyle\delta_{i, j}
    \frac{\phi_{i, k + l}^{+} - \phi_{i, k + l}^{-}}{q_i - q_i^{-1}}$
    for $i, j \in I$ and $k, l \in \Z$,
    \item $\displaystyle\sum_{\sigma \in \mathfrak{S}_b}\sum_{k = 0}^{b}
    (-1)^k \left[\begin{matrix} b \\ k \end{matrix}\right]_{q_i}
    x_{i, r_{\sigma(1)}}^{\pm}\cdots x_{i, r_{\sigma(k)}}^{\pm}
    x_{j, s}^{\pm} x_{i, r_{\sigma(k + 1)}}^{\pm}\cdots
    x_{i, r_{\sigma(b)}}^{\pm} = 0$,
    where $b = 1 - c_{ij}$, for $i, j \in I$ with $i \neq j$
    and $s, r_1, \ldots, r_b \in \Z$,
\end{itemize}
where
\[
\phi_i^{\pm}(z) \coloneq
\sum_{r = 0}^{\infty}\phi_{i, \pm r}^{\pm}z^{\pm r} \coloneq
k_i^{\pm 1}\exp\left(\pm (q_i - q_i^{-1})
\sum_{l = 1}^{\infty}h_{i, \pm l}z^{\pm l}\right).
\]
\end{Def}

\subsection{Finite-dimensional representations}
\label{ssec : finite-dimensional representations}

A finite-dimensional $U_q(L\mathfrak{g})$-module is said to be of type
$\mathbf{1}$
if the element $k_i$ acts as a diagonal $\kk$-linear operator
whose eigenvalues belong to the set $\{q^k \mid k \in \Z\}$
for each $i \in I$.
It is well-known that the study of finite-dimensional representations
of $U_q(L\mathfrak{g})$ reduces essentially to the study of the category
$\mathscr{C}_{\mathfrak{g}}$ of type $\mathbf{1}$ finite-dimensional
$U_q(L\mathfrak{g})$-modules.
It is also well known that $U_q(L\mathfrak{g})$ is a subquotient of the
corresponding
quantum affine algebra and the study of type $\mathbf{1}$
finite-dimensional modules of quantum affine algebra reduces
to the study of $\mathscr{C}_{\mathfrak{g}}$.
Since $U_q(L\mathfrak{g})$ has a Hopf algebra structure,
$\mathscr{C}_{\mathfrak{g}}$ is a monoidal category.

Let $V \in \mathscr{C}_{\mathfrak{g}}$.
Since the elements
$\{k_i^{\pm 1}, h_{i, l} \mid i \in I, l \in \Z \setminus \{0\}\}$
mutually commute, there is a decomposition
\[
V = \bigoplus_{\gamma \in
(\kk[\![z]\!] \times \kk[\![z^{-1}]\!])^I}V_{\gamma}
\]
where, for each
$\gamma = (\gamma_i^{+}(z), \gamma_i^{-}(z))_{i \in I}$,
we define $V_{\gamma}$ as the subspace of $V$ on which each coefficient
of the series $\phi_i^{\pm}(z) - \gamma_i^{\pm}(z)$ acts nilpotently
for any $i \in I$.
If $V_{\gamma} \neq 0$, it is called \emph{an $\ell$-weight space}
and $\gamma$ is called \emph{$\ell$-weight}.

For a simple module $L$ in $\mathscr{C}_{\mathfrak{g}}$,
there uniquely exists an $\ell$-weight space $L_{\gamma_0}$
such that $x_{i, k}^{+}L_{\gamma_0} = 0$
for all $i \in I$ and $k \in \Z$.
Such $\ell$-weight $\gamma_0$ is called
\emph{the highest $\ell$-weight} of $L$. 
The isomorphism classes of simple objects
in $\mathscr{C}_{\mathfrak{g}}$ is classified by highest $\ell$-weights :

\begin{Thm}[{\cite[Theorem 3.3]{CP95}, \cite[Theorem 12.2.6]{CP94}}]
\label{Thm : Drinfeld poly}
If $\gamma = (\gamma_i^{+}(z), \gamma_i^{-}(z))_{i \in I}$
is an highest $\ell$-weight of a simple $U_q(L\mathfrak{g})$-module
in $\mathscr{C}_{\mathfrak{g}}$, for each $i \in I$ there exists a unique
polynomial
$P_i(z) \in 1 + z\kk[z]$ for $i \in I$ such that
\begin{equation}
\label{eq : ellhwt}
\gamma_i^{\pm}(z) = q_i^{\deg(P_i)}
\frac{P_i(zq_i^{-1})}{P_i(zq_i)}.
\end{equation}
Conversely, for any $(P_i(z))_{i \in I} \in (1 + z\kk[z])^I$,
we have a simple $U_q(L\mathfrak{g})$-module in $\mathscr{C}_{\mathfrak{g}}$
whose highest $\ell$-weight $\gamma$ is given by
{\rm (\ref{eq : ellhwt})}.
\end{Thm}

This $I$-tuple of polynomials $(P_i(z))_{i \in I}$
associated with a simple module $L$
is called \emph{the Drinfeld polynomials of $L$}.

Let $\Irr\mathscr{C}_{\mathfrak{g}}$ denote the classes of simple objects of
$\mathscr{C}_{\mathfrak{g}}$.
By theorem~\ref{Thm : Drinfeld poly}, $\Irr\mathscr{C}_{\mathfrak{g}}$
is parametrized by $(1 + z\kk[z])^I$.

Let $Y_{i, a}\ (i \in I, a \in \kk^{\times})$ be an indeterminate and
\[
\mathcal{M}_{+} \coloneq
\biggl\{m = \prod_{i \in I, a \in \kk^{\times}}Y_{i, a}^{u_{i, a}}
\ \bigg|\ u_{i, a} \in \Z_{\geq 0},
u_{i, a} = 0\ \text{except for finitely many }(i, a)\biggr\}.
\]
Then, there is a bijection
\[
\mathcal{M}_{+} \to (1 + z\kk[z])^I; \quad
\prod_{i \in I, a \in \kk^{\times}}Y_{i, a}^{u_{i, a}}
\mapsto \left(\prod_{a \in \kk^{\times}}(1 - az)^{u_{i, a}}\right)_{i \in I}.
\]
Using this bijection, $\Irr\mathscr{C}_{\mathfrak{g}}$ is also parametrized
by $\mathcal{M}_{+}$.
Let $L(m)$ denote the simple module corresponding to the monomial $m$.
The monomial $m$ is called the \emph{highest $\ell$-weight monomial}
of the simple module $L(m)$.

\begin{Prop}[{\cite[Corollary 3.5]{CP95}}]
\label{Prop : subquot}
Let $m_1, \ldots, m_k \in \mathcal{M}_{+}$.
Then, $L(m_1\cdots m_k)$ is isomorphic to a subquotient
of $\bigotimes_{i = 1}^{k}L(m_i)$.
\end{Prop}

For each $a \in \kk^{\times}$, we consider a $\kk$-algebra automorphism
$\tau_a$ of $U_q(L_{\mathfrak{g}})$ given by
\[
\tau_a(k_i) = k_i, \quad
\tau_a(x_{i, k}^{\pm}) = a^k x_{i, k}^{\pm}, \quad
\tau_a(h_{i, l}) = a^l h_{i, l},
\]
for $i \in I, k \in \Z, l \in \Z \setminus \{0\}$.
This automorphism $\tau_a$ is a $q$-analog of the loop rotation
$z \mapsto az$ for the loop algebra $L\mathfrak{g} = \mathfrak{g}
\otimes_{\C}\C[z^{\pm}]$.
For a finite-dimensional representation $V$ in $\mathscr{C}_{\mathfrak{g}}$,
the twisted representation $\tau_a^*V$ by $\tau_a$
again belongs to $\mathscr{C}$.
Thus the assignment $V \mapsto \tau_a^*V$ defines a monoidal
auto-equivalence $T_a$ of $\mathscr{C}_{\mathfrak{g}}$.
If we denote by $m \mapsto m_a$ the monoid automorphism
of $\mathcal{M}_{+}$ given by $Y_{i, b} \mapsto Y_{i, ab}$,
we have $T_a(L(m)) \cong L(m_a)$
for any simple module $L(m) \in \mathscr{C}_{\mathfrak{g}}$.

In $\mathscr{C}_{\mathfrak{g}}$,
$V \otimes W$ is not necessarily isomorphic to $W \otimes V$
for $V, W \in \mathscr{C}_{\mathfrak{g}}$.
However, using the $q$-character,
which is a type of character in the quantum loop algebra,
we have the following fact.

\begin{Prop}[{\cite[Theorem 3]{FR98}}]
\label{Prop : q-character}
Let $K(\mathscr{C}_{\mathfrak{g}})$ be the Grothendieck ring
of the monoidal category $\mathscr{C}_{\mathfrak{g}}$.
Then, there is an injective algebra homomorphism
$\chi_q \colon K(\mathscr{C}_{\mathfrak{g}}) \to
\Z[Y_{i, a}^{\pm 1} \mid i \in I, a \in \kk^{\times}]$.
In particular, although $\mathscr{C}_{\mathfrak{g}}$ is not braided,
$K(\mathscr{C}_{\mathfrak{g}})$ is commutative.
\end{Prop}

\subsection{Some subcategories of $\mathscr{C}_{\mathfrak{g}}$}
\label{ssec : subcategories}

For $i, j \in I$, write
\[
d_{ij} \coloneq \min\{d_i, d_j\}.
\]
A function $\epsilon \colon I \to \{0, 1\}$
is called \emph{a parity function} for $\mathfrak{g}$
if it satisfies the condition
\[
\epsilon_i \equiv \epsilon_j + d_{ij} \pmod{2}
\]
for any adjacent nodes $i, j \in I$.
In what follows, we choose and fix a parity function $\epsilon$.
Let $\widehat{I} \coloneq
\{(i, p) \in I \times \Z \mid p - \epsilon_i \in 2\Z\}$.

\begin{Def}[{\cite[$\S$3.7]{HL10}, \cite[Definition 2.6]{FHOO22}}]
\label{Def : Cg0}
Let $\mathscr{C}_{\mathfrak{g}}^{0}$ be the smallest full subcategory of
$\mathscr{C}_{\mathfrak{g}}$
containing $L(Y_{i, q^p})$ for $(i, p) \in \widehat{I}$
and trivial module, and stable under taking subquotients,
extensions, and tensor products.
\end{Def}

\begin{Prop}[{\cite[$\S$3.7]{HL10}, \cite[Proposition 2.8]{FHOO22}}]
\label{Prop : Cg factorization}
Let $L \in \Irr\mathscr{C}_{\mathfrak{g}}$.
Then $L$ is factorized into a commutative tensor product
\[
L \cong T_{a_1}(L_1) \otimes \cdots \otimes T_{a_d}(L_d)
\]
with $L_k \in \Irr\mathscr{C}_{\mathfrak{g}}^{0}$ such that $L_k$ is non-trivial and $a_k \in \kk^{\times}$
for $1 \leq k \leq d$ such that
$a_k/a_l \notin q^{2\Z}$ for all $k \neq l$.
Moreover, the set $\{T_{a_k}(L_k)\}_k$ is unique. 
\end{Prop}

For a complex finite-dimensional simple Lie algebra $\mathfrak{g}$,
we associate a Lie algebra $\mathfrak{g}_{\mathrm{fin}}$
of simply laced type as follows :
\[
\begin{array}{c||c|c|c|c|c|c|c}
\mathfrak{g} & A_n\ (n \geq 1) & B_n\ (n \geq 2) & C_n\ (n \geq 3)
& D_n\ (n \geq 4) & E_{6, 7, 8} & F_4 & G_2 \\
\hline
\mathfrak{g}_{\mathrm{fin}} & A_n & A_{2n - 1} & D_{n + 1}
& D_n & E_{6, 7, 8} & E_6 & D_4
\end{array}
\]
Let $\Delta_{\mathfrak{g}_{\mathrm{fin}}}$ be the Dynkin diagram
for $\mathfrak{g}_{\mathrm{fin}}$.
Then there exists a Dynkin diagram automorphism $\sigma$
of $\Delta_{\mathfrak{g}_{\mathrm{fin}}}$ whose orbit set
$\Delta_{\mathfrak{g}_{\mathrm{fin}}}^{\sigma}$
yields the Dynkin diagram $\Delta_{\mathfrak{g}}$.
We denote by $I_{\mathrm{fin}} = \{1, 2, \ldots, |I_{\mathrm{fin}}|\}$
the index set of $\Delta_{\mathfrak{g}_{\mathrm{fin}}}$.
To distinguish the index set $I_{\mathrm{fin}}$ and $I$,
we use $\iota$ for indices in $I_{\mathrm{fin}}$
and $i$ for indices in $I$.
Let $\pi \colon I_{\mathrm{fin}} \to I$
be the projection induced by $\sigma$.
For $\iota \in I_{\mathrm{fin}}$,
we write $d_{\iota} \coloneq |\pi^{-1}\pi(\iota)|$.

\begin{Def}[{\cite[Definition 6.9]{KKOP24}}]
\label{Def : admissible sequence}
We say that an infinite sequence
$\mathfrak{s} \coloneq ((\iota_k, p_k))_{k \in \Z}$
in $I_{\mathrm{fin}} \times \Z$ is \emph{admissible}
if the following conditions are satisfied :
\begin{enumerate}
    \item[(1)] $p_{s^{+}} = p_s + 2d_{\iota_s}$,
    \item[(2)] $p_t = p_s + \min\{d_{\iota_s}, d_{\iota_t}\}$
    if $\iota_s$ and $\iota_t$ are adjacent
    and $t^{-} < s < t < s^{+}$,
    \item[(3)] $p_t - (p_s + 2) \in 2d_{\iota_s}\Z$
    if $s, t \in \Z$ satisfy $\iota_t = \sigma(\iota_s)$,
    \item[(4)] $s_{\iota_k}s_{\iota_{k + 1}}\cdots
    s_{\iota_{k + l - 1}} = w_0$ for every $k \in \Z$,
    where $l$ denotes the length of the longest element $w_0$
    of the Weyl group of $\mathfrak{g}_{\mathrm{fin}}$.
\end{enumerate}
Here, for $s \in \Z$,
$s^{+} \coloneq \min\{t > s \mid \iota_t = \iota_s\}$
and $s^{-} \coloneq \max\{t < s \mid \iota_t = \iota_s\}$.
\end{Def}

In what follows, we only consider admissible sequences
$\mathfrak{s} = ((\iota_k, p_k))_{k \in \Z}$
which satisfies $p_k \equiv \epsilon_{\pi(\iota_k)} \pmod{2}$.

\begin{Def}[{\cite[Definition 6.16]{KKOP24}}]
\label{Def : category}
Let $[a, b]$ be an interval with $-\infty \leq a \leq b \leq \infty$.
Let $\mathfrak{s} = ((\iota_k, p_k))_{k \in \Z}
\subset I_{\mathrm{fin}} \times \Z$ be an admissible sequence.
We define the category $\mathscr{C}_{\mathfrak{g}}^{[a, b], \mathfrak{s}}$
as the smallest full subcategory of $\mathscr{C}_{\mathfrak{g}}^{0}$ which
contains
$S_k^{\mathfrak{s}} \coloneq L(Y_{\pi(\iota_k), q^{p_k}})$
$(k \in [a, b])$ and trivial module,
and is closed under taking subquotients, extensions, and tensor products.
\end{Def}
We call $S_k$ a \emph{cuspidal module}. 

\begin{Prop}[{\cite[$\S$6.3]{KKOP24}}]
\label{Prop : whole category}
For any admissible sequence $\mathfrak{s}$,
$\mathscr{C}_{\mathfrak{g}}^{[-\infty, \infty], \mathfrak{s}} =
\mathscr{C}_{\mathfrak{g}}^{0}$.
\end{Prop}

\begin{Prop}[{\cite[$\S$6.3]{PBW24}}]
\label{Prop : PBW parameterization}
Let $[a, b]$ be an interval and
$\mathfrak{s} = ((\iota_k, p_k))_{k \in \Z}$
be an admissible sequence.
\begin{enumerate}
    \item[(i)] For any sequence $(c_k)_k \in \Z_{\geq 0}^{[a, b] \cap \Z}$
    with $c_k = 0$ for all but finitely many $k$'s,
    the head of $S_b^{c_b} \otimes S_{b - 1}^{c_{b - 1}}
    \otimes \cdots \otimes S_a^{c_a}$
    is in $\Irr\mathscr{C}_{\mathfrak{g}}^{[a, b], \mathfrak{s}}$.
    Moreover, this correspondence form a bijection between
    \[
    \bigl\{(c_k)_k \in \Z_{\geq 0}^{[a, b] \cap \Z}
    \mid c_k = 0\ \text{for all but finitely many }k\text{'s}\bigr\}
    \]
    and $\Irr\mathscr{C}_{\mathfrak{g}}^{[a, b], \mathfrak{s}}$.
    \item[(ii)] The Grothendieck ring
    $K(\mathscr{C}_{\mathfrak{g}}^{[a, b], \mathfrak{s}})$ is the polynomial ring
    generated by $\{[S_s] \mid s \in [a, b]\}$
    and the set of highest $\ell$-weights
    of $\Irr\mathscr{C}_{\mathfrak{g}}^{[a, b], \mathfrak{s}}$ is the set of
    monomials
    in $Y_{\pi(\iota_k), q^{p_k}}$ $(k \in [a, b])$.
\end{enumerate}
\end{Prop}

For $L \in \Irr\mathscr{C}_{\mathfrak{g}}^{[a, b], \mathfrak{s}}$,
the corresponding sequence $(c_k)_k$ is called \emph{PBW parameter}
and we write it as $\mathrm{PBW}(L)$.

\section{Quiver Hecke algebras}
\label{Sec : klr}

\subsection{Quantum group and dual canonical basis}
\label{ssec : quantum group and dual canonical basis}

Let $\mathfrak{g}$ be a symmetric Kac--Moody Lie algebra
with Dynkin index set $I$.
Let $P$ and $P^\vee$ be its weight and coweight lattices,
and let $\{\alpha_i\}_{i \in I}$ and $\{h_i\}_{i \in I}$
be the sets of simple roots and simple coroots.
We set
\[
Q \coloneq \bigoplus_{i \in I}\Z\alpha_i, \qquad
Q_+ \coloneq \bigoplus_{i \in I}\Z_{\geq 0}\alpha_i, \qquad
Q_- \coloneq -Q_+.
\]
Let $W$ be the Weyl group of $\mathfrak{g}$,
and let $s_i$ denote the simple reflection corresponding to $i \in I$.

Let $U_q(\mathfrak{g})$ be the quantum group associated with
$\mathfrak{g}$, that is, the $\Q(q)$-algebra generated by
$e_i, f_i$ $(i \in I)$ and $q^h$ $(h \in P^\vee)$
with the usual defining relations.
We denote by $U_q^+(\mathfrak{g})$ and $U_q^-(\mathfrak{g})$
the subalgebras generated by $\{e_i \mid i \in I\}$
and $\{f_i \mid i \in I\}$, respectively.
For $n \in \Z_{\geq 0}$, set
\[
e_i^{(n)} \coloneq \frac{e_i^n}{[n]_{q_i}!}, \qquad
f_i^{(n)} \coloneq \frac{f_i^n}{[n]_{q_i}!}.
\]
For $\beta \in Q_-$, let $U_q^-(\mathfrak{g})_\beta$
be the homogeneous component of $U_q^-(\mathfrak{g})$
of weight $\beta$.

For each $i \in I$, consider Lusztig's braid group automorphism
$T_i : U_q(\mathfrak{g}) \to U_q(\mathfrak{g})$.
(See \cite{Lusbook}.
Note that, we write $v$ in \cite{Lusbook} as $q$
and adopt $T'_{i, -1}$ for $T_i$.)

Let $w \in W$ and fix a reduced expression
$\widetilde{w} = s_{i_1}\cdots s_{i_r}$.
For $1 \leq k \leq r$, set
\[
\beta_k \coloneq s_{i_1}\cdots s_{i_{k - 1}}(\alpha_{i_k})
\]
and, for $c \in \Z_{\geq 0}$, define
\[
F(c\beta_k) \coloneq T_{i_1}\cdots T_{i_{k - 1}}
\bigl(f_{i_k}^{(c)}\bigr).
\]
For $\mathbf{c} = (c_1, \ldots, c_r) \in \Z_{\geq 0}^r$, set
\[
F(\mathbf{c}, \widetilde{w}) \coloneq
F(c_r\beta_r)F(c_{r - 1}\beta_{r - 1})\cdots F(c_1\beta_1).
\]
It is well-known that the $\Q(q)$-linear subspace spanned by
$\{F(\mathbf{c}, \widetilde{w}) \mid
\mathbf{c} \in \Z_{\geq 0}^r\}$
does not depend on the choice of the reduced word $\widetilde{w}$.
We denote the subspace by $U_q^-(w)$.
Then $\{F(\mathbf{c}, \widetilde{w}) \mid
\mathbf{c} \in \Z_{\geq 0}^r\}$
is a $\Q(q)$-basis of $U_q^-(w)$ and we call it a PBW-basis.

There exists Kashiwara's bilinear form $(\ ,\ )_K$
on $U_q^-(\mathfrak{g})$, which is symmetric and non-degenerate
\cite[Section 3.4]{Kas91}.
By \cite[Lemma 1.4.4, Proposition 38.2.3]{Lusbook},
for $\mathbf{c}, \mathbf{c}' \in \Z_{\geq 0}^r$, we have
\[
\bigl(F(\mathbf{c}, \widetilde{w}),
F(\mathbf{c}', \widetilde{w})\bigr)_K = 0
\quad \text{if } \mathbf{c} \neq \mathbf{c}',
\]
and $\bigl(F(\mathbf{c}, \widetilde{w}),
F(\mathbf{c}, \widetilde{w})\bigr)_K$ is non-zero.
Hence, if we set
\[
F^{\mathrm{up}}(\mathbf{c}, \widetilde{w}) \coloneq
\frac{F(\mathbf{c}, \widetilde{w})}
{\bigl(F(\mathbf{c}, \widetilde{w}),
F(\mathbf{c}, \widetilde{w})\bigr)_K},
\]
then $\{F^{\mathrm{up}}(\mathbf{c}, \widetilde{w}) \mid
\mathbf{c} \in \Z_{\geq 0}^r\}$ is the basis of $U_q^-(w)$
dual to the PBW basis with respect to $(\ ,\ )_K$.
We call it the \emph{dual PBW basis}.

For $\mathbf{c}, \mathbf{c}' \in \Z_{\geq 0}^r$,
we write $\mathbf{c}' < \mathbf{c}$ if there exists $p$
with $1 \leq p \leq r$ such that
\[
c'_1 = c_1, \ldots, c'_{p - 1} = c_{p - 1}, \qquad c'_p < c_p.
\]
Let $\overline{\phantom{x}} : U_q^-(\mathfrak{g})
\to U_q^-(\mathfrak{g})$ be a $\Q$-algebra homomorphism
such that $f_i \mapsto f_i$, $q \mapsto q^{-1}$.
Define a $\Q$-linear map $\sigma : U_q^-(\mathfrak{g})
\to U_q^-(\mathfrak{g})$ by
\[
(\sigma(x), y)_K = \overline{(x, \overline{y})_K}
\qquad \text{for } x, y \in U_q^-(\mathfrak{g}).
\]

\begin{Prop}[{\cite[Theorem 4.29]{Kim12}}]
\label{Prop : dual canonical basis}
There exists a unique $\Q(q)$-basis
\[
\{B^{\mathrm{up}}(\mathbf{c}, \widetilde{w}) \mid
\mathbf{c} \in \Z_{\geq 0}^r\}
\]
of $U_q^-(w)$ satisfying
\[
\sigma\bigl(B^{\mathrm{up}}(\mathbf{c}, \widetilde{w})\bigr)
= B^{\mathrm{up}}(\mathbf{c}, \widetilde{w})
\]
and
\[
F^{\mathrm{up}}(\mathbf{c}, \widetilde{w})
= B^{\mathrm{up}}(\mathbf{c}, \widetilde{w})
+ \sum_{\mathbf{c}' < \mathbf{c}}\phi_{\mathbf{c}, \mathbf{c}'}
B^{\mathrm{up}}(\mathbf{c}', \widetilde{w}), \qquad
\phi_{\mathbf{c}, \mathbf{c}'} \in q\Z[q].
\]
This basis is called the \emph{dual canonical basis} of $U_q^-(w)$.
\end{Prop}

Let $U_q^-(\mathfrak{g})_{\Z[q^{\pm 1}]}$
be a $\Z[q^{\pm 1}]$-algebra generated by
$f_i^{(n)}$ $(i \in I, n \in \Z_{\geq 0})$, and
\[
U_q^-(\mathfrak{g})_{\Z[q^{\pm 1}]}^\vee \coloneq
\{x \in U_q^-(\mathfrak{g}) \mid
(x, U_q^-(\mathfrak{g})_{\Z[q^{\pm 1}]})_K
\subset \Z[q^{\pm 1}]\}.
\]
It is a $\Z[q^{\pm 1}]$-subalgebra of $U_q^-(\mathfrak{g})$
and it has a $\Z[q^{\pm 1}]$-basis called upper global basis
\cite{Kas91}.
Let
\[
U_q^-(w)_{\Z[q^{\pm 1}]}^\vee \coloneq
\bigoplus_{\mathbf{c} \in \Z_{\geq 0}^r}
\Z[q^{\pm 1}]B^{\mathrm{up}}(\mathbf{c}, \widetilde{w}).
\]
Then $U_q^-(w)_{\Z[q^{\pm 1}]}^\vee$
is independent of the choice of reduced word
and it is a subalgebra of $U_q^-(\mathfrak{g})_{\Z[q^{\pm 1}]}^\vee$.
Moreover, its basis
$\{B^{\mathrm{up}}(\mathbf{c}, \widetilde{w}) \mid
\mathbf{c} \in \Z_{\geq 0}^r\}$
is a subset of the upper global basis.
For detail, see \cite[$\S$4.7]{Kim12}.

For $1 \leq s \leq r$, we set
\[
\begin{aligned}
s^- &\coloneq
\max(\{t \in [1, s - 1] \mid i_t = i_s\} \cup \{0\}), \\
s^+ &\coloneq
\min(\{t \in [s + 1, r] \mid i_t = i_s\} \cup \{r + 1\}).
\end{aligned}
\]
We also set
\[
s^{\min} \coloneq \min\{t \in [1, r] \mid i_t = i_s\}, \qquad
s^{\max} \coloneq \max\{t \in [1, r] \mid i_t = i_s\}.
\]
For $1 \leq s \leq t \leq r$ satisfying $i_s = i_t$, we define
\[
D_{\widetilde{w}}[s, t] \coloneq
B^{\mathrm{up}}(\varepsilon_s + \varepsilon_{s^+} + \cdots
+ \varepsilon_{t^-} + \varepsilon_t, \widetilde{w}).
\]
When there is no risk of confusion, we simply write $D[s, t]$.

\subsection{Quiver Hecke algebras}
\label{ssec : quiver Hecke algebras}

Let $\mathfrak{g}$ be a symmetric Kac--Moody Lie algebra
with Cartan matrix $A = (a_{ij})_{i, j \in I}$,
and let $\kk$ be a field of characteristic zero.
For $i, j \in I$, choose a polynomial $Q_{ij}(u, v) \in \kk[u, v]$
of the form
\[
Q_{ij}(u, v) =
\begin{cases}
0 & \text{if } i = j, \\
\displaystyle\sum_{\substack{p, q \in \Z_{\geq 0} \\
p + q = -a_{ij}}} t_{i, j; p, q}u^p v^q & \text{if } i \neq j,
\end{cases}
\]
where $t_{i, j; p, q} \in \kk$ satisfy
\[
Q_{ij}(u, v) = Q_{ji}(v, u), \qquad
t_{i, j; -a_{ij}, 0} \in \kk^\times.
\]

Let $\mathfrak{S}_n = \langle s_1, \ldots, s_{n - 1}\rangle$
be the symmetric group, where $s_k$ is the transposition
interchanging $k$ and $k + 1$.
It acts naturally on $I^n$.
For $\beta \in Q_+$ with $|\beta| = n$, set
\[
I^\beta \coloneq
\{\nu = (\nu_1, \ldots, \nu_n) \in I^n \mid
\alpha_{\nu_1} + \cdots + \alpha_{\nu_n} = \beta\}.
\]

\begin{Def}
\label{Def : quiver Hecke algebra}
For $\beta \in Q_+$ with $|\beta| = n$,
the \emph{quiver Hecke algebra} $R(\beta)$ is the $\kk$-algebra
generated by
\[
\{e(\nu) \mid \nu \in I^\beta\}, \qquad
\{x_k \mid 1 \leq k \leq n\}, \qquad
\{\tau_m \mid 1 \leq m \leq n - 1\},
\]
subject to the following relations:
\[
\begin{gathered}
e(\nu)e(\nu') = \delta_{\nu, \nu'}e(\nu), \qquad
\sum_{\nu \in I^\beta}e(\nu) = 1, \\
x_k x_m = x_m x_k, \qquad x_k e(\nu) = e(\nu)x_k, \\
\tau_m e(\nu) = e(s_m(\nu))\tau_m, \qquad
\tau_k\tau_m = \tau_m\tau_k \quad \text{if } |k - m| > 1, \\
\tau_k^2 e(\nu) = Q_{\nu_k, \nu_{k + 1}}(x_k, x_{k + 1})e(\nu),
\end{gathered}
\]
\[
(\tau_k x_m - x_{s_k(m)}\tau_k)e(\nu) =
\begin{cases}
-e(\nu) & \text{if } m = k \text{ and } \nu_k = \nu_{k + 1}, \\
e(\nu) & \text{if } m = k + 1 \text{ and } \nu_k = \nu_{k + 1}, \\
0 & \text{otherwise},
\end{cases}
\]
and
\[
\begin{aligned}
&(\tau_{k + 1}\tau_k\tau_{k + 1}
- \tau_k\tau_{k + 1}\tau_k)e(\nu) \\
&\quad =
\begin{cases}
\displaystyle
\frac{Q_{\nu_k, \nu_{k + 1}}(x_k, x_{k + 1})
- Q_{\nu_k, \nu_{k + 1}}(x_{k + 2}, x_{k + 1})}
{x_k - x_{k + 2}}e(\nu)
& \text{if } \nu_k = \nu_{k + 2}, \\
0 & \text{otherwise}.
\end{cases}
\end{aligned}
\]
The algebra $R(\beta)$ is $\Z$-graded by
\[
\deg e(\nu) = 0, \qquad
\deg(x_k e(\nu)) = 2, \qquad
\deg(\tau_k e(\nu)) = -a_{\nu_k, \nu_{k + 1}}.
\]
\end{Def}

Let $R(\beta)\text{-}\mathrm{mod}$ be the category of finite-dimensional
$R(\beta)$-modules on which the elements $x_k$ act nilpotently,
and let $R(\beta)\text{-}\mathrm{gmod}$ be the category
of finite-dimensional graded $R(\beta)$-modules.
We set
\[
R\text{-}\mathrm{gmod} \coloneq
\bigoplus_{\beta \in Q_+}R(\beta)\text{-}\mathrm{gmod}.
\]
If $M$ is an $R(\beta)$-module, we set
$\operatorname{wt}(M) \coloneq -\beta \in Q_-$
and call it the \emph{weight} of $M$.
For $\beta, \gamma \in Q_+$, define
\[
e(\beta, \gamma) \coloneq
\sum_{\substack{\nu \in I^\beta \\ \nu' \in I^\gamma}}
e(\nu * \nu') \in R(\beta + \gamma),
\]
where $\nu * \nu'$ denotes the concatenation of $\nu$ and $\nu'$.
For $M \in R(\beta)\text{-}\mathrm{gmod}$
and $N \in R(\gamma)\text{-}\mathrm{gmod}$,
their \emph{convolution product} is
\[
M \circ N \coloneq R(\beta + \gamma)e(\beta, \gamma)
\otimes_{R(\beta) \otimes R(\gamma)}(M \otimes N).
\]

Let $\psi$ be the anti-involution of $R(\beta)$
fixing $e(\nu)$, $x_k$, and $\tau_k$.
For $M \in R(\beta)\text{-}\mathrm{gmod}$, its \emph{graded dual} is
\[
M^* \coloneq \operatorname{Hom}_{\kk}(M, \kk), \qquad
(M^*)_d = \operatorname{Hom}_{\kk}(M_{-d}, \kk),
\]
with the $R(\beta)$-action given by $(xf)(m) = f(\psi(x)m)$.
For $s \in \Z$, the grading shift $q^s M$ is defined by
\[
(q^s M)_d = M_{d - s}.
\]
A simple module $M \in R\text{-}\mathrm{gmod}$
is called \emph{self-dual} if $M \simeq M^*$.

Let $K(R\text{-}\mathrm{gmod})$ be the Grothendieck group
of $R\text{-}\mathrm{gmod}$.
The grading shift and convolution product make it
a $\Z[q, q^{-1}]$-algebra via
\[
[qM] = q[M], \qquad [M][N] = [M \circ N].
\]
The quiver Hecke algebra $R$ is called \emph{symmetric}
if $Q_{ij}(u, v)$ is a polynomial in $u - v$ for every $i, j \in I$.
In what follows, we assume that $R$ is symmetric.

\begin{Thm}[{\cite{KL09, Ro08, VV11}}]
\label{Thm : KLR categorification}
There exists a $\Z[q, q^{-1}]$-algebra isomorphism
\[
\operatorname{ch} : K(R\text{-}\mathrm{gmod})
\xrightarrow{\sim} U_q^-(\mathfrak{g})_{\Z[q^{\pm 1}]}^\vee
\]
under which the isomorphism classes of self-dual simple modules
correspond bijectively to the upper global basis
of $U_q^-(\mathfrak{g})_{\Z[q^{\pm 1}]}^\vee$.
\end{Thm}

Let $w \in W$ and fix a reduced word
$\widetilde{w} = s_{i_1}\cdots s_{i_r}$.
For $1 \leq s \leq t \leq r$ satisfying $i_s = i_t$,
since $D[s, t]$ is a member of the upper global basis,
it follows from Theorem~\ref{Thm : KLR categorification}
that there exists a unique self-dual simple module
$M[s, t] \in R\text{-}\mathrm{gmod}$
satisfying $\operatorname{ch}(M[s, t]) = D[s, t]$.

\begin{Def}[{\cite[$\S$2.2]{strata}}]
\label{Def : Cw}
Let $\mathscr{C}_w$ be the smallest monoidal abelian full subcategory
of $R\text{-}\mathrm{gmod}$ satisfying the following properties:
\begin{enumerate}
    \item $\mathscr{C}_w$ is stable under the subquotients, extensions,
    and grading shifts,
    \item $\mathscr{C}_w$ contains $M[s, s]$ for all $1 \leq s \leq r$.
\end{enumerate}
\end{Def}

Then by \cite{GLSqua}, $M \in R\text{-}\mathrm{gmod}$
belongs to $\mathscr{C}_w$ if and only if $\operatorname{ch}(M)$
belongs to $U_q^-(w)_{\Z[q^{\pm 1}]}^\vee$.
Hence, the following holds:

\begin{Prop}
\label{Prop : Cw categorification}
The isomorphism $\operatorname{ch}$ induces
a $\Z[q, q^{-1}]$-algebra isomorphism
\[
K(\mathscr{C}_w) \xrightarrow{\sim} U_q^-(w)_{\Z[q^{\pm 1}]}^\vee
\]
under which the isomorphism classes of self-dual simple modules
correspond bijectively to the dual canonical basis
of $U_q^-(w)_{\Z[q^{\pm 1}]}^\vee$.
\end{Prop}

For each $1 \leq k \leq r$, let $S_k \coloneq M[k, k]$
and we call it a cuspidal module.

\begin{Prop}[{\cite[$\S$2]{TW16}, \cite[Proposition 2.15]{strata}}]
\label{Prop : KLR PBW parameterization}
\begin{enumerate}
\item [(i)]
For $(c_k)_{k \in [1, r]} \in \Z_{\geq 0}^r$,
the head of $S_r^{\circ c_r} \circ S_{r - 1}^{\circ c_{r - 1}}
\circ \cdots \circ S_1^{\circ c_1}$
is a simple module corresponding to
$B^{\mathrm{up}}(\mathbf{c}, \widetilde{w})$ up to grading shift.
Moreover, if we ignore grading, this correspondence forms
a bijection between $\Z_{\geq 0}^r$ and $\Irr\mathscr{C}_w$.
\item[(ii)] The ring
    $K(\mathscr{C}_w)|_{q=1}$ is the polynomial ring
    generated by $\{[S_s] \mid s \in [1, r]\}$. 
    \end{enumerate}
\end{Prop}

For $L \in \Irr\mathscr{C}_w$, the corresponding sequence $(c_k)_k$
is called \emph{PBW parameter} and we write it as $\mathrm{PBW}(L)$.

\section{Invariants of monoidal categories}
\label{Sec : invariants of monoidal categories}
\subsection{Settings in this paper}\label{ssec : settings in this paper}
From now on, we take one of the following choices
for a monoidal category $\mathscr{C}$, an integer $N$,
modules $S_k$ and $T_k^{t_0}$, a matrix $A$, and a ring $K(\mathscr{C})_{\Z}$:
\begin{enumerate}
    \item[(i)] Let $\mathfrak{g}$ be a finite-dimensional simple Lie algebra
    and $\mathscr{C}_{\mathfrak{g}}$ be the category
    of finite-dimensional representations of $U_q(L\mathfrak{g})$. 
    Let $A = (a_{\iota \jmath})_{\iota, \jmath \in I_{\mathrm{fin}}}$ 
    be the Cartan matrix of $\mathfrak{g}_{\mathrm{fin}}$. 
    Let $[a, b]$ be an interval in $\Z$ of length $N$
    and $\mathfrak{s} = ((\iota_k, p_k))_{k \in \Z}$
    be an admissible sequence.
    Let $\mathscr{C}$ be the monoidal category
    $\mathscr{C}_{\mathfrak{g}}^{[a, b], \mathfrak{s}}$.
    Replacing $\mathfrak{s}$ by the admissible sequence
    $((\iota_{k + a - 1}, p_{k + a - 1}))_{k \in \Z}$,
    we may assume without loss of generality that $[a, b] = [1, N]$.
    Let $S_k = S_k^{\mathfrak{s}}$ be the cuspidal module
    in $\mathscr{C}_{\mathfrak{g}}^{[1, N], \mathfrak{s}}$
    for $1 \leq k \leq N$.
    Let $T_k^{t_0}$ be the simple module whose PBW parameter is
    $\sum_{\substack{1\leq s \leq k \\ \iota_s = \iota_k}} e_s$.
    Let
    \[
    K(\mathscr{C})_{\Z} \coloneq K(\mathscr{C})
    \simeq \Z\bigl[[S_k] \mid 1 \leq k \leq N\bigr]
    \]
    be the Grothendieck ring
    of $\mathscr{C}_{\mathfrak{g}}^{[1, N], \mathfrak{s}}$.
    \item[(ii)] Let $\mathfrak{g}$ be a symmetric Kac--Moody Lie algebra
    and $R$ be a symmetric quiver Hecke algebra
    corresponding to $\mathfrak{g}$
    over a base field of characteristic $0$. 
    Let $A = (a_{\iota, \jmath})$ be the Cartan matrix of $\mathfrak{g}$. 
    Let $w \in W$ and
    $\widetilde{w} = s_{\iota_1}\cdots s_{\iota_N}$
    be a reduced word.
    Let $\mathscr{C}$ be the monoidal category $\mathscr{C}_w$,
    and let $S_k$ be the cuspidal module for $1 \leq k \leq N$.
    Let $T_k^{t_0} \coloneq M[k^{\min}, k]$.
    Let
    \[
    K(\mathscr{C})_{\Z} \coloneq K(\mathscr{C}_w)|_{q = 1}
    \simeq \Z\bigl[[S_k] \mid 1 \leq k \leq N\bigr]
    \]
    be the Grothendieck ring of $\mathscr{C}_w$
    specialized at $q = 1$.
\end{enumerate}

In what follows, for $L, M \in \mathscr{C}$,
even if $\mathscr{C} = \mathscr{C}_w$,
we write the tensor product or the convolution product
of $L$ and $M$ as $L \otimes M$.
Also, we ignore the grading shifts when $\mathscr{C} = \mathscr{C}_w$.

\subsection{Invariants}
\begin{Def}
\label{Def : commutativity reality and primeness}
\begin{enumerate}
    \item[(i)] Let $L, M \in \Irr\mathscr{C}$.
    We say that $L$ and $M$ \emph{commute}
    if $L \otimes M \simeq M \otimes L$.
    We say that $L$ and $M$ \emph{strongly commute}
    if $L \otimes M \in \Irr\mathscr{C}$.
    \item[(ii)] Let $L \in \Irr\mathscr{C}$.
    We call $L$ \emph{real} if $L \otimes L \in \Irr\mathscr{C}$.
    We call $L$ \emph{imaginary}
    if $L \otimes L \notin \Irr\mathscr{C}$.
    We call $L$ \emph{prime} if $L$ is non-trivial and,
    if $L \simeq L_1 \otimes L_2$, then one of them is trivial.
\end{enumerate}
\end{Def}

Since $K(\mathscr{C})_{\Z}$ is commutative,
strong commutativity implies commutativity.

Let $\Omega$ be the smallest set of non-zero objects in $\mathscr{C}$
which contains $\Irr\mathscr{C}$
and is closed under taking non-zero subquotients and tensor products.
In \cite[$\S$3]{KKOP20} and \cite[$\S$3]{KKKO18},
the invariant $\Lambda \colon \Omega \times \Omega \to \Z$
is defined using $R$-matrices.
Also, the symmetric invariant
$\mathfrak{d} \colon \Omega \times \Omega \to \Z$ is defined by
\[
\mathfrak{d}(M, N)
\coloneq \frac{1}{2}\bigl(\Lambda(M, N) + \Lambda(N, M)\bigr).
\]
Actually, it is known that $\mathfrak{d}(M, N) \in \Z_{\geq 0}$
for $M, N \in \Omega$
\cite[$\S$3]{KKOP20}, \cite[$\S$3]{KKKO18}.
When $\mathscr{C} = \mathscr{C}_{\mathfrak{g}}^{[1, N], \mathfrak{s}}$,
the invariant $\Lambda^{\infty} \colon \Omega \times \Omega \to \Z$
is defined in \cite[Definition 3.6]{KKOP20}.
When $\mathscr{C} = \mathscr{C}_w$,
let $\Lambda^{\infty} \colon \Omega \times \Omega \to \Z$
be the invariant such that
\[
\Lambda^{\infty}(M, N)
= -(\operatorname{wt} M, \operatorname{wt} N).
\]
Let
\[
\widetilde{\Lambda}(M, N)
\coloneq \frac{1}{2}
\bigl(\Lambda(M, N) - \Lambda^{\infty}(M, N)\bigr).
\]
It is known that $\widetilde{\Lambda}(M, N) \in \Z$
\cite[Definition 3.6]{KKOP20}, \cite[Lemma 3.1.2]{KKKO18}.

\begin{Rem}
\label{Rem : convention for tilde Lambda}
This $\widetilde{\Lambda}$ is compatible with the definition
of $\widetilde{\Lambda}$ in \cite[Definition 3.1.1]{KKKO18}
for the quiver Hecke case.
However, in \cite[Definition 3.6]{KKOP20}
for the quantum affine case,
the original $\widetilde{\Lambda}$ is defined by
$\frac{1}{2}(\Lambda + \Lambda^{\infty})$.
Since the $\Lambda^{\infty}$-invariant for quantum affine algebras
corresponds to $-(\operatorname{wt}, \operatorname{wt})$
for quiver Hecke algebras via quantum Schur--Weyl duality
\cite[Theorem 4.12]{PBW24},
we define $\widetilde{\Lambda}$ for quantum affine algebras by
$\frac{1}{2}(\Lambda - \Lambda^{\infty})$.
\end{Rem}

\begin{Prop}[{\cite[Corollary 3.17]{KKOP20},
\cite[Lemma 3.2.3]{KKKO18}}]
\label{Prop : d vanishing and commutativity}
Let $L, M \in \Irr\mathscr{C}$ and suppose one of them is real.
Then the following are equivalent:
\begin{enumerate}
    \item[(i)] $\mathfrak{d}(L, M) = 0$,
    \item[(ii)] $L$ and $M$ commute,
    \item[(iii)] $L$ and $M$ strongly commute.
\end{enumerate}
\end{Prop}

\begin{Lem}[{\cite[Lemma 3.10, Proposition 3.11]{KKOP20},
\cite[Lemma 3.1.5, Proposition 3.2.10]{KKKO18}}]
\label{Lem : Lambda tensor products and subquotients}
Let $L, M, N \in \Irr\mathscr{C}$.
\begin{enumerate}
    \item[(i)] We have
    \[
    \begin{aligned}
    \Lambda(L, M \otimes N) &= \Lambda(L, M) + \Lambda(L, N), \\
    \Lambda(M \otimes N, L) &= \Lambda(M, L) + \Lambda(N, L).
    \end{aligned}
    \]
    \item[(ii)] Let $S$ be a simple subquotient of $M \otimes N$.
    Then
    \[
    \begin{aligned}
    \Lambda(L, S) &\leq \Lambda(L, M) + \Lambda(L, N), \\
    \Lambda(S, L) &\leq \Lambda(M, L) + \Lambda(N, L).
    \end{aligned}
    \]
\end{enumerate}
\end{Lem}

By Proposition~\ref{Prop : d vanishing and commutativity}
and Lemma~\ref{Lem : Lambda tensor products and subquotients},
if $L$ is real, then $L^{\otimes n}$ is real for every $n \in \Z_{\geq 0}$.

\begin{Prop}[{\cite[Proposition 4.2]{KKOP20},
\cite[Proposition 3.2.10]{KKKO18}}]
\label{Prop : commutativity with simple subquotients}
Let $L, M, N \in \Irr\mathscr{C}$.
When $L$ is real, the following conditions are equivalent:
\begin{enumerate}
    \item[(a)] $L$ strongly commutes with $M$ and $N$.
    \item[(b)] Any simple subquotient $S$ of $M \otimes N$
    commutes with $L$ and satisfies
    \[
    \Lambda(S, L) = \Lambda(M, L) + \Lambda(N, L).
    \]
    \item[(c)] Any simple subquotient $S$ of $M \otimes N$
    commutes with $L$ and satisfies
    \[
    \Lambda(L, S) = \Lambda(L, M) + \Lambda(L, N).
    \]
\end{enumerate}
\end{Prop}

\begin{Lem}[{\cite[Lemma 3.10, Proposition 3.11]{KKOP20}}]
\label{Lem : Lambda infinity tensor products and subquotients}
Let $L, M, N \in \Irr\mathscr{C}$.
Then
\[
\begin{aligned}
\Lambda^{\infty}(L, M \otimes N)
&= \Lambda^{\infty}(L, M) + \Lambda^{\infty}(L, N), \\
\Lambda^{\infty}(M \otimes N, L)
&= \Lambda^{\infty}(M, L) + \Lambda^{\infty}(N, L).
\end{aligned}
\]
Let $S$ be a simple subquotient of $M \otimes N$.
Then
\[
\begin{aligned}
\Lambda^{\infty}(L, S)
&= \Lambda^{\infty}(L, M) + \Lambda^{\infty}(L, N), \\
\Lambda^{\infty}(S, L)
&= \Lambda^{\infty}(M, L) + \Lambda^{\infty}(N, L).
\end{aligned}
\]
\end{Lem}

\begin{Prop}[{\cite[Theorems 3.2 and 3.12]{headsimplicity}}]
\label{Prop : simplicity of heads and socles}
Let $L, M \in \Irr\mathscr{C}$
and assume that at least one of them is real.
Then the head and the socle of $L \otimes M$ are simple.
\end{Prop}

We write the head of $L \otimes M$ as $L \nabla M$
and the socle as $L \Delta M$. 

\begin{Lem}[{\cite[Corollary~3.7, Corollary~3.14]{headsimplicity}}]
\label{Lem : hd soc injectivity}
  Let $L \in \Irr \mathscr{C}$ be a real module. 
  Then the maps 
  \[
  L \nabla - , L \Delta -: \Irr \mathscr{C} \to  \Irr \mathscr{C}
  \]
  are both injective. 
\end{Lem}

\begin{Lem}[{\cite[Theorem 4.11]{KKOP20},
\cite[Theorem 4.1.1]{KKKO18}}]
\label{Lem : Lambda inequalities for heads and socles}
Let $M, N \in \Irr\mathscr{C}$ and assume that $M$ is real.
\begin{enumerate}
    \item[(i)] $[M \otimes N] = [M \nabla N] + \sum_k[S_k]$
    with simple modules $S_k$ such that
    \[
    \Lambda(M, S_k) < \Lambda(M, M \nabla N) = \Lambda(M, N);
    \]
    \item[(ii)] $[M \otimes N] = [M \Delta N] + \sum_k[S_k]$
    with simple modules $S_k$ such that
    \[
    \Lambda(S_k, M) < \Lambda(M \Delta N, M) = \Lambda(N, M);
    \]
    \item[(iii)] $[N \otimes M] = [N \nabla M] + \sum_k[S_k]$
    with simple modules $S_k$ such that
    \[
    \Lambda(S_k, M) < \Lambda(N \nabla M, M) = \Lambda(N, M);
    \]
    \item[(iv)] $[N \otimes M] = [N \Delta M] + \sum_k[S_k]$
    with simple modules $S_k$ such that
    \[
    \Lambda(M, S_k) < \Lambda(M, N \Delta M) = \Lambda(M, N).
    \]
\end{enumerate}
In particular, $M \nabla N$ as well as $M \Delta N$
appears only once in the Jordan--H\"older series
of $M \otimes N$ in $\mathscr{C}$.
\end{Lem}

\begin{Lem}[{\cite[Corollary 3.18]{LocalizationsI},
\cite[Proposition 2.17]{PBW24}}]
\label{Lem : decrease of d}
Let $L, M \in \Irr\mathscr{C}$ and assume that $L$ is real.
\begin{enumerate}
    \item[(i)] If $\mathfrak{d}(L, M) > 0$,
    then for any simple subquotient $S$ of $L \otimes M$,
    \[
    \mathfrak{d}(L, S) < \mathfrak{d}(L, M).
    \]
    \item[(ii)] There exists $n \in \Z_{\geq 0}$ such that
    \[
    \mathfrak{d}(L, S) = 0
    \]
    for any simple subquotient $S$ of $L^{\otimes n}\otimes M$. 
\end{enumerate}
\end{Lem}

\begin{Lem}[{\cite[Proposition 4.7]{KKOP20},
\cite[Proposition 3.2.17]{KKKO18}, \cite[Lemma 2.28]{PBW24}}]
\label{Lem : d one exact sequence}
Let $L, M \in \Irr\mathscr{C}$ and assume that one of them is real.
Assume that $\mathfrak{d}(L, M) = 1$.
Then
\[
0 \longrightarrow L \Delta M
\longrightarrow L \otimes M
\longrightarrow L \nabla M \longrightarrow 0
\]
is an exact sequence. 
In addition, if $L, M$ are both real, then $L \Delta M$ and $L \nabla M$ are real modules.
\end{Lem}

\begin{Lem}[{\cite[Lemma 2.23]{KKOP24}}]
\label{Lem : reality of head commuting with a factor}
Let $L, M \in \Irr\mathscr{C}$ be real modules.
Assume that $\mathfrak{d}(L, L \nabla M) = 0$.
Then $L \nabla M$ is real.
\end{Lem}

\begin{Prop}[{\cite[Theorem 4.13]{KKOP20},
\cite[Theorem 4.1]{KK19}}]
\label{Prop : divisibility by a real simple module}
Let $X, L \in \Irr\mathscr{C}$.
Assume that $L$ is real.
If $[X] = [L]\phi$ for some $\phi \in K(\mathscr{C})_{\Z}$,
then $X \simeq L \otimes Y$ for some $Y \in \Irr\mathscr{C}$
such that $\mathfrak{d}(L, Y) = 0$.
\end{Prop}

For $k \in [1, N]$, let $\beta_k \coloneq s_{\iota_1}\cdots s_{\iota_{k - 1}}\alpha_{\iota_k}$. 

\begin{Lem}[{\cite[Proposition 5.4, Proposition 5.11]{Laurentfamily}}]
\label{Lem : Lambda and PBW}
We set
\[
\lambda_{st} \coloneq (-1)^{\delta(s > t)}\delta(s \neq t)
(\beta_s, \beta_t)
\qquad \text{for }s, t \in [1, N].
\]
Here, for a proposition $P$,
\[
\delta(P) =
\begin{cases}
1 & \text{if }P\text{ is true}, \\
0 & \text{if }P\text{ is false}.
\end{cases}
\]
Let $L, M \in \Irr \mathscr{C}$. 
Assume that at least one of them is real and they commute. 
Then
\[
\Lambda(L, M) = \mathrm{PBW}(L)^{T}(\lambda_{kl})_{1 \leq k, l\leq N}\mathrm{PBW}(M). 
\]
\end{Lem}

\begin{Lem}[{\cite[Lemma 3.2]{simplylaced}}]
\label{Lem : Lambda infinity of cuspidal modules}
Let $k, l \in [1, N]$.
Then
\[
\Lambda^{\infty}(S_k, S_l) = -(\beta_k, \beta_l).
\]
\end{Lem}
For $\mathbf{a}, \mathbf{b} \in \Z^N$, we write $\mathbf{a}<_{\mathrm{lex}}\mathbf{b}$ by 
if $a_1=b_1, \ldots, a_{k-1}=b_{k-1}$ and $a_k<b_k$ for some $1 \leq k\leq N$. 
This is a total order. 

\begin{Lem}[{\cite[Theorem 6.12]{PBW24}, \cite[Proposition 4.31]{Kim12}}]\label{Lem : tensor PBW}
\begin{enumerate}
  \item Let $L \in \Irr \mathscr{C}$ and $\mathbf{a} = \mathrm{PBW}(M)$. 
  Then 
  \[
  [L] \in  \bigl[\bigotimes_{k=1}^N S_k^{\otimes a_k}\bigr] 
  + \sum_{\mathbf{b}<_{\mathrm{lex}}\mathbf{a}}\Z \bigl[\bigotimes_{k=1}^N S_k^{\otimes b_k}\bigr]
  \]
  \item Let $L_1, L_2 \in \Irr \mathscr{C}$ and $\mathbf{a}_1, \mathbf{a}_2$ be their PBW parameters. 
  Then 
  \[
  [L_1\otimes L_2] \in [L]+ \sum_{\mathrm{PBW}(M)<_{\mathrm{lex}}\mathbf{a}_1+\mathbf{a}_2}  \Z [M]
  \]
  in $K(\mathscr{C})_{\Z}$ where $\mathrm{PBW}(L) = \mathbf{a}_1 + \mathbf{a}_2$. 
\end{enumerate}
\end{Lem}

\begin{Prop}
\label{Prop : tensor factors in category}
Let $M,N$ be nonzero objects in $\mathscr{C}_{\mathfrak{g}}$ for a quantum affine case and 
in $R\text{-gmod}$ for a quiver Hecke case. 
If $M\otimes N\in\mathscr{C}$, then $M,N\in\mathscr{C}$. 
In particular, if $L \in \Irr\mathscr{C}$ is non-trivial and does not have non-trivial tensor factorization 
in $\mathscr{C}$, then $L$ is a prime module in the entire category 
$\mathscr{C}_{\mathfrak{g}}$ or $R\text{-gmod}$. 
\end{Prop}

\begin{proof}
First, we consider the case
$\mathscr{C}=\mathscr{C}_{\mathfrak{g}}^{[1, N], \mathfrak{s}}$.
Since an object belongs to $\mathscr{C}$ if and only if all its
composition factors belong to $\mathscr{C}$, we may assume that
$M,N\in\Irr \mathscr{C}_{\mathfrak{g}}$.
Let $m_1,m_2$ be the highest $\ell$-weight monomials of $M,N$. 
By Proposition~\ref{Prop : PBW parameterization}, the assumption
$M\otimes N\in\mathscr{C}$ implies that $m_1m_2$ is a monomial in
$\{Y_{\pi(\iota_k),p_k}\mid k\in[1,N]\}$.
Thus $m_1$ and $m_2$ are also monomials in
$\{Y_{\pi(\iota_k),p_k}\mid k\in[1,N]\}$, and hence $M,N\in\mathscr{C}$.

Next, we consider the case $\mathscr{C}=\mathscr{C}_w$.
For $X\in R(\beta)\text{-gmod}$, set
\[
W(X)\coloneq
\{\gamma\in Q_+\cap(\beta-Q_+)\mid e(\gamma,\beta-\gamma)X\neq0\}.
\]
By \cite[(2.5)]{strata}, we have $X\in\mathscr{C}_w$ if and only if
$W(X)\subseteq Q_+\cap wQ_-$.
Moreover, by \cite[Lemma~2.2]{strata},
\[
W(M\circ N)=W(M)+W(N).
\]
Since $0\in W(M),W(N)$, we obtain
\[
W(M),W(N)\subseteq W(M\circ N)\subseteq Q_+\cap wQ_-.
\]
Hence $M,N\in\mathscr{C}_w$.
\end{proof}

\section{Monoidal categorification of cluster algebras}
\label{Sec : monoidal categorification of cluster algebras}

\subsection{Cluster algebras}
\label{ssec : cluster algebras}

Let $J = J_{\mathrm{uf}} \sqcup J_{\mathrm{fr}}$ be a finite set.
An integer matrix
\[
\widetilde{B}
= \begin{pmatrix} B \\ P \end{pmatrix}
= (b_{ij})_{(i, j) \in J \times J_{\mathrm{uf}}}
\]
is called an \emph{exchange matrix}
if its principal part $B$ is skew-symmetric.
Let $\Lambda = (\lambda_{ij})_{i, j \in J}$
be a skew-symmetric integer matrix.
Let $D = \operatorname{diag}(d_i)_{i \in J_{\mathrm{uf}}}$
be a diagonal integer matrix with $d_i > 0$.
The pair $(\widetilde{B}, \Lambda)$
is called a \emph{compatible pair of type $D$} if
\[
\widetilde{B}^{\mathsf{T}}\Lambda = (D\ 0),
\qquad
(D\ 0) = (\delta_{ij}d_i)_{(i, j) \in J_{\mathrm{uf}} \times J}.
\]
For a compatible pair $(\widetilde{B}, \Lambda)$ of type $D$
and $k \in J_{\mathrm{uf}}$, we define
\[
\mu_k\widetilde{B}
= (b'_{ij})_{(i, j) \in J \times J_{\mathrm{uf}}},
\qquad
\mu_k\Lambda = (\lambda'_{ij})_{i, j \in J}
\]
by
\[
b'_{ij} =
\begin{cases}
    -b_{ij}, & \text{if $i = k$ or $j = k$}, \\
    b_{ij} + [b_{ik}]_+[b_{kj}]_+
    - [-b_{ik}]_+[-b_{kj}]_+, & \text{otherwise},
\end{cases}
\]
\[
\lambda'_{ij} =
\begin{cases}
    -\lambda_{ik}
    + \displaystyle\sum_{\ell \in J}[-b_{\ell k}]_+\lambda_{i\ell},
    & \text{if $j = k$, $i \neq k$}, \\
    -\lambda_{kj}
    + \displaystyle\sum_{\ell \in J}[-b_{\ell k}]_+\lambda_{\ell j},
    & \text{if $i = k$, $j \neq k$}, \\
    \lambda_{ij}, & \text{otherwise},
\end{cases}
\]
where $[a]_+ \coloneq \max\{a, 0\}$ for any $a \in \mathbb{R}$.
It is known that
\begin{itemize}
    \item
    $\mu_k(\widetilde{B}, \Lambda)
    \coloneq (\mu_k\widetilde{B}, \mu_k\Lambda)$
    is also a compatible pair of type $D$,
    i.e., the principal part of $\mu_k\widetilde{B}$ is skew-symmetric,
    $\mu_k\Lambda$ is skew-symmetric, and
    $(\mu_k\widetilde{B})^{\mathsf{T}}\mu_k\Lambda = (D\ 0)$,
    \item $\mu_k^2 = \id$, 
\end{itemize}
see \cite{BZ05}.

Let $\mathbb{F}$ be the rational function field
in $|J|$ variables over $\Q$.
A \emph{seed} is a pair $(\mathbf{x}, \widetilde{B})$ where
\begin{itemize}
    \item $\mathbf{x} = (x_i)_{i \in J}$ satisfies that
    $\{x_i \mid i \in J\}$ is a free generating set
    of $\mathbb{F}$ over $\Q$,
    \item $\widetilde{B}$ is an exchange matrix.
\end{itemize}
For a seed $(\mathbf{x}, \widetilde{B})$ and $k \in J_{\mathrm{uf}}$,
the \emph{mutation} of $(\mathbf{x}, \widetilde{B})$ at $k$
is the pair $(\mu_k\mathbf{x}, \mu_k\widetilde{B})$ given by
\[
(\mu_k\mathbf{x})_i =
\begin{cases}
    x_i, & \text{if $i \neq k$}, \\
    x_k^{-1}\left(
    \displaystyle\prod_{j \in J}x_j^{[b_{jk}]_+}
    + \displaystyle\prod_{j \in J}x_j^{[-b_{jk}]_+}
    \right), & \text{if $i = k$}.
\end{cases}
\]
Then
$\mu_k(\mathbf{x}, \widetilde{B})
\coloneq (\mu_k\mathbf{x}, \mu_k\widetilde{B})$
is also a seed and $\mu_k^2 = \id$.
Let $\mathbb{T}$ denote the $|J_{\mathrm{uf}}|$-regular tree.
We label the edges of $\mathbb{T}$ by $J_{\mathrm{uf}}$
such that the $|J_{\mathrm{uf}}|$ different edges adjacent to the same vertex
of $\mathbb{T}$ receive different labels.

\begin{Def}
\label{Def : cluster patterns}
A \emph{cluster pattern}
$\mathcal{S}_{\mathbf{x}}
= \{(\mathbf{x}_t, \widetilde{B}_t) \mid t \in \mathbb{T}\}$
is an assignment of a seed $(\mathbf{x}_t, \widetilde{B}_t)$
to every vertex $t \in \mathbb{T}$ such that
\[
(\mathbf{x}_{t'}, \widetilde{B}_{t'})
= \mu_k(\mathbf{x}_t, \widetilde{B}_t)
\]
whenever $\begin{xy}(0,1)*+{t}="A",(10,1)*+{t'}="B",\ar@{-}^k"A";"B" \end{xy}$ in $\mathbb{T}$.
A \emph{$\Lambda$-cluster pattern}
$\mathcal{S}_{\mathbf{x}}^{\Lambda}
= \{(\mathbf{x}_t, \widetilde{B}_t, \Lambda_t)
\mid t \in \mathbb{T}\}$ of type $D$
is an assignment of a triple $(\mathbf{x}_t, \widetilde{B}_t, \Lambda_t)$
to every vertex $t \in \mathbb{T}$ such that
$\{(\mathbf{x}_t, \widetilde{B}_t) \mid t \in \mathbb{T}\}$
is a cluster pattern,
$(\widetilde{B}_t, \Lambda_t)$ is a compatible pair of type $D$, and
\[
(\mathbf{x}_{t'}, \widetilde{B}_{t'}, \Lambda_{t'})
= (\mu_k\mathbf{x}_t, \mu_k\widetilde{B}_t, \mu_k\Lambda_t)
\]
whenever $\begin{xy}(0,1)*+{t}="A",(10,1)*+{t'}="B",\ar@{-}^k"A";"B" \end{xy}$ in $\mathbb{T}$.
We call $(\mathbf{x}_t, \widetilde{B}_t, \Lambda_t)$
a \emph{$\Lambda$-seed}.
\end{Def}

For a $\Lambda$-seed $(\mathbf{x}_t, \widetilde{B}_t, \Lambda_t)$,
we usually write
\[
\mathbf{x}_t = (x_{i;t})_{i \in J},
\qquad
\widetilde{B}_t = (b_{ij}^{t}),
\qquad
\Lambda_t = (\lambda_{ij}^{t}).
\]
We call each $x_{i;t}$ $(i \in J)$ a \emph{cluster variable}.
More precisely, we call $x_{i;t}$ $(i \in J_{\mathrm{uf}})$
an \emph{unfrozen cluster variable} and
$x_{i;t}$ $(i \in J_{\mathrm{fr}})$ a \emph{frozen cluster variable}.
A monomial in $\{x_{i;t} \mid i \in J\}$
is called a \emph{cluster monomial}.

\begin{Def}
\label{Def : cluster algebras}
Let
$\mathcal{S}_{\mathbf{x}}^{\Lambda}
= \{(\mathbf{x}_t, \widetilde{B}_t, \Lambda_t)
\mid t \in \mathbb{T}\}$
be a $\Lambda$-cluster pattern.
The \emph{(partially compactified) cluster algebra}
$\overline{\mathcal{A}}(\mathcal{S}_{\mathbf{x}})$
is the $\Z$-subalgebra of $\mathbb{F}$
generated by $\{x_{i;t} \mid i \in J,\ t \in \mathbb{T}\}$.
The \emph{upper cluster algebra}
$\mathcal{U}(\mathcal{S}_{\mathbf{x}})$
is the $\Z$-subalgebra of $\mathbb{F}$ such that
\[
\mathcal{U}(\mathcal{S}_{\mathbf{x}})
= \bigcap_{t \in \mathbb{T}}
\Z[x_{i;t}^{\pm 1} \mid i \in J].
\]
The \emph{(partially compactified) $\Lambda$-cluster algebra}
$\overline{\mathcal{A}}(\mathcal{S}_{\mathbf{x}}^{\Lambda})$
is a cluster algebra
$\overline{\mathcal{A}}(\mathcal{S}_{\mathbf{x}})$
endowed with the additional data
$\Lambda = \{\Lambda_t \mid t \in \mathbb{T}\}$.
The \emph{$\Lambda$-upper cluster algebra}
$\mathcal{U}(\mathcal{S}_{\mathbf{x}}^{\Lambda})$
is defined in the same manner.
\end{Def}

It is known that
$\overline{\mathcal{A}}(\mathcal{S}_{\mathbf{x}})
\subset \mathcal{U}(\mathcal{S}_{\mathbf{x}})$ \cite{FZ02}.
For an exchange matrix $\widetilde{B}$,
we can associate a quiver $Q_{\widetilde{B}}$
characterized by the following properties:
\begin{enumerate}
    \item[(1)] its vertices are labeled by $J$,
    and its frozen vertices are labeled by $J_{\mathrm{fr}}$,
    \item[(2)] it has no loops, $2$-cycles,
    and arrows between frozen vertices,
    \item[(3)] for any $i \in J$ and $j \in J_{\mathrm{uf}}$,
    \[
    b_{ij}
    = (\text{the number of arrows from $i$ to $j$})
    - (\text{the number of arrows from $j$ to $i$}).
    \]
\end{enumerate}
Conversely, every quiver satisfying conditions (1) and (2)
determines an integer matrix $\widetilde{B}$ by condition (3). 

\subsection{Pointed elements}
\label{ssec : pointed elements}

Let
$\mathcal{S}_{\mathbf{x}}^{\Lambda}
= \{(\mathbf{x}_t, \widetilde{B}_t, \Lambda_t)
\mid t \in \mathbb{T}\}$
be a $\Lambda$-cluster pattern.
For $t \in \mathbb{T}$ and $i \in J_{\mathrm{uf}}$, we set
\[
\widehat{y}_{i;t} \coloneq \mathbf{x}_t^{\widetilde{B}_t e_i}.
\]

\begin{Def}[{\cite[Definition 3.1.4]{Qin17},
\cite[Definition 3.2.2]{Qin24}}]
\label{Def : pointed and bipointed elements}
Let $u \in \mathcal{U}(\mathcal{S}_{\mathbf{x}}^{\Lambda})$
and $t \in \mathbb{T}$.
\begin{enumerate}
    \item[(i)] We say that $u$ is \emph{pointed at $t$}
    if there exist $\mathbf{g} \in \Z^J$ and a polynomial
    \[
    F(u_i \mid i \in J_{\mathrm{uf}})
    = \sum_{\mathbf{v} \in \Z_{\geq 0}^{J_{\mathrm{uf}}}}c_{\mathbf{v}}\mathbf{u}^{\mathbf{v}}
    \in \Z[u_i \mid i \in J_{\mathrm{uf}}],
    \qquad c_0 = 1
    \]
    such that
    \[
    u = \mathbf{x}_t^{\mathbf{g}}
    \cdot F(\widehat{y}_{i;t} \mid i \in J_{\mathrm{uf}}).
    \]
    We write $\deg^t u \coloneq \mathbf{g}$
    and call it the \emph{degree} of $u$.
    We write
    \[
    F_u^t(u_i \mid i \in J_{\mathrm{uf}})
    \coloneq F(u_i \mid i \in J_{\mathrm{uf}})
    \]
    and call it the \emph{$F$-polynomial} of $u$.
    \item[(ii)] We say that $u$ is \emph{bipointed at $t$}
    if $u$ is pointed and the $F$-polynomial is of the form
    \[
    F_u^t(u_i \mid i \in J_{\mathrm{uf}})
    = \sum_{0 \leq \mathbf{v} \leq \mathbf{v}_0}c_{\mathbf{v}}\mathbf{u}^{\mathbf{v}},
    \qquad c_0 = c_{\mathbf{v}_0} = 1
    \]
    for some $\mathbf{v}_0 \in \Z_{\geq 0}^{J_{\mathrm{uf}}}$.
    Here, $0 \leq \mathbf{v} \leq \mathbf{v}_0$ means that
    $\mathbf{v}, \mathbf{v}_0 - \mathbf{v} \in \Z_{\geq 0}^{J_{\mathrm{uf}}}$.
    We write
    \[
    \operatorname{codeg}^t u
    \coloneq \deg^t u + \widetilde{B}_t\mathbf{v}_0
    \]
    and call it the \emph{codegree} of $u$.
\end{enumerate}
\end{Def}

\begin{Def}[{\cite{FG09}, \cite{FZ07},
\cite[Definition 3.3]{Cao}}]
\label{Def : tropical transformations and tropical points}
\begin{enumerate}
    \item[(i)] Let $t \in \mathbb{T}$ and $k \in J_{\mathrm{uf}}$.
    We set
    \[
    \phi^Y_{\mu_k t, t} \colon \Z^J \longrightarrow \Z^J;
    \qquad \mathbf{g}^t \longmapsto \mathbf{g}^{\mu_k t},
    \]
    where
    \[
    g_i^{\mu_k t} =
    \begin{cases}
        -g_k^t, & \text{if $i = k$}, \\
        g_i^t + [-b_{ik}^t]_+g_k^t + b_{ik}^t[g_k^t]_+,
        & \text{if $i \neq k$}.
    \end{cases}
    \]
    We set
    \[
    S_Y(\Z^{\max}) \coloneq
    \left\{
    [\mathbf{g}] \coloneq \{\mathbf{g}^t \mid t \in \mathbb{T}\}
    \;\middle|\;
    \begin{gathered}
    \mathbf{g}^{\mu_k t} = \phi^Y_{\mu_k t, t}(\mathbf{g}^t) \\
    \text{for all $t \in \mathbb{T}$, $k \in J_{\mathrm{uf}}$}
    \end{gathered}
    \right\}.
    \]
    \item[(ii)] Let $t \in \mathbb{T}$ and $k \in J_{\mathrm{uf}}$.
    We set
    \[
    \phi^X_{\mu_k t, t} \colon \Z^J \longrightarrow \Z^J;
    \qquad \mathbf{a}^t \longmapsto \mathbf{a}^{\mu_k t},
    \]
    where
    \[
    a_i^{\mu_k t} =
    \begin{cases}
        -a_k^t + \max\left\{
        \displaystyle\sum_{j \in J}[b_{jk}^t]_+a_j^t,
        \displaystyle\sum_{j \in J}[-b_{jk}^t]_+a_j^t
        \right\}, & \text{if $i = k$}, \\
        a_i^t, & \text{if $i \neq k$}.
    \end{cases}
    \]
    We set
    \[
    S_X(\Z^{\max}) \coloneq
    \left\{
    [\mathbf{a}] \coloneq \{\mathbf{a}^t \mid t \in \mathbb{T}\}
    \;\middle|\;
    \begin{gathered}
    \mathbf{a}^{\mu_k t} = \phi^X_{\mu_k t, t}(\mathbf{a}^t) \\
    \text{for all $t \in \mathbb{T}$, $k \in J_{\mathrm{uf}}$}
    \end{gathered}
    \right\}.
    \]
\end{enumerate}
\end{Def}

\begin{Def}[{\cite[Definition 3.4.2]{Qin24}}]
\label{Def : compatibly pointed elements}
Let $u \in \mathcal{U}(\mathcal{S}_{\mathbf{x}}^{\Lambda})$.
We say that $u$ is \emph{compatibly pointed}
if $u$ is pointed at any seed and
\[
\{\deg^t u \mid t \in \mathbb{T}\} \in S_Y(\Z^{\max}).
\]
\end{Def}

It is known that each cluster monomial is compatibly pointed
(see \cite{DWZ10, GHKK18}).

\begin{Lem}[{\cite[Proposition 4.3]{Cao}}]
\label{Lem : tropical points and Lambda}
Let $\mathcal{S}_{\mathbf{x}}^{\Lambda}$
be a $\Lambda$-cluster pattern.
For $[\mathbf{g}] = \{\mathbf{g}^t \mid t \in \mathbb{T}\}
\in S_Y(\Z^{\max})$, we have
\[
\{\Lambda_t\mathbf{g}^t \mid t \in \mathbb{T}\}
\in S_X(\Z^{\max}).
\]
\end{Lem}

\begin{Def}[{\cite{Qin17}}]
\label{Def : injective reachable seeds}
Let
$\mathcal{S}_{\mathbf{x}}^{\Lambda}
= \{(\mathbf{x}_t, \widetilde{B}_t, \Lambda_t)
\mid t \in \mathbb{T}\}$
be a $\Lambda$-cluster pattern.
Let $\operatorname{pr}_{\mathrm{uf}} \colon
\Z^J \longrightarrow \Z^{J_{\mathrm{uf}}}$
be the natural projection.
Then $t \in \mathbb{T}$ is called \emph{injective-reachable}
if there exist $t' \in \mathbb{T}$
and a permutation $\sigma$ of $J_{\mathrm{uf}}$ such that
\[
\operatorname{pr}_{\mathrm{uf}}\deg^t x_{k;t'}
= -f_{\sigma(k)}
\qquad \text{for all $k \in J_{\mathrm{uf}}$},
\]
where $f_{\sigma(k)} \in \Z^{J_{\mathrm{uf}}}$
is the $\sigma(k)$-th unit vector.
\end{Def}

It is known that $t \in \mathbb{T}$ is injective-reachable
if and only if every $t' \in \mathbb{T}$ is injective-reachable
(see \cite{Qin17}).

\begin{Prop}[{\cite[Proposition 3.4.8]{Qin24}}]
\label{Prop : compatible pointedness and bipointedness}
Let $\mathcal{S}_{\mathbf{x}}^{\Lambda}$
be a $\Lambda$-cluster pattern.
Assume that $t \in \mathbb{T}$ is injective-reachable.
If $u \in \mathcal{U}(\mathcal{S}_{\mathbf{x}}^{\Lambda})$
is compatibly pointed, then $u$ is bipointed at any seed.
Moreover, if
$u, u' \in \mathcal{U}(\mathcal{S}_{\mathbf{x}}^{\Lambda})$
are compatibly pointed with the same degree, then
\[
\operatorname{codeg}^{t'}u = \operatorname{codeg}^{t'}u'
\qquad \text{for all $t' \in \mathbb{T}$}.
\]
\end{Prop}

\subsection{Definitions of monoidal categorification of cluster algebras}
\label{ssec : definitions of monoidal categorification}

Let $\mathscr{C}$ be a monoidal category
in Subsection~\ref{ssec : settings in this paper}.
Let $J_{\mathrm{uf}} \subset J$ be finite sets.

\begin{Def}[{\cite[Definition 6.1]{KKOP20}}]
\label{Def : monoidal seed}
A \emph{monoidal seed} in $\mathscr{C}$ is a pair
$\mathscr{S} = (\{M_i\}_{i \in J}, \widetilde{B})$
consisting of a strongly commuting family $\{M_i\}_{i \in J}$
of real modules in $\mathscr{C}$ and an exchange matrix
$\widetilde{B} = (b_{ij})_{(i, j) \in J \times J_{\mathrm{uf}}}$.
\end{Def}

For a monoidal seed $\mathscr{S}$, we set
$\Lambda^{\mathscr{S}} \coloneq (\Lambda(M_i, M_j))_{i, j \in J}$.
By Proposition~\ref{Prop : d vanishing and commutativity},
$\Lambda^{\mathscr{S}}$ is a skew-symmetric integer matrix.

\begin{Def}[{\cite[Definitions 6.3, 6.5]{KKOP20}}]
\label{Def : Lambda admissible monoidal seed}
A monoidal seed
$\mathscr{S} = (\{M_i\}_{i \in J}, \widetilde{B})$
is \emph{$\Lambda$-admissible}
if it satisfies the following conditions:
\begin{enumerate}
    \item[(i)] For any $k \in J_{\mathrm{uf}}$,
    there exists a simple module $M'_k$ in $\mathscr{C}$
    such that there is an exact sequence
    \[
    0 \longrightarrow
    \bigotimes_{i : b_{ik} > 0}M_i^{\otimes b_{ik}}
    \longrightarrow M_k \otimes M'_k
    \longrightarrow
    \bigotimes_{i : b_{ik} < 0}M_i^{\otimes(-b_{ik})}
    \longrightarrow 0,
    \]
    and $M'_k$ commutes with $M_i$ for any $i \neq k$.
    \item[(ii)] $(\widetilde{B}, \Lambda^{\mathscr{S}})$
    is a compatible pair of type $D = 2I_{|J_{\mathrm{uf}}|}$.
\end{enumerate}
\end{Def}

If $\mathscr{S}$ is a $\Lambda$-admissible monoidal seed,
$M'_k$ is uniquely determined by Lemma~\ref{Lem : hd soc injectivity}
For a $\Lambda$-admissible monoidal seed $\mathscr{S}$
and $k \in J_{\mathrm{uf}}$, we set
\[
\mu_k\mathscr{S}
\coloneq
\bigl(\{M_i\}_{i \neq k} \cup \{M'_k\}, \mu_k\widetilde{B}\bigr).
\]

\begin{Prop}[{\cite[Proposition 6.4]{KKOP20}}]
\label{Prop : mutation of Lambda admissible monoidal seeds}
Let $\mathscr{S} = (\{M_i\}_{i \in J}, \widetilde{B})$
be a $\Lambda$-admissible monoidal seed in $\mathscr{C}$
and $k \in J_{\mathrm{uf}}$.
\begin{enumerate}
    \item[(i)] $\mu_k\mathscr{S}$ is a monoidal seed i.e. $M'_k$ is real,
    \item[(ii)] $\Lambda^{\mu_k\mathscr{S}}
    = \mu_k\Lambda^{\mathscr{S}}$.
\end{enumerate}
\end{Prop}

By Proposition~\ref{Prop : mutation of Lambda admissible monoidal seeds}
and the property that compatibility does not change by mutation,
$(\mu_k\widetilde{B}, \Lambda^{\mu_k\mathscr{S}})$
is a compatible pair of type $2I_{|J_{\mathrm{uf}}|}$.

\begin{Def}[{\cite[Definition 7.6]{KKOP24}}]
\label{Def : completely Lambda admissible monoidal seed}
A monoidal seed
$\mathscr{S} = (\{M_i\}_{i \in J}, \widetilde{B})$
is \emph{completely $\Lambda$-admissible}
if $\mathscr{S}$ is $\Lambda$-admissible and,
for any finite sequence $(k_1, \ldots, k_\ell)$ in $J_{\mathrm{uf}}$,
$\mu_{k_\ell}\cdots\mu_{k_1}\mathscr{S}$
is well-defined and a $\Lambda$-admissible monoidal seed.
\end{Def}

\begin{Def}[{\cite[Definition 6.7]{KKOP20}}]
\label{Def : Lambda monoidal categorification}
Let $\overline{\mathcal{A}}(\mathcal{S}_{\mathbf{x}}^{\Lambda})$
be a $\Lambda$-cluster algebra with initial $\Lambda$-seed
$(\mathbf{x}_{t_0}, \widetilde{B}_{t_0}, \Lambda_{t_0})$.
We call $\mathscr{C}$ a \emph{$\Lambda$-monoidal categorification}
of $\overline{\mathcal{A}}$ if it satisfies the following conditions:
\begin{enumerate}
    \item[(1)] There exists a $\Z$-algebra isomorphism
    \[
    \iota \colon K(\mathscr{C})_{\Z}
    \xrightarrow{\sim}
    \overline{\mathcal{A}}(\mathcal{S}_{\mathbf{x}}^{\Lambda}).
    \]
    \item[(2)] There exists a completely $\Lambda$-admissible
    monoidal seed
    $\mathscr{S}_0 \coloneq (\{M_i\}_{i \in J}, \widetilde{B})$
    such that
    \[
    \bigl(\{\iota([M_i])\}_{i \in J},
    \widetilde{B}, \Lambda^{\mathscr{S}_0}\bigr)
    = (\mathbf{x}_{t_0}, \widetilde{B}_{t_0}, \Lambda_{t_0}).
    \]
\end{enumerate}
\end{Def}

By Definition~\ref{Def : Lambda admissible monoidal seed}
and Proposition~\ref{Prop : mutation of Lambda admissible monoidal seeds},
for any finite sequence $(k_1, \ldots, k_m)$ of $J_{\mathrm{uf}}$ and
$\mu_{k_m}\cdots\mu_{k_1}\mathscr{S}_0
= (\{M'_i\}_i, \widetilde{B}')$, we have
\[
\bigl(\{\iota([M'_i])\}_i, \widetilde{B}',
\Lambda^{\mu_{k_m}\cdots\mu_{k_1}\mathscr{S}_0}\bigr)
= (\mathbf{x}_t, \widetilde{B}_t, \Lambda_t),
\]
where $t = \mu_{k_m}\cdots\mu_{k_1}t_0$. 

For each $t \in \mathbb{T}$, the real module corresponding to a cluster variable $x_{i;t}$ is denoted by 
$M_{i;t}$, and the set  $\{M_{i;t} \mid i \in J\}$ is called a \emph{monoidal cluster}. 

\begin{Prop}
\label{Prop : cluster monomials and real prime modules}
Let $\mathscr{C}$ be a $\Lambda$-monoidal categorification of
$\overline{\mathcal{A}}(\mathcal{S}_{\mathbf{x}}^{\Lambda})$
and let
\[
\iota \colon K(\mathscr{C})_{\Z}
\xrightarrow{\sim}
\overline{\mathcal{A}}(\mathcal{S}_{\mathbf{x}}^{\Lambda})
\]
be the corresponding isomorphism.
\begin{enumerate}
    \item[(i)] For any cluster monomial
    $\mathbf{x}_t^{\mathbf{a}} \in \overline{\mathcal{A}}$,
    there exists a real module $M \in \mathscr{C}$
    such that $\iota^{-1}(\mathbf{x}_t^{\mathbf{a}}) = [M]$.
    \item[(ii)] For any cluster variable
    $x_{i;t} \in \overline{\mathcal{A}}$,
    there exists a real prime module $M \in \mathscr{C}$
    such that $\iota^{-1}(x_{i;t}) = [M]$.
\end{enumerate}
\end{Prop}

\begin{proof}
(i) For $t = \mu_{k_m}\cdots\mu_{k_1}t_0$,
$\iota([M_{i;t}]) = x_{i;t}$ for any $i \in J$, where
\[
\mu_{k_m}\cdots\mu_{k_1}\mathscr{S}_0
= (\{M_{i;t}\}_i, \widetilde{B}_t).
\]
Since $\{M_{i;t}\}$ is a real commuting family,
$\bigotimes_{i \in J}M_{i;t}^{\otimes a_i}$
is a real simple module and
\[
\iota^{-1}(\mathbf{x}_t^{\mathbf{a}})
= \left[\bigotimes_{i \in J}M_{i;t}^{\otimes a_i}\right].
\]

(ii) It is enough to show that each $M_{i;t}$ is prime.
By Proposition~\ref{Prop : tensor factors in category},
it is enough to see that $[M_{i;t}]$ in $K(\mathscr{C})_{\Z}$
is a prime element.
By \cite{GLSclu}, each $x_{i;t}$ is an irreducible element
in $\overline{\mathcal{A}}(\mathcal{S}_{\mathbf{x}}^{\Lambda})$.
Since
$K(\mathscr{C})_{\Z}
\simeq \Z\bigl[[S_i] \mid 1 \leq i \leq N\bigr]$,
$K(\mathscr{C})_{\Z}$ is a unique factorization domain.
Now we have that $[M_{i;t}]$ in $K(\mathscr{C})_{\Z}$
is a prime element.
\end{proof}

A simple module corresponding to a cluster monomial is called \emph{reachable}. 

\subsection{$\Lambda$-cluster algebra structure on $K(\mathscr{C})_{\Z}$}
\label{ssec : Lambda cluster algebra structure}

Let $\mathscr{C}$, $N$, $S_k$, $T_k^{t_0}$ $(1 \leq k \leq N)$
be as in Subsection~\ref{ssec : settings in this paper}.
Let $J \coloneq \{1, \ldots, N\}$,
$J_{\mathrm{fr}} \coloneq \{k \in J \mid k^+ > N\}$,
and $J_{\mathrm{uf}} \coloneq J \setminus J_{\mathrm{fr}}$.
Let $Q_{t_0}$ be the quiver with vertex set $J$,
frozen vertex set $J_{\mathrm{fr}}$, and the following arrows:
\begin{itemize}
    \item ordinary arrows:
    $s \xrightarrow{|a_{\iota_s\iota_t}|} t$
    if $1 \leq s < t < s^+ < t^+$ and $s^+ \leq N$,
    \item horizontal arrows:
    $s \longrightarrow s^-$
    if $1 \leq s^- < s \leq N$.
\end{itemize}
We denote by $\widetilde{B}_{t_0}$
the exchange matrix corresponding to $Q_{t_0}$.

\begin{Thm}[{\cite[Theorem 8.1]{KKOP24},
\cite[Theorem 11.2.3]{KKKO18}}]
\label{Thm : Lambda monoidal categorification}
The pair
\[
\mathscr{S}_{t_0}
\coloneq (\{T_k^{t_0}\}_{k \in J}, \widetilde{B}_{t_0})
\]
is a completely $\Lambda$-admissible monoidal seed in $\mathscr{C}$.
For a $\Lambda$-cluster pattern $\mathcal{S}_{\mathbf{x}}^{\Lambda}$
with initial $\Lambda$-seed
$(\mathbf{x}_{t_0}, \widetilde{B}_{t_0},
\Lambda^{\mathscr{S}_{t_0}})$,
$\mathscr{C}$ provides a $\Lambda$-monoidal categorification by
\[
\iota \colon K(\mathscr{C})_{\Z}
\xrightarrow{\sim}
\overline{\mathcal{A}}(\mathcal{S}_{\mathbf{x}}^{\Lambda}),
\qquad
[T_k^{t_0}] \longmapsto x_{k;t_0}.
\]
\end{Thm}

\begin{Rem}
\label{Rem : initial monoidal clusters and i boxes}
In \cite[Theorem 7.20]{KKOP24}, the initial monoidal cluster for quantum affine case is
\[
\left\{
L\left(\prod_{\substack{k \leq s \leq N \\ \iota_s = \iota_k}}
Y_{\pi(\iota_s), q^{p_s}}\right)
\;\middle|\; k \in J
\right\},
\]
not $\{T_k^{t_0}\}_{k \in J}$.
However, we can obtain $\mathscr{S}_{t_0}$
from the monoidal seed in \cite[Theorem 7.20]{KKOP24}
by performing the mutation sequence
proposed in \cite[$\S$13]{GLSKac} and relabeling.
The family $\{T_k^{t_0}\}_{k \in J}$ corresponds to
a family of $\mathbf{i}$-boxes
$\{[k^{\min}, k] \mid k \in J\}$,
and we can calculate the exchange matrix
corresponding to the $\mathbf{i}$-boxes by \cite{KK26}. 
By the same argument, we can check that $t_0 \in \mathbb{T}$ is injective reachable. 
\end{Rem}

\subsection{Laurent families}
\label{ssec : Laurent families}

Let $\mathscr{C}$ be as in
Subsection~\ref{ssec : settings in this paper}, and let
\[
\iota \colon K(\mathscr{C})_{\Z}
\xrightarrow{\sim}
\overline{\mathcal{A}}(\mathcal{S}_{\mathbf{x}}^{\Lambda})
\]
be a $\Lambda$-monoidal categorification as in
Subsection~\ref{ssec : Lambda cluster algebra structure}.

For a real commuting family $\{L_i \mid 1 \leq i \leq r\}$
and $\mathbf{a} \in \Z_{\geq 0}^r$, let
\[
L^{\mathbf{a}} \coloneq \bigotimes_{i=1}^r L_i^{\otimes a_i}.
\]

\begin{Def}[{\cite[Definition 3.2]{Laurentfamily}}]
\label{Def : Laurent family}
A finite set of real modules $\{L_i \mid 1 \leq i \leq r\}$
is called a \emph{Laurent family} in $\mathscr{C}$
if it satisfies the following conditions:
\begin{enumerate}
    \item[(i)] $\mathfrak{d}(L_i, L_j) = 0$
    for $1 \leq i, j \leq r$,
    \item[(ii)] it is independent, i.e.,
    if $L^{\mathbf{a}} \simeq L^{\mathbf{b}}$
    for $\mathbf{a}, \mathbf{b} \in \Z_{\geq 0}^r$,
    then $\mathbf{a} = \mathbf{b}$,
    \item[(iii)] if $X \in \Irr\mathscr{C}$ commutes with all $L_i$,
    then there exists $\mathbf{a} \in \Z_{\geq 0}^r$
    such that $X \simeq L^{\mathbf{a}}$.
\end{enumerate}
\end{Def}

\begin{Prop}[{\cite[Proposition 3.6]{Laurentfamily}}]
\label{Prop : monoidal clusters are Laurent families}
Any monoidal cluster in $\mathscr{C}$ is a Laurent family.
\end{Prop}

\begin{Lem}[{\cite[Lemma 2.14]{Laurentfamily}}]
\label{Lem : algebraic independence of Laurent families}
Let $\mathcal{M} = \{L_i \mid 1 \leq i \leq r\}$
be a Laurent family.
Then $\{[L_i] \mid 1 \leq i \leq r\} \subset K(\mathscr{C})_{\Z}$
is algebraically independent over $\Z$.
\end{Lem}

\begin{Prop}[{\cite[Corollary 3.4]{Laurentfamily}}]
\label{Prop : Laurent expansions with respect to a Laurent family}
Let $\mathcal{M} = \{L_i \mid 1 \leq i \leq r\}$
be a Laurent family.
Then
\[
K(\mathscr{C})_{\Z}
\subset \Z\bigl[[L_i]^{\pm 1} \mid 1 \leq i \leq r\bigr]
\]
as subrings of $\operatorname{Frac}K(\mathscr{C})_{\Z}$.
\end{Prop}

\begin{Prop}
\label{Prop : cardinality of Laurent families}
Any Laurent family has exactly $N$ elements.
\end{Prop}

\begin{proof}
Let $\mathcal{M} = \{L_i \mid 1 \leq i \leq r\}$
be a Laurent family.
Then
\[
\Z\bigl[[L_i] \mid 1 \leq i \leq r\bigr]
\subset K(\mathscr{C})_{\Z}
\subset \Z\bigl[[L_i]^{\pm 1} \mid 1 \leq i \leq r\bigr]
\]
holds by
Proposition~\ref{Prop : Laurent expansions with respect to a Laurent family}.
By Lemma~\ref{Lem : algebraic independence of Laurent families},
the transcendence degrees of
\[
\operatorname{Frac}\Z\bigl[[L_i] \mid 1 \leq i \leq r\bigr],
\qquad
\operatorname{Frac}\Z\bigl[[L_i]^{\pm 1} \mid 1 \leq i \leq r\bigr]
\]
over $\Q$ are both $r$.
However, since $K(\mathscr{C})_{\Z}$ is a polynomial ring
in $N$ variables, the transcendence degree of
$\operatorname{Frac}K(\mathscr{C})_{\Z}$ over $\Q$ is $N$.
Then we have $r = N$.
\end{proof}

Let $\mathcal{M} = \{L_i \mid 1 \leq i \leq N\}$
be a Laurent family.
For $X \in \Irr\mathscr{C}$,
there exists $\mathbf{a} \in \Z_{\geq 0}^N$ such that
\[
[L^{\mathbf{a}} \otimes X]
= \sum_{s=1}^{\ell}c_s[L^{\mathbf{b}_s}]
\qquad
\text{for }c_s \in \Z_{>0},\quad
\mathbf{b}_s \in \Z_{\geq 0}^N
\]
by
Proposition~\ref{Prop : Laurent expansions with respect to a Laurent family}.
Since $L^{\mathbf{a}}$ is real,
$L^{\mathbf{a}} \nabla X$ and $L^{\mathbf{a}} \Delta X$ are simple,
and their multiplicities in $L^{\mathbf{a}} \otimes X$ are both $1$.
We take $\mathbf{h}, \mathbf{s} \in \Z_{\geq 0}^N$ such that
\[
L^{\mathbf{a}} \nabla X \simeq L^{\mathbf{h}},
\qquad
L^{\mathbf{a}} \Delta X \simeq L^{\mathbf{s}}.
\]

\begin{Lem}[{\cite[Lemma 3.1]{Laurentfamily}}]
\label{Lem : well defined degrees for Laurent families}
In this setting, $\mathbf{h} - \mathbf{a}$ and
$\mathbf{s} - \mathbf{a}$ in $\Z^N$ are well-defined.
\end{Lem}

\begin{Def}[{\cite[Definition 3.8]{KK19},
\cite[Definition 3.7]{Laurentfamily}}]
\label{Def : degrees with respect to a Laurent family}
Let $\mathcal{M} = \{L_i \mid 1 \leq i \leq N\}$
be a Laurent family and $X \in \Irr\mathscr{C}$.
We take $\mathbf{h}, \mathbf{s}, \mathbf{a} \in \Z_{\geq 0}^N$
as above.
We write
\[
\deg^{\mathcal{M}}X \coloneq \mathbf{h} - \mathbf{a},
\qquad
\operatorname{codeg}^{\mathcal{M}}X
\coloneq \mathbf{s} - \mathbf{a}.
\]
\end{Def}

For a Laurent family
$\mathcal{M} = \{L_i \mid 1 \leq i \leq N\}$, we set
\[
\begin{aligned}
\Lambda_{\mathcal{M}}
&\coloneq (\Lambda(L_i, L_j))_{1 \leq i, j \leq N}, \\
\Lambda^{\infty}_{\mathcal{M}}
&\coloneq (\Lambda^{\infty}(L_i, L_j))_{1 \leq i, j \leq N}, \\
\widetilde{\Lambda}_{\mathcal{M}}
&\coloneq (\widetilde{\Lambda}(L_i, L_j))_{1 \leq i, j \leq N}.
\end{aligned}
\]

\begin{Prop}[{\cite[Lemma 3.6]{KK19},
\cite[Lemma 5.12]{Laurentfamily}}]
\label{Prop : Lambda inequalities for Laurent expansions}
Let $\mathcal{M} = \{L_i \mid 1 \leq i \leq N\}$
be a Laurent family and $X \in \Irr\mathscr{C}$.
Let
\[
[L^{\mathbf{a}} \otimes X]
= \sum_{s=1}^{\ell}c_s[L^{\mathbf{b}_s}]
\qquad\text{for }\mathbf{b}_s \in \Z_{\geq 0}^N,
\]
and let
\[
L^{\mathbf{a}} \nabla X \simeq L^{\mathbf{h}},
\qquad
L^{\mathbf{a}} \Delta X \simeq L^{\mathbf{s}}.
\]
\begin{enumerate}
    \item[(i)] For any $s$, we have
    \[
    \Lambda^{\infty}_{\mathcal{M}}(\mathbf{h} - \mathbf{b}_s) = 0.
    \]
    \item[(ii)] 
    For $1\leq s \leq \ell$ such that $\mathbf{b}_s \neq \mathbf{h}$, 
    \[
    \Lambda_{\mathcal{M}}(\mathbf{h} - \mathbf{b}_s)\in \Z_{\geq 0}^N\setminus \{0\}. 
    \]
    For $1\leq s \leq \ell$ such that $\mathbf{b}_s \neq \mathbf{s}$, 
    \[
    \Lambda_{\mathcal{M}}(\mathbf{b}_s - \mathbf{s})\in \Z_{\geq 0}^N\setminus \{0\}. 
    \]
\end{enumerate}
\end{Prop}

\begin{Lem}[{\cite[Lemma 3.9]{Laurentfamily}}]
\label{Lem : injectivity of degrees for Laurent families}
Let $\mathcal{M} = \{L_i \mid 1 \leq i \leq N\}$
be a Laurent family.
Then
\[
\deg^{\mathcal{M}},\ \operatorname{codeg}^{\mathcal{M}}
\colon \Irr\mathscr{C} \longrightarrow \Z^N
\]
are injective.
\end{Lem}

\begin{Prop}[{\cite{KK19}, \cite[Proposition 5.24]{CFM}}]
\label{Prop : support of Laurent expansions at monoidal clusters}
Let $\mathcal{M}_t = \{M_{i;t} \mid 1 \leq i \leq N\}$
be a monoidal cluster corresponding to $t \in \mathbb{T}$
and $X \in \Irr\mathscr{C}$.
Let
\[
[M_t^{\mathbf{a}} \otimes X]
= \sum_{s=1}^{\ell}c_s[M_t^{\mathbf{b}_s}]
\qquad\text{for }\mathbf{a}, \mathbf{b}_s \in \Z_{\geq 0}^N,
\]
and let
\[
M_t^{\mathbf{a}} \nabla X \simeq M_t^{\mathbf{h}},
\qquad
M_t^{\mathbf{a}} \Delta X \simeq M_t^{\mathbf{s}}.
\]
Then
\[
\mathbf{b}_s - \mathbf{h},\quad \mathbf{s} - \mathbf{b}_s
\in \widetilde{B}_t\Z_{\geq 0}^{J_{\mathrm{uf}}}.
\]
\end{Prop}

By Proposition~\ref{Prop : support of Laurent expansions at monoidal clusters},
Definition~\ref{Def : pointed and bipointed elements}, and
Lemma~\ref{Lem : Lambda inequalities for heads and socles},
for $t \in \mathbb{T}$,
$\iota([X]) \in \mathcal{U}(\mathcal{S}_{\mathbf{x}}^{\Lambda})$
is bipointed at $t$, and
\[
\deg^{\mathcal{M}_t}X = \deg^t\iota([X]),
\qquad
\operatorname{codeg}^{\mathcal{M}_t}X
= \operatorname{codeg}^t\iota([X]).
\]

\begin{Prop}[{\cite[Theorem 3.16]{KK19}, \cite{Qin17}}]
\label{Prop : compatible pointedness of simple modules}
For $X \in \Irr\mathscr{C}$,
$\iota([X]) \in \mathcal{U}(\mathcal{S}_{\mathbf{x}}^{\Lambda})$
is compatibly pointed.
\end{Prop}

\begin{Prop}
\label{Prop : denominators and commutativity for Laurent families}
Let $\mathcal{M} = \{L_i \mid 1 \leq i \leq N\}$
be a Laurent family and $X \in \Irr\mathscr{C}$.
Let
\[
[X] = \sum_{s=1}^{\ell}c_s\prod_{i=1}^N[L_i]^{\mathbf{b}(s)_i}
\]
be a Laurent expansion, where $c_s \in \Z_{>0}$
and $\mathbf{b}(s) \in \Z^N$.
For $1 \leq k \leq N$, let
\[
\mathbf{d}(X)_k
\coloneq -\min\{\mathbf{b}(s)_k \mid 1 \leq s \leq \ell\}.
\]
Assume that each $[L_i]$ is a prime element in $K(\mathscr{C})_{\Z}$.
Then $\mathbf{d}(X)_k \leq 0$ if and only if $X$ commutes with $L_k$.
\end{Prop}

\begin{proof}
If $\mathbf{d}(X)_k \leq 0$, then we can prove that
$X$ commutes with $L_k$ by the same proof as
\cite[Lemma 4.2]{KK19}.

Assume that $X$ commutes with $L_k$.
By Lemma~\ref{Lem : decrease of d}, there exists $n \in \Z_{\geq 0}$
such that any simple subquotient of
\[
\left(\bigotimes_{i\ne k}L_i^{\otimes n}\right)\otimes X
\]
commutes with $L_i$ for all $i \ne k$.
Since
\[
\mathfrak{d}\left(
\left(\bigotimes_{i\ne k}L_i^{\otimes n}\right)\otimes X, L_k
\right)
= \mathfrak{d}(X, L_k) = 0,
\]
any simple subquotient of
$\left(\bigotimes_{i\ne k}L_i^{\otimes n}\right)\otimes X$
also commutes with $L_k$
by Lemma~\ref{Lem : Lambda tensor products and subquotients}.
By the definition of a Laurent family, any simple subquotient of
$\left(\bigotimes_{i\ne k}L_i^{\otimes n}\right)\otimes X$
is a monomial of $\mathcal{M}$.
It means $\mathbf{d}(X)_k \leq 0$.
\end{proof}

\begin{Prop}\label{Prop : frozen commutativity}
  Let $t \in \mathbb{T}$ and $j \in J_{\mathrm{fr}}$. 
  Then the real prime module $M_{j ; t}$ corresponding to a frozen cluster variable commutes with 
  all simple modules in $\mathscr{C}$. 
\end{Prop}
\begin{proof}
  By \cite{GLSclu}, each $[M_{i ; t}]$ for $i \in J$ is a prime element in $K(\mathscr{C})_{\Z}$. 
  The claim follows by Proposition~\ref{Prop : denominators and commutativity for Laurent families} 
  and the fact that 
  \[
  K(\mathscr{C})_{\Z} \subset 
  \Z[[M_{k ; t}]^{\pm 1}, [M_{i ; t}] \mid k \in J_{\mathrm{uf}}, i \in J_{\mathrm{fr}}]. 
  \]
\end{proof}

\begin{Prop}
\label{Prop : degrees of Laurent families form an integral basis}
Let $\mathcal{M} = \{L_i \mid 1 \leq i \leq N\}$
be a Laurent family and $t \in \mathbb{T}$.
Then
\[
\{\deg^t\iota([L_i]) \mid 1 \leq i \leq N\} \subset \Z^N
\]
is a $\Z$-basis of $\Z^N$.
\end{Prop}

\begin{proof}
We set $\mathbf{g}_i \coloneq \deg^t\iota([L_i])$ and
$\mathcal{M}_t \coloneq \{M_{i;t} \mid 1 \leq i \leq N\}$.
First we prove that $\{\mathbf{g}_i\}$ is linearly independent.
Assume
\[
\sum_{i=1}^N a_i\mathbf{g}_i
= \sum_{i=1}^N b_i\mathbf{g}_i
\qquad\text{for }a_i, b_i \in \Z_{\geq 0},\quad a_i b_i = 0.
\]
Then $\deg^{\mathcal{M}_t}L^{\mathbf{a}}
= \deg^{\mathcal{M}_t}L^{\mathbf{b}}$ holds.
Since $\deg^{\mathcal{M}_t}$ is injective,
$L^{\mathbf{a}} \simeq L^{\mathbf{b}}$.
Since $\mathcal{M}$ is independent, $\mathbf{a} = \mathbf{b}$ holds.

Next we prove that $\{\mathbf{g}_i\}$ generates $\Z^N$.
Let
\[
[L^{\mathbf{a}}\otimes M_{k;t}]
= \sum_{s=1}^{\ell}c_s[L^{\mathbf{b}_s}]
\qquad\text{for }\mathbf{a}, \mathbf{b}_s \in \Z_{\geq 0}^N,
\quad c_s \in \Z_{>0}.
\]
Since each $\iota([L_i])$ for $1 \leq i \leq N$
and $\iota([M_{i;t}])$ is pointed at $t$,
so are $\iota([L^{\mathbf{a}}\otimes M_{k;t}])$
and $\iota([L^{\mathbf{b}_s}])$ for $1 \leq s \leq \ell$.
By the positivity of $c_s$, there exists $1 \leq s_0 \leq \ell$
such that
\[
\deg^t\iota([L^{\mathbf{a}}\otimes M_{k;t}])
= \deg^t\iota([L^{\mathbf{b}_{s_0}}]).
\]
Thus we have
\[
\sum_{i=1}^N a_i\mathbf{g}_i + f_k
= \sum_{i=1}^N(\mathbf{b}_{s_0})_i\mathbf{g}_i,
\]
where $f_k$ is the $k$-th unit vector of $\Z^N$.
This means that $\{\mathbf{g}_i \mid 1 \leq i \leq N\}$
generates $\Z^N$.
\end{proof}

\section{Monoidal analog of cluster-tilting theory}
\label{Sec : monoidal analog of cluster tilting theory}

Let $\mathscr{C}$ and $N$ be as in
Subsection~\ref{ssec : settings in this paper}.

\subsection{Maximal real sets}
\label{ssec : maximal real set}

\begin{Def}
\label{Def : maximal real set}
Let $X$ be a set.
A set $\{L_\lambda\}_{\lambda \in X} \subset \Irr\mathscr{C}$
is called a \emph{maximal real set}
if it satisfies the following conditions:
\begin{enumerate}
    \item[(i)] $L_\lambda$ is a real prime module
    for all $\lambda \in X$,
    \item[(ii)] $\mathfrak{d}(L_{\lambda_1}, L_{\lambda_2}) = 0$
    for all $\lambda_1, \lambda_2 \in X$,
    \item[(iii)] if a real prime module $L \in \Irr\mathscr{C}$
    satisfies $\mathfrak{d}(L_\lambda, L) = 0$
    for all $\lambda \in X$,
    then $L \in \{L_\lambda\}_{\lambda \in X}$.
\end{enumerate}
\end{Def}

This is a monoidal analog of a maximal rigid object. 

Note that for any commuting family $\mathcal{S}$ of real prime modules,
there is a maximal real set which contains $\mathcal{S}$
by Zorn's lemma.

\begin{Lem}
\label{Lem : finiteness of maximal real sets}
A maximal real set in $\mathscr{C}$ is a finite set.
\end{Lem}

\begin{proof}
Let $\{L_\lambda\}_{\lambda \in X}$ be a maximal real set.
We choose $\lambda_1, \ldots, \lambda_k \in X$ such that
$\{\mathrm{PBW}(L_{\lambda_i}) \mid 1 \leq i \leq k\}$
forms a basis of
$\sum_{\lambda \in X}\Q\mathrm{PBW}(L_\lambda) \subset \Q^N$.
For $\lambda \in X$, we can write
\[
m\mathrm{PBW}(L_\lambda)
+ \sum_{i=1}^k m_i\mathrm{PBW}(L_{\lambda_i})
= \sum_{i=1}^k p_i\mathrm{PBW}(L_{\lambda_i}),
\]
where $m \in \Z_{>0}$, $m_i, p_i \in \Z_{\geq 0}$,
and $m_i p_i = 0$ for all $i$.
Since $\{L_\lambda\} \cup \{L_{\lambda_i} \mid 1 \leq i \leq k\}$
is a commuting set,
\[
L_\lambda^{\otimes m}
\otimes \bigotimes_{i=1}^k L_{\lambda_i}^{\otimes m_i},
\qquad
\bigotimes_{i=1}^k L_{\lambda_i}^{\otimes p_i}
\]
are both simple.
Then these are simple modules with the same PBW parameter,
thus we have
\begin{equation}
\label{eq : PBW relation in a maximal real set}
L_\lambda^{\otimes m}
\otimes \bigotimes_{i=1}^k L_{\lambda_i}^{\otimes m_i}
\simeq \bigotimes_{i=1}^k L_{\lambda_i}^{\otimes p_i}.
\end{equation}
Since $K(\mathscr{C})_{\Z}$ is a unique factorization domain,
$[\bigotimes_{i=1}^k L_{\lambda_i}] \in K(\mathscr{C})_{\Z}$
has a factorization into prime elements as a ring-theoretic meaning.
We write $\{q_1, \ldots, q_\ell\}$ for the set of pairwise
non-associate prime elements appearing in this factorization.
For each $\lambda \in X$, there exists a unique
$\mathbf{a}_\lambda \in \Z_{\geq 0}^{\ell}$ such that
$[L_\lambda]$ is associate to
$\prod_{j=1}^{\ell}q_j^{a_{\lambda, j}}$.
Suppose that there exist $\lambda \neq \lambda' \in X$ such that
$\mathbf{a}_\lambda - \mathbf{a}_{\lambda'}
\in \Z_{\geq 0}^{\ell}$.
Then there is some $\phi \in K(\mathscr{C})_{\Z}$ such that
$[L_\lambda] = [L_{\lambda'}]\phi$.
Since $L_{\lambda'}$ is real,
$L_\lambda \simeq L_{\lambda'} \otimes Y$
for some $Y \in \Irr\mathscr{C}$
by Proposition~\ref{Prop : divisibility by a real simple module}.
It contradicts the primeness of $L_\lambda$.
Thus $\mathbf{a}_\lambda - \mathbf{a}_{\lambda'} \notin \Z_{\geq 0}^{\ell}$ 
holds for any $\lambda \neq \lambda'$.
To finish the proof, it is enough to prove that
for any infinite sequence
$(\mathbf{a}_n)_{n \in \Z_{\geq 0}} \subset \Z_{\geq 0}^m$,
there exist $n_1 \neq n_2 \in \Z_{\geq 0}$ such that
$\mathbf{a}_{n_2} - \mathbf{a}_{n_1} \in \Z_{\geq 0}^m$.
For any infinite sequence in $\Z_{\geq 0}$,
it is easy to see that there exists a non-decreasing infinite subsequence.
In $(\mathbf{a}_n)_{n \in \Z_{\geq 0}}$,
we take an infinite subsequence whose first component is non-decreasing,
and then take an infinite subsequence whose next component
is non-decreasing.
Repeating this up to the $m$-th component, we obtain a subsequence
$(\mathbf{a}_{n_i})_{i \in \Z_{\geq 0}}$
which is non-decreasing in every component.
In particular,
$\mathbf{a}_{n_2} - \mathbf{a}_{n_1} \in \Z_{\geq 0}^m$.
\end{proof}

We can regard a maximal real set as a monoidal analog
of a cluster-tilting object by the following proposition.

\begin{Prop}
\label{Prop : simple modules commuting with a maximal real set}
Let $\{L_i \mid 1 \leq i \leq r\}$ be a maximal real set.
Let $X \in \Irr\mathscr{C}$.
If $X$ commutes with all $L_i$,
then there exists $\mathbf{a} \in \Z_{\geq 0}^r$
such that $X \simeq L^{\mathbf{a}}$.
\end{Prop}

\begin{proof}
By the definition of a maximal real set,
it is enough to show that if $X$ commutes with all $L_i$,
then $X$ is real.
Let $\{T_k \mid 1 \leq k \leq N\}$ be one of the monoidal clusters.
By Lemma~\ref{Lem : decrease of d}
and the reality of the monomials in $\{T_k\}$,
there exists $\mathbf{b} \in \Z_{\geq 0}^N$ such that
$X \nabla T^{\mathbf{b}}$ commutes with all $T_k$.
Since $\{T_k \mid 1 \leq k \leq N\}$ is a Laurent family
by Proposition~\ref{Prop : monoidal clusters are Laurent families},
there exists $\mathbf{c} \in \Z_{\geq 0}^N$ such that
$X \nabla T^{\mathbf{b}} \simeq T^{\mathbf{c}}$.
By the same argument, there exists $\mathbf{d} \in \Z_{\geq 0}^r$
such that $L^{\mathbf{d}} \nabla T^{\mathbf{b}}$ and
$L^{\mathbf{d}} \nabla T^{\mathbf{c}}$ commute with all $L_i$.
Since $L^{\mathbf{d}}$, $T^{\mathbf{b}}$, and $T^{\mathbf{c}}$
are real and
\[
\mathfrak{d}(L^{\mathbf{d}},
L^{\mathbf{d}} \nabla T^{\mathbf{b}})
= \mathfrak{d}(L^{\mathbf{d}},
L^{\mathbf{d}} \nabla T^{\mathbf{c}}) = 0,
\]
we have that $L^{\mathbf{d}} \nabla T^{\mathbf{b}}$ and
$L^{\mathbf{d}} \nabla T^{\mathbf{c}}$ are real
by Lemma~\ref{Lem : reality of head commuting with a factor}.
Then there exist $\mathbf{e}, \mathbf{f} \in \Z_{\geq 0}^r$
such that
\[
L^{\mathbf{d}} \nabla T^{\mathbf{b}} \simeq L^{\mathbf{e}},
\qquad
L^{\mathbf{d}} \nabla T^{\mathbf{c}} \simeq L^{\mathbf{f}},
\]
because $\{L_i \mid 1 \leq i \leq r\}$ is a maximal real set.
Now we consider
$L^{\mathbf{d}} \otimes X \otimes T^{\mathbf{b}} \in \mathscr{C}$.
Since $L^{\mathbf{d}} \otimes X \in \Irr\mathscr{C}$
and $T^{\mathbf{b}}$ is real,
it has a simple head by
Proposition~\ref{Prop : simplicity of heads and socles}.
There are two surjections:
\[
\begin{aligned}
L^{\mathbf{d}} \otimes X \otimes T^{\mathbf{b}}
&\twoheadrightarrow
L^{\mathbf{d}} \otimes (X \nabla T^{\mathbf{b}})
\simeq L^{\mathbf{d}} \otimes T^{\mathbf{c}}
\twoheadrightarrow L^{\mathbf{f}}, \\
L^{\mathbf{d}} \otimes X \otimes T^{\mathbf{b}}
&\simeq X \otimes L^{\mathbf{d}} \otimes T^{\mathbf{b}}
\twoheadrightarrow X \otimes L^{\mathbf{e}}.
\end{aligned}
\]
Since $L^{\mathbf{f}}$ and $X \otimes L^{\mathbf{e}}$ are simple,
both are isomorphic to the head of
$L^{\mathbf{d}} \otimes X \otimes T^{\mathbf{b}}$,
thus we have $L^{\mathbf{f}} \simeq X \otimes L^{\mathbf{e}}$.
By the reality of $L^{\mathbf{f}}$, it follows that $X$ is real.
\end{proof}

\begin{Lem}
\label{Lem : localization at a maximal real set}
Let $\{L_i \mid 1 \leq i \leq r\}$ be a maximal real set, and let
\[
G \coloneq
\langle [L_i]^{\pm 1} \mid 1 \leq i \leq r \rangle
\subset \bigl(\operatorname{Frac} K(\mathscr{C})_{\Z}\bigr)^\times.
\]
Then
\[
\Z[G] = K(\mathscr{C})_{\Z}
\bigl[[L_i]^{-1} \mid 1 \leq i \leq r\bigr]
\]
as subalgebras of $\operatorname{Frac} K(\mathscr{C})_{\Z}$.
Moreover, $G$ is a free $\Z$-module.
\end{Lem}

\begin{proof}
The inclusion
$\Z[G] \subset K(\mathscr{C})_{\Z}
\bigl[[L_i]^{-1} \mid 1 \leq i \leq r\bigr]$
is clear.
For $X \in \Irr\mathscr{C}$,
there exists $\mathbf{a} \in \Z_{\geq 0}^r$
such that any simple subquotient of $L^{\mathbf{a}} \otimes X$
commutes with all $L_i$ by Lemma~\ref{Lem : decrease of d}.
By Proposition~\ref{Prop : simple modules commuting with a maximal real set},
any simple subquotient of $L^{\mathbf{a}} \otimes X$
is a monomial in $\{L_i \mid 1 \leq i \leq r\}$.
It follows that $[X]$ is a Laurent polynomial in
$\{[L_i] \mid 1 \leq i \leq r\}$, and
$\Z[G] \supset K(\mathscr{C})_{\Z}
\bigl[[L_i]^{-1} \mid 1 \leq i \leq r\bigr]$.

Next, we prove that $G$ is free.
Since $G$ is finitely generated,
it is enough to prove that $G$ is torsion-free.
Suppose $([L^{\mathbf{a}}]/[L^{\mathbf{b}}])^k = 1$ for $\mathbf{a}, \mathbf{b} \in \Z_{\geq 0}^r$ and 
$k \in \Z_{\geq 1}$. 
Then we have $[L^{k\mathbf{a}}] = [L^{k\mathbf{b}}]$ in $K(\mathscr{C})_{\Z}$. 
Then it holds that $k\mathrm{PBW}(L^{\mathbf{a}}) = k\mathrm{PBW}(L^{\mathbf{b}})$, 
hence $L^{\mathbf{a}}$ and $L^{\mathbf{b}}$ have the same PBW parameter and they are isomorphic. 
Therefore $[L^{\mathbf{a}}]/[L^{\mathbf{b}}] = 1$ holds and 
$G$ is torsion free. 
\end{proof}
Note that $\Z[G]$ as a subalgebra of $\operatorname{Frac} K(\mathscr{C})_{\Z}$ is 
isomorphic to $\Z[G]$ as a group ring, 
that because simple classes are linearly independent in $K(\mathscr{C})_{\Z}$.

\begin{Prop}
\label{Prop : real prime modules give prime elements}
Let $L \in \Irr\mathscr{C}$ be a real prime module.
Then $[L] \in K(\mathscr{C})_{\Z}$ is a prime element.
\end{Prop}

\begin{proof}
Since $K(\mathscr{C})_{\Z}$ is a unique factorization domain,
it is enough to prove that $[L] \in K(\mathscr{C})_{\Z}$
is irreducible.
We take $x, y \in K(\mathscr{C})_{\Z}$ such that $[L] = xy$.
Let $\{L = L_1, \allowbreak L_2, \ldots, L_r\}$
be a maximal real set,
and let $G$ be as in
Lemma~\ref{Lem : localization at a maximal real set}.
Since $[L]$ is invertible in $\Z[G]$, so are $x$ and $y$.
By Lemma~\ref{Lem : localization at a maximal real set}, we have
\[
\Z[G]^\times
= \left\{
\pm \frac{[L^{\mathbf{b}}]}{[L^{\mathbf{a}}]}
\;\middle|\; \mathbf{a}, \mathbf{b} \in \Z_{\geq 0}^r
\right\}.
\]
Therefore,
\[
x = \varepsilon_1
\frac{[L^{\mathbf{b}_1}]}{[L^{\mathbf{a}_1}]},
\qquad
y = \varepsilon_2
\frac{[L^{\mathbf{b}_2}]}{[L^{\mathbf{a}_2}]}
\]
for $\varepsilon_1, \varepsilon_2 \in \{\pm 1\}$
and $\mathbf{a}_1, \mathbf{b}_1, \mathbf{a}_2, \mathbf{b}_2
\in \Z_{\geq 0}^r$.
Then we have
$[L^{\mathbf{b}_1}] = [L^{\mathbf{a}_1}]\varepsilon_1 x$
in $K(\mathscr{C})_{\Z}$.
By Proposition~\ref{Prop : divisibility by a real simple module},
there exists $X \in \Irr\mathscr{C}$ such that
$\varepsilon_1 x = [X]$.
By the same argument, there exists $Y \in \Irr\mathscr{C}$
such that $\varepsilon_2 y = [Y]$.
Then we have $[L] = \varepsilon_1\varepsilon_2[X \otimes Y]$
in $K(\mathscr{C})_{\Z}$.
By the primeness of $L$, at least one of $X$ and $Y$ is trivial,
and we have that at least one of $x$ and $y$ is $\pm 1$.
\end{proof}

By Proposition~\ref{Prop : real prime modules give prime elements},
the following corollary is immediate.

\begin{Cor}
\label{Cor : unique prime tensor factorization in C}
For a real module $L \in \mathscr{C}$,
tensor factorization of $L$ into prime modules
is unique up to permutation.
\end{Cor}

\begin{Cor}
\label{Cor : unique prime tensor factorization in Cg}
For a real module $L \in \mathscr{C}_{\mathfrak{g}}$,
tensor factorization of $L$ into prime modules
is unique up to permutation.
\end{Cor}

\begin{proof}
By Proposition~\ref{Prop : Cg factorization},
it is enough to prove that any real module
in $\mathscr{C}_{\mathfrak{g}}^0$
has a unique tensor factorization into prime modules.
By Proposition~\ref{Prop : whole category}
and Proposition~\ref{Prop : PBW parameterization},
any module $L \in \mathscr{C}_{\mathfrak{g}}^0$
is in $\mathscr{C}_{\mathfrak{g}}^{[c, d], \mathfrak{s}}$
for some finite interval $[c, d]$
and an admissible sequence $\mathfrak{s}$.
Now the claim is clear
by Corollary~\ref{Cor : unique prime tensor factorization in C}.
\end{proof}

\begin{Lem}
\label{Lem : real prime modules give prime elements over C}
Let $L \in \Irr\mathscr{C}$ be a real prime module.
Then $[L] \in \C \otimes_{\Z} K(\mathscr{C})_{\Z}$
is a prime element.
\end{Lem}

\begin{proof}
Let $\{L = L_1, L_2, \ldots, L_r\}$
be a maximal real set in $\mathscr{C}$.
Since $\C \otimes_{\Z} K(\mathscr{C})_{\Z}$
is a polynomial ring over $\C$,
it is enough to prove that $[L]$ is an irreducible element.
We take $x, y \in \C \otimes_{\Z} K(\mathscr{C})_{\Z}$
such that $[L] = xy$.
Since $[L]$ is invertible in $\C \otimes_{\Z} \Z[G]$,
so are $x$ and $y$.
Since
\[
\C[G]^\times
= \left\{
\alpha\frac{[L^{\mathbf{b}}]}{[L^{\mathbf{a}}]}
\;\middle|\; \alpha \in \C^\times,\ 
\mathbf{a}, \mathbf{b} \in \Z_{\geq 0}^r
\right\},
\]
we have
\[
x = \alpha_1\frac{[L^{\mathbf{b}_1}]}{[L^{\mathbf{a}_1}]},
\qquad
y = \alpha_2\frac{[L^{\mathbf{b}_2}]}{[L^{\mathbf{a}_2}]}
\]
for $\alpha_1, \alpha_2 \in \C^\times$
and $\mathbf{a}_1, \mathbf{b}_1, \mathbf{a}_2, \mathbf{b}_2
\in \Z_{\geq 0}^r$.
Let $x = \sum_{i=1}^k \beta_i[M_i]$,
where $\beta_i \in \C^\times$, $M_i \in \Irr\mathscr{C}$, and
\[
\mathrm{PBW}(M_1) >_{\mathrm{lex}} \mathrm{PBW}(M_2)
>_{\mathrm{lex}} \cdots >_{\mathrm{lex}} \mathrm{PBW}(M_k).
\]
Then
\[
[L^{\mathbf{b}_1}]
= \sum_{i=1}^k \alpha_1^{-1}\beta_i
[L^{\mathbf{a}_1} \otimes M_i].
\]
By Lemma~\ref{Lem : tensor PBW},
$\alpha_1^{-1}\beta_1 = 1$ and
\[
\mathrm{PBW}(L^{\mathbf{b}_1})
= \mathrm{PBW}(L^{\mathbf{a}_1}) + \mathrm{PBW}(M_1).
\]
We prove $\alpha_1^{-1}\beta_s \in \Z$ by induction.
Assume $\alpha_1^{-1}\beta_1, \ldots,
\alpha_1^{-1}\beta_{s-1} \in \Z$.
Then
\[
[L^{\mathbf{b}_1}]
- \sum_{i=1}^{s-1}\alpha_1^{-1}\beta_i
[L^{\mathbf{a}_1} \otimes M_i]
= \sum_{i=s}^k \alpha_1^{-1}\beta_i
[L^{\mathbf{a}_1} \otimes M_i],
\]
and the left-hand side is in $K(\mathscr{C})_{\Z}$.
By Lemma~\ref{Lem : tensor PBW},
the coefficient of the simple module whose PBW parameter is
$\mathrm{PBW}(L^{\mathbf{a}_1}) + \mathrm{PBW}(M_s)$
in the right-hand side is $\alpha_1^{-1}\beta_s$.
Then we have $\alpha_1^{-1}\beta_s \in \Z$.
As a result, we have
\[
\alpha_1^{-1}x
= \frac{[L^{\mathbf{b}_1}]}{[L^{\mathbf{a}_1}]}
\in K(\mathscr{C})_{\Z}.
\]
By the same argument, we have
\[
\alpha_2^{-1}y
= \frac{[L^{\mathbf{b}_2}]}{[L^{\mathbf{a}_2}]}
\in K(\mathscr{C})_{\Z}.
\]
Then $\alpha_1^{-1}x = [X]$ and $\alpha_2^{-1}y = [Y]$
for some $X, Y \in \Irr\mathscr{C}$
by Proposition~\ref{Prop : divisibility by a real simple module}.
By the primeness of $L$, we conclude that
at least one of $X$ and $Y$ is trivial.
\end{proof}

\begin{Thm}
\label{Thm : Laurent families and maximal real sets}
Let $\{L_i \mid 1 \leq i \leq r\}$
be a real commuting family in $\mathscr{C}$.
Then it is a Laurent family if and only if it is a maximal real set.
\end{Thm}

\begin{proof}
Assume that $\{L_i\}$ is a Laurent family.
First, we prove that any $L_i$ is a prime module.
Let
\[
L_i \simeq \bigotimes_{k=1}^{\ell}M_k^{\otimes c_k}
\]
be a prime factorization.
Since $L_i$ commutes with any $L_j$ for $1 \leq j \leq r$,
so does each $M_k$.
Then, for each $1 \leq k \leq \ell$,
$M_k \simeq L^{\mathbf{a}_k}$
for some $\mathbf{a}_k \in \Z_{\geq 0}^r \setminus \{0\}$
by Definition~\ref{Def : Laurent family}(iii).
Then $L_i \simeq L^{\sum_{k=1}^{\ell}c_k\mathbf{a}_k}$.
By Definition~\ref{Def : Laurent family}(ii),
we have $\sum_{k=1}^{\ell}c_k\mathbf{a}_k = e_i$,
and we have that $L_i$ is prime.

Next, we prove the maximality.
Let $L \in \Irr\mathscr{C}$ be a real prime module
which commutes with any $L_i$.
Then $L \simeq L^{\mathbf{a}}$
for some $\mathbf{a} \in \Z_{\geq 0}^r$
by Definition~\ref{Def : Laurent family}(iii).
By the primeness of $L$, $L \simeq L_i$
for some $1 \leq i \leq r$.
Then we have that $\{L_i\}$ is a maximal real set.

If $\{L_i\}$ is a maximal real set, then it is a Laurent family
by Corollary~\ref{Cor : unique prime tensor factorization in C}
and Proposition~\ref{Prop : simple modules commuting with a maximal real set}.
\end{proof}

The following theorem is a monoidal analog of
\cite[Theorem 3.3]{BMRRT} and \cite[Theorem 2.4]{DK08}.

\begin{Thm}
\label{Thm : cardinality of maximal real sets}
Any maximal real set has exactly $N$ elements. 
Moreover, for a maximal real set $\{L_i \mid 1 \leq i \leq N\}$ and $t \in \mathbb{T}$, 
\[
\{\deg^t\iota([L_i]) \mid 1 \leq i \leq N\} \subset \Z^N
\]
is a $\Z$-basis of $\Z^N$.

\end{Thm}

\begin{proof}
It is immediate from Proposition~\ref{Prop : cardinality of Laurent families}, 
Propositoin~\ref{Prop : degrees of Laurent families form an integral basis}, 
and Theorem~\ref{Thm : Laurent families and maximal real sets}. 
\end{proof}

\subsection{Exchange pair of maximal real sets}
\label{ssec : exchange pairs of maximal real sets}

\begin{Def}
\label{Def : almost complete maximal real sets}
A commuting family $\overline{\mathcal M}$ of real prime modules in
$\mathscr{C}$ is called an \emph{almost complete maximal real set} if
$|\overline{\mathcal M}| = N - 1$.
\end{Def}

\begin{Lem}
\label{Lem : exchange polynomial for maximal real sets}
Let $\overline{\mathcal M} = \{L_i \mid 2 \leq i \leq N\}$ be an almost
complete maximal real set.
Suppose that there exist real prime modules $L, L'$ such that
$L \not\simeq L'$ and $\{L\} \cup \overline{\mathcal M}$,
$\{L'\} \cup \overline{\mathcal M}$ are both maximal real sets.
Then there exists
$P(u_2, \ldots, u_N) \in \Z_{\geq 0}[u_2, \ldots, u_N]$ such that
\[
[L][L'] = P([L_2], \ldots, [L_N]).
\]
\end{Lem}

\begin{proof}
Since $\{L\} \cup \overline{\mathcal M}$ is a Laurent family,
$[L']$ has a Laurent expansion in
\[
\{[L], [L_2], \ldots, [L_N]\}.
\]
By Proposition~\ref{Prop : denominators and commutativity for Laurent families},
we have
\[
[L'] = \frac{Q([L], [L_2], \ldots, [L_N])}{[L]^k}
\]
for $k \in \Z_{\geq 0}$ and $Q \in \Z_{\geq 0}[u_1, \ldots, u_N]$
which is indivisible by $u_1$.
First we prove that $k = 1$.
If $k = 0$, then $L'$ commutes with $L$ by
Proposition~\ref{Prop : denominators and commutativity for Laurent families}.
Since $L'$ commutes with $L_2, \ldots, L_N$ by the assumption,
$L'$ commutes with $\{L\} \cup \overline{\mathcal M}$.
The maximality implies that $L'$ is a member of
$\{L\} \cup \overline{\mathcal M}$ and it is a contradiction.
Hence we have $k \geq 1$.
Consider the extension degree
\[
\bigl[\Q([L_i] \mid 2 \leq i \leq N)([L]) :
\Q([L_i] \mid 2 \leq i \leq N)([L'])\bigr].
\]
Since $[L']$ is a rational function in $\{[L], [L_2], \ldots, [L_N]\}$
whose denominator has degree $k$ as a polynomial in $[L]$,
the extension degree is at least $k$.
However, the extension degree is $1$ because both fields are
$\operatorname{Frac}(K(\mathscr{C})_{\Z})$ by
Proposition~\ref{Prop : Laurent expansions with respect to a Laurent family}.
Then we have $k = 1$.
By the same argument, we have
\[
[L] = \frac{Q'([L'], [L_2], \ldots, [L_N])}{[L']}
\]
for $Q' \in \Z_{\geq 0}[u_1, u_2, \ldots, u_N]$
which is indivisible by $u_1$.
Then
\[
Q([L], [L_2], \ldots, [L_N])
= Q'([L'], [L_2], \ldots, [L_N])
\]
holds.
By Corollary~\ref{Cor : unique prime tensor factorization in C},
it follows that $Q, Q' \in \Z_{\geq 0}[u_2, \ldots, u_N]$.
\end{proof}

\begin{Prop}
\label{Prop : at most two complements of an almost complete maximal real set}
Let $\overline{\mathcal M} = \{L_2, \ldots, L_N\}$ be an almost complete
maximal real set.
Then there are at most two real prime modules outside
$\overline{\mathcal M}$ which commute with every member of
$\overline{\mathcal M}$.
\end{Prop}

\begin{proof}
Suppose that there exist mutually non-isomorphic real prime modules
$L, L', L''$ which commute with every $L_i \in \overline{\mathcal M}$.

By Lemma~\ref{Lem : exchange polynomial for maximal real sets},
there exist $P', P'' \in \Z_{\geq 0}[u_2, \ldots, u_N]$ such that
\[
[L][L'] = P'([L_i] \mid 2 \leq i \leq N)
\quad\text{and}\quad
[L][L''] = P''([L_i] \mid 2 \leq i \leq N).
\]
Then we have
\[
[L']P''([L_i] \mid 2 \leq i \leq N)
= [L'']P'([L_i] \mid 2 \leq i \leq N).
\]
By the assumption, $L'$ is not isomorphic to $L''$ or $L_i$
for $2 \leq i \leq N$.
By Corollary~\ref{Cor : unique prime tensor factorization in C},
it is a contradiction.
\end{proof}

\begin{Def}
\label{Def : exchange pairs of maximal real sets}
Let $L, L' \in \Irr\mathscr{C}$ be non-isomorphic real prime modules.
We call $(L, L')$ an \emph{exchange pair} if there is an almost complete
maximal real set $\overline{\mathcal M}$ such that
$\{L\} \cup \overline{\mathcal M}$ and
$\{L'\} \cup \overline{\mathcal M}$ are both maximal real sets.
\end{Def}

\begin{Lem}
\label{Lem : unimodularity of PBW and tilde Lambda matrices}
Let $\mathcal M = \{L_i \mid 1 \leq i \leq N\}$ be a maximal real set.
Let
\[
P_{\mathcal M} \coloneq
(\mathrm{PBW}(L_i))_{1 \leq i \leq N} \in M_N(\Z).
\]
\begin{enumerate}
\item[(i)] $P_{\mathcal M} \in \mathrm{GL}_N(\Z)$.
\item[(ii)] $\widetilde{\Lambda}_{\mathcal M} \in \mathrm{GL}_N(\Z)$.
\end{enumerate}
\end{Lem}

\begin{proof}
(i) We can prove this similarly by replacing $\deg^{\mathcal M_t}$
by $\mathrm{PBW}$ and $M_{i;t}$ by $S_i$ in the proof of
Proposition~\ref{Prop : degrees of Laurent families form an integral basis}.

(ii) By Lemma~\ref{Lem : Lambda and PBW}, we have
\[
\Lambda_{\mathcal M}
= P_{\mathcal M}^{T}(\lambda_{kl})_{1 \leq k, l \leq N}P_{\mathcal M}.
\]
By Lemma~\ref{Lem : Lambda infinity tensor products and subquotients},
\[
\begin{aligned}
\Lambda^{\infty}(L_i, L_j)
&= \Lambda^{\infty}\left(
\bigotimes_{k=1}^{N} S_k^{\otimes\mathrm{PBW}(L_i)_k},
\bigotimes_{k=1}^{N} S_k^{\otimes\mathrm{PBW}(L_j)_k}
\right) \\
&= \mathrm{PBW}(L_i)^{T}
(\Lambda^{\infty}(S_k, S_l))_{1 \leq k, l \leq N}
\mathrm{PBW}(L_j).
\end{aligned}
\]
We have
\[
\Lambda^{\infty}_{\mathcal M}
= P_{\mathcal M}^{T}(- (\beta_k, \beta_l))_{1 \leq k, l \leq N}
P_{\mathcal M}
\]
by Lemma~\ref{Lem : Lambda infinity of cuspidal modules}.
Then we obtain that
\[
\widetilde{\Lambda}_{\mathcal M}
= P_{\mathcal M}^{T}
(\delta_{kl} + \delta(k < l)(\beta_k, \beta_l))_{1 \leq k, l \leq N}
P_{\mathcal M}.
\]
Since $(\delta_{kl} + \delta(k < l)(\beta_k, \beta_l))_{1 \leq k, l \leq N}$
is unitriangular, it is in $\mathrm{GL}_N(\Z)$.
Then $\widetilde{\Lambda}_{\mathcal M} \in \mathrm{GL}_N(\Z)$ by (i).
\end{proof}

\begin{Thm}
\label{Thm : exact sequence for an exchange pair}
Let $(L, L')$ be an exchange pair associated to an almost complete
maximal real set
$\overline{\mathcal{M}} = \{L_i \mid 2 \leq i \leq N\}$.
Let $\mathcal{M} \coloneq \{L_1 \coloneq L, L_2, \ldots, L_N\}$.
\begin{enumerate}
    \item[(i)] $\mathfrak{d}(L, L') = 1$.
    \item[(ii)] Let
    $\mathbf{v}_1 \coloneq -\widetilde{\Lambda}_{\mathcal{M}}^{-1}e_1$.
    Then there is an exact sequence
    \[
    0 \longrightarrow L^{[\mathbf{v}_1]_+}
    \longrightarrow L \otimes L'
    \longrightarrow L^{[-\mathbf{v}_1]_+}
    \longrightarrow 0.
    \]
\end{enumerate}
\end{Thm}

\begin{proof}
By Lemma~\ref{Lem : exchange polynomial for maximal real sets},
\[
[L \otimes L'] = \sum_{s=1}^{\ell}c_s[L^{\mathbf{b}(s)}]
\]
for $c_s \in \Z_{>0}$, $\mathbf{b}(s) \in \Z_{\geq 0}^N$,
and $\mathbf{b}(s)_1 = 0$.
Let $L \nabla L' \simeq L^{\mathbf{h}}$ and
$L \Delta L' \simeq L^{\mathbf{s}}$.
Note that the coefficients of $[L^{\mathbf{h}}]$ and $[L^{\mathbf{s}}]$
are $1$ by Lemma~\ref{Lem : Lambda inequalities for heads and socles}.
By Proposition~\ref{Prop : d vanishing and commutativity}
and Lemma~\ref{Lem : Lambda inequalities for heads and socles},
\[
\begin{aligned}
\bigl(\Lambda_{\mathcal{M}}(\mathbf{s} - \mathbf{h})\bigr)_1
&= \Lambda(L, L \Delta L') - \Lambda(L, L \nabla L') \\
&= -\Lambda(L \Delta L', L) - \Lambda(L, L') \\
&= -\Lambda(L', L) - \Lambda(L, L') \\
&= -2\mathfrak{d}(L, L').
\end{aligned}
\]
Since $L_i$ commutes with $L$ and $L'$ for any $2 \leq i \leq N$,
\[
\begin{aligned}
\bigl(\Lambda_{\mathcal{M}}(\mathbf{s} - \mathbf{h})\bigr)_i
&= \Lambda(L_i, L \Delta L') - \Lambda(L_i, L \nabla L') \\
&= \bigl(\Lambda(L_i, L) + \Lambda(L_i, L')\bigr)
 - \bigl(\Lambda(L_i, L) + \Lambda(L_i, L')\bigr) \\
&= 0
\end{aligned}
\]
by Proposition~\ref{Prop : commutativity with simple subquotients}.
By Lemma~\ref{Lem : Lambda infinity tensor products and subquotients},
$\Lambda^{\infty}_{\mathcal{M}}(\mathbf{s} - \mathbf{h}) = 0$.
Hence
\[
\widetilde{\Lambda}_{\mathcal{M}}(\mathbf{s} - \mathbf{h})
= -\mathfrak{d}(L, L')e_1.
\]
By Proposition~\ref{Prop : Lambda inequalities for Laurent expansions},
for $1 \leq s \leq \ell$ such that
$\mathbf{b}(s) \neq \mathbf{s}, \mathbf{h}$,
there exists $a_s \in \Z$ with $0 < a_s < \mathfrak{d}(L, L')$
such that
\[
\widetilde{\Lambda}_{\mathcal{M}}(\mathbf{b}(s) - \mathbf{h})
= -a_s e_1.
\]
Multiplying by $\widetilde{\Lambda}_{\mathcal{M}}^{-1}$, we have
\[
\mathbf{s} - \mathbf{h} = \mathfrak{d}(L, L')\mathbf{v}_1,
\qquad
\mathbf{b}(s) - \mathbf{h} = a_s\mathbf{v}_1.
\]
Note that $(\mathbf{v}_1)_1 = 0$
because $(\mathbf{s} - \mathbf{h})_1 = 0$.
Let
\[
P(t) \coloneq 1
+ \sum_{\mathbf{b}(s) \neq \mathbf{s}, \mathbf{h}}c_s t^{a_s}
+ t^{\mathfrak{d}(L, L')}
\in \Z_{\geq 0}[t].
\]
Then
\[
[L \otimes L'] = [L^{\mathbf{h}}]
P\left(\prod_{i=1}^N[L_i]^{(\mathbf{v}_1)_i}\right)
\]
by the above argument.
Let
\[
P(t) = \prod_{j=1}^{\mathfrak{d}(L, L')}(t - \alpha_j)
\qquad \text{in $\C[t]$}.
\]
Then
\[
[L][L'] = [L^{\mathbf{h}}]
\prod_{j=1}^{\mathfrak{d}(L, L')}
\left(\prod_{i=1}^N[L_i]^{(\mathbf{v}_1)_i} - \alpha_j\right)
\]
in $\C \otimes K(\mathscr{C})_{\Z}
[[L_i]^{-1} \mid 2 \leq i \leq N]$.

(i) Since $(L, L')$ is an exchange pair,
$\mathfrak{d}(L, L') > 0$ is clear.
Since $[L']$ and $[L_i]$ for $1 \leq i \leq N$ are prime elements
in $\C \otimes K(\mathscr{C})_{\Z}$ and
$[L']$ is not associate to $[L_i]$ for any $1 \leq i \leq N$
by Lemma~\ref{Lem : real prime modules give prime elements over C},
$[L']$ is a prime element in
$\C \otimes K(\mathscr{C})_{\Z}
[[L_i]^{-1} \mid 1 \leq i \leq N]$.
In $\C \otimes K(\mathscr{C})_{\Z}
[[L_i]^{-1} \mid 1 \leq i \leq N]$, we have
\[
[L'] = \frac{[L^{\mathbf{h}}]}{[L_1]}
\prod_{j=1}^{\mathfrak{d}(L, L')}
\left(\prod_{i=1}^N[L_i]^{(\mathbf{v}_1)_i} - \alpha_j\right).
\]
Since $\mathbf{v}_1 \neq 0$ and $\alpha_j \in \C^{\times}$,
$\prod_{i=1}^N[L_i]^{(\mathbf{v}_1)_i} - \alpha_j$
is not invertible in
$\C \otimes K(\mathscr{C})_{\Z}
[[L_i]^{-1} \mid 1 \leq i \leq N]$.
Then the right-hand side is a product of at least
$\mathfrak{d}(L, L')$ prime elements.
Therefore $\mathfrak{d}(L, L') \leq 1$.

(ii) We have
\[
[L \otimes L'] = [L^{\mathbf{h}}]
\left(1 + \prod_{i=1}^N[L_i]^{(\mathbf{v}_1)_i}\right)
\]
by (i).
Since $\mathbf{h} + \mathbf{v}_1 = \mathbf{s}$,
$\mathbf{h}, \mathbf{h} + \mathbf{v}_1 \in \Z_{\geq 0}^N$.
Note that $h_1 = 0$
by Lemma~\ref{Lem : exchange polynomial for maximal real sets}.
Suppose that $\min\{h_i, (\mathbf{v}_1 + \mathbf{h})_i\} > 0$
for some $2 \leq i \leq N$.
Then $[L_i] \mid [L][L']$ in $K(\mathscr{C})_{\Z}$.
However, $[L_i]$, $[L]$, and $[L']$ are non-associate prime elements
in $K(\mathscr{C})_{\Z}$, which is a contradiction.
Therefore $\min\{h_i, (\mathbf{v}_1 + \mathbf{h})_i\} = 0$
for any $1 \leq i \leq N$.
Thus we have $\mathbf{h} = [-\mathbf{v}_1]_+$ and
$\mathbf{v}_1 + \mathbf{h} = [\mathbf{v}_1]_+$.
Then we obtain the exact sequence
\[
0 \longrightarrow L^{[\mathbf{v}_1]_+}
\longrightarrow L \otimes L'
\longrightarrow L^{[-\mathbf{v}_1]_+}
\longrightarrow 0.
\]
\end{proof}

The following theorem is a monoidal analog of
\cite[Theorem 7.5]{BMRRT}, \cite[Theorem 5.8]{FG}.

\begin{Thm}
\label{Thm : characterization of exchange pairs}
Let $L, L' \in \Irr\mathscr{C}$ be non-isomorphic real prime modules.
Then $(L, L')$ is an exchange pair if and only if
$\mathfrak{d}(L, L') = 1$.
\end{Thm}

\begin{proof}
The only if implication is shown in
Theorem~\ref{Thm : exact sequence for an exchange pair}.
Assume that $\mathfrak{d}(L, L') = 1$.
By Lemma~\ref{Lem : d one exact sequence},
\[
0 \longrightarrow L \Delta L'
\longrightarrow L \otimes L'
\longrightarrow L \nabla L'
\longrightarrow 0
\]
is exact and $L \Delta L'$ and $L \nabla L'$ are both real.
We have
\[
\mathfrak{d}(L, L \Delta L')
= \mathfrak{d}(L, L \nabla L')
= \mathfrak{d}(L', L \Delta L')
= \mathfrak{d}(L', L \nabla L') = 0
\]
by Lemma~\ref{Lem : decrease of d}.
In addition, we have
$\mathfrak{d}(L \Delta L', L \nabla L') = 0$
by Lemma~\ref{Lem : Lambda tensor products and subquotients}.
Thus $L \otimes (L \Delta L') \otimes (L \nabla L')$ is a real module. 
Let 
\[
\mathcal{S} \coloneq
\{\text{prime factors of }
L \otimes (L \Delta L') \otimes (L \nabla L')\}
\]
be a commuting family of real prime modules in $\mathscr{C}$.
Let $\mathcal{M} = \{L_1 \coloneq L, L_2, \ldots, L_N\}$
be a maximal real set which contains $\mathcal{S}$.
Since $L \Delta L'$ and $L \nabla L'$ are monomials in $\mathcal{M}$,
\[
[L'] = \frac{[L^{\mathbf{a}}] + [L^{\mathbf{b}}]}{[L]}
\qquad \text{for $\mathbf{a}, \mathbf{b} \in \Z_{\geq 0}^N$}.
\]
By Proposition~\ref{Prop : denominators and commutativity for Laurent families},
$L'$ commutes with any $L_i$ for $2 \leq i \leq N$.
Thus $\{L'\} \cup \{L_i \mid 2 \leq i \leq N\}$
is a maximal real set.
Then we have that $(L, L')$ is an exchange pair.
\end{proof}

\subsection{$\Lambda$-admissible structure of maximal real sets}
\label{ssec : Lambda admissible structures of maximal real sets}

Let $J \coloneq \{1, \ldots, N\}$, and let
$J_{\mathrm{uf}}$ and $J_{\mathrm{fr}}$ be as in
Subsection~\ref{ssec : Lambda cluster algebra structure}.
By Proposition~\ref{Prop : frozen commutativity} and
Proposition~\ref{Prop : cluster monomials and real prime modules},
any maximal real set contains the real prime modules corresponding to
the frozen variables.
In what follows, by changing the labeling if necessary, we may assume
that each maximal real set $\mathcal M = \{L_i \mid i \in J\}$ satisfies
$L_j \simeq M_{j;t_0}$ for all $j \in J_{\mathrm{fr}}$.

For a maximal real set $\mathcal M = \{L_i \mid i \in J\}$, let
\[
\widetilde B_{\mathcal M}
\coloneq -\widetilde\Lambda_{\mathcal M}^{-1}
\begin{pmatrix} I_{|J_{\mathrm{uf}}|} \\ 0 \end{pmatrix}.
\]
Here,
\[
\begin{pmatrix} I_{|J_{\mathrm{uf}}|} \\ 0 \end{pmatrix}
= (\delta_{ij})_{i \in J,\,j \in J_{\mathrm{uf}}}.
\]

\begin{Prop}
\label{Prop : exchange matrices of monoidal clusters}
Let $t \in \mathbb T$, and let
$\mathcal M_t = \{M_{i;t} \mid i \in J\}$ be the corresponding
monoidal cluster.
Then $\widetilde B_t = \widetilde B_{\mathcal M_t}$.
\end{Prop}

\begin{proof}
Let $k \in J_{\mathrm{uf}}$.
Since $(M_{k;t}, M_{k;\mu_k t})$ is an exchange pair corresponding to
$\{M_{i;t} \mid i \in J \setminus \{k\}\}$,
\[
0 \longrightarrow
M_t^{[\widetilde B_{\mathcal M_t}e_k]_+}
\longrightarrow M_{k;t} \otimes M_{k;\mu_k t}
\longrightarrow M_t^{[-\widetilde B_{\mathcal M_t}e_k]_+}
\longrightarrow 0
\]
is an exact sequence by
Theorem~\ref{Thm : exact sequence for an exchange pair}.
By Definition~\ref{Def : Lambda admissible monoidal seed},
\[
0 \longrightarrow M_t^{[\widetilde B_t e_k]_+}
\longrightarrow M_{k;t} \otimes M_{k;\mu_k t}
\longrightarrow M_t^{[-\widetilde B_t e_k]_+}
\longrightarrow 0
\]
is also an exact sequence.
By Proposition~\ref{Prop : simplicity of heads and socles},
these exact sequences are non-split.
Then we see that
\[
M_t^{[\widetilde B_{\mathcal M_t}e_k]_+}
\simeq M_t^{[\widetilde B_t e_k]_+},
\qquad
M_t^{[-\widetilde B_{\mathcal M_t}e_k]_+}
\simeq M_t^{[-\widetilde B_t e_k]_+}.
\]
By Corollary~\ref{Cor : unique prime tensor factorization in C},
$\widetilde B_{\mathcal M_t}e_k = \widetilde B_t e_k$ holds.
\end{proof}

\begin{Prop}
\label{Prop : compatibility for maximal real sets}
Let $\mathcal M = \{L_i \mid i \in J\}$ be a maximal real set.
\begin{enumerate}
    \item[(i)] $\widetilde B_{\mathcal M}$ is an exchange matrix.
    \item[(ii)] $(\widetilde B_{\mathcal M}, \Lambda_{\mathcal M})$
    is a compatible pair of type $2I_{|J_{\mathrm{uf}}|}$.
    \item[(iii)]
    $\Lambda^{\infty}_{\mathcal M}\widetilde B_{\mathcal M} = 0$.
\end{enumerate}
\end{Prop}

\begin{proof}
Let $P_{\mathcal M} = (\mathrm{PBW}(L_i))_{i \in J}$ and
$P_{t_0} = (\mathrm{PBW}(M_{i;t_0}))_{i \in J}$.
Note that the claim is clear if $\mathcal M$ is a monoidal cluster.

(i) By the proof of
Lemma~\ref{Lem : unimodularity of PBW and tilde Lambda matrices},
we have
\[
\widetilde\Lambda_{\mathcal M}
= P_{\mathcal M}^{T}(P_{t_0}^{T})^{-1}
\widetilde\Lambda_{t_0}P_{t_0}^{-1}P_{\mathcal M}.
\]
Thus
\[
\begin{aligned}
\left.\widetilde B_{\mathcal M}\right|_
{J_{\mathrm{uf}} \times J_{\mathrm{uf}}}
&= (I_{|J_{\mathrm{uf}}|}\ 0)\widetilde B_{\mathcal M} \\
&= -(I_{|J_{\mathrm{uf}}|}\ 0)\widetilde\Lambda_{\mathcal M}^{-1}
\begin{pmatrix} I_{|J_{\mathrm{uf}}|} \\ 0 \end{pmatrix} \\
&= -\bigl((I_{|J_{\mathrm{uf}}|}\ 0)P_{\mathcal M}^{-1}P_{t_0}\bigr)
\widetilde\Lambda_{t_0}^{-1}
\bigl((I_{|J_{\mathrm{uf}}|}\ 0)P_{\mathcal M}^{-1}P_{t_0}\bigr)^{T}.
\end{aligned}
\]
Since $L_j \simeq M_{j;t_0}$ for $j \in J_{\mathrm{fr}}$,
$P_{t_0}e_j = P_{\mathcal M}e_j$ for $j \in J_{\mathrm{fr}}$.
Therefore
\[
\left.(P_{\mathcal M}^{-1}P_{t_0})\right|_
{J \times J_{\mathrm{fr}}}
= \begin{pmatrix} 0 \\ I_{|J_{\mathrm{fr}}|} \end{pmatrix}.
\]
Since $P_{\mathcal M}^{-1}P_{t_0} \in \mathrm{GL}_N(\Z)$, we have
\[
A_{\mathrm{uf}} \coloneq
\left.(P_{\mathcal M}^{-1}P_{t_0})\right|_
{J_{\mathrm{uf}} \times J_{\mathrm{uf}}}
\in \mathrm{GL}_{|J_{\mathrm{uf}}|}(\Z).
\]
Then
\[
(I_{|J_{\mathrm{uf}}|}\ 0)(P_{\mathcal M}^{-1}P_{t_0})
= (A_{\mathrm{uf}}\ 0)
= A_{\mathrm{uf}}(I_{|J_{\mathrm{uf}}|}\ 0).
\]
Hence it holds that
\[
\begin{aligned}
\left.\widetilde B_{\mathcal M}\right|_
{J_{\mathrm{uf}} \times J_{\mathrm{uf}}}
&= -\bigl(A_{\mathrm{uf}}(I_{|J_{\mathrm{uf}}|}\ 0)\bigr)
\widetilde\Lambda_{t_0}^{-1}
\begin{pmatrix} I_{|J_{\mathrm{uf}}|} \\ 0 \end{pmatrix}
A_{\mathrm{uf}}^{T} \\
&= A_{\mathrm{uf}}(I_{|J_{\mathrm{uf}}|}\ 0)
\widetilde B_{t_0}A_{\mathrm{uf}}^{T}.
\end{aligned}
\]
Because $\widetilde B_{t_0}$ is an exchange matrix,
$(I_{|J_{\mathrm{uf}}|}\ 0)\widetilde B_{t_0}$ is skew-symmetric,
and thus
$\left.\widetilde B_{\mathcal M}\right|_
{J_{\mathrm{uf}} \times J_{\mathrm{uf}}}$ is skew-symmetric.

(ii) By the proof of
Lemma~\ref{Lem : unimodularity of PBW and tilde Lambda matrices},
we have
\[
\Lambda_{\mathcal M}
= P_{\mathcal M}^{T}(P_{t_0}^{T})^{-1}
\Lambda_{t_0}P_{t_0}^{-1}P_{\mathcal M}.
\]
Thus it holds that
\[
\begin{aligned}
\widetilde B_{\mathcal M}^{T}\Lambda_{\mathcal M}
&= -(I_{|J_{\mathrm{uf}}|}\ 0)
(\widetilde\Lambda_{\mathcal M}^{T})^{-1}\Lambda_{\mathcal M}
\\
&= -(I_{|J_{\mathrm{uf}}|}\ 0)(P_{\mathcal M}^{-1}P_{t_0})
(\widetilde\Lambda_{t_0}^{T})^{-1}
(P_{\mathcal M}^{-1}P_{t_0})^{T} 
 (P_{t_0}^{-1}P_{\mathcal M})^{T}
\Lambda_{t_0}(P_{t_0}^{-1}P_{\mathcal M}) \\
&= -(I_{|J_{\mathrm{uf}}|}\ 0)(P_{\mathcal M}^{-1}P_{t_0})
(\widetilde\Lambda_{t_0}^{T})^{-1}
\Lambda_{t_0}(P_{t_0}^{-1}P_{\mathcal M}).
\end{aligned}
\]
Note that $P_{\mathcal M}^{-1}P_{t_0}$ is of the form
\[
\begin{pmatrix} A_{\mathrm{uf}} & 0 \\ * & I_{|J_{\mathrm{fr}}|} \end{pmatrix}.
\]
Hence we have
\[
\begin{aligned}
&-(I_{|J_{\mathrm{uf}}|}\ 0)(P_{\mathcal M}^{-1}P_{t_0})
(\widetilde\Lambda_{t_0}^{T})^{-1}
\Lambda_{t_0}(P_{t_0}^{-1}P_{\mathcal M}) \\
&\quad= -A_{\mathrm{uf}}(I_{|J_{\mathrm{uf}}|}\ 0)
(\widetilde\Lambda_{t_0}^{T})^{-1}\Lambda_{t_0}
\begin{pmatrix}
A_{\mathrm{uf}}^{-1} & 0 \\ * & I_{|J_{\mathrm{fr}}|}
\end{pmatrix} \\
&\quad= A_{\mathrm{uf}}\widetilde B_{t_0}^{T}\Lambda_{t_0}
\begin{pmatrix}
A_{\mathrm{uf}}^{-1} & 0 \\ * & I_{|J_{\mathrm{fr}}|}
\end{pmatrix} \\
&\quad= 2(A_{\mathrm{uf}}\ 0)
\begin{pmatrix}
A_{\mathrm{uf}}^{-1} & 0 \\ * & I_{|J_{\mathrm{fr}}|}
\end{pmatrix} \\
&\quad= (2I_{|J_{\mathrm{uf}}|}\ 0).
\end{aligned}
\]

(iii) By the proof of
Lemma~\ref{Lem : unimodularity of PBW and tilde Lambda matrices},
we have
\[
\Lambda^{\infty}_{\mathcal M}
= (P_{t_0}^{-1}P_{\mathcal M})^{T}\Lambda^{\infty}_{t_0}
(P_{t_0}^{-1}P_{\mathcal M}).
\]
The proof is the same as (ii), except using the fact that
$\Lambda^{\infty}_{t_0}\widetilde B_{t_0} = 0$.
\end{proof}

The goal of this subsection is to prove that
$(\mathcal M, \widetilde B_{\mathcal M})$ is a completely
$\Lambda$-admissible monoidal seed for all maximal real sets
$\mathcal M = \{L_i \mid i \in J\}$.
For $k \in J_{\mathrm{uf}}$, let
$\mathbf v_k \coloneq \widetilde B_{\mathcal M}e_k$.
What remains is to prove the existence of a real prime module $L'_k$
such that
\[
0 \longrightarrow L^{[\mathbf v_k]_+}
\longrightarrow L_k \otimes L'_k
\longrightarrow L^{[-\mathbf v_k]_+}
\longrightarrow 0
\]
is exact.
By Proposition~\ref{Prop : compatibility for maximal real sets},
$(\mathbf v_k)_k = 0$ and
\[
\Lambda(L_k, L^{[\mathbf v_k]_+})
- \Lambda(L_k, L^{[-\mathbf v_k]_+})
= e_k^{T}\Lambda_{\mathcal M}\widetilde B_{\mathcal M}e_k
= -2.
\]

\begin{Lem}
\label{Lem : binomial divisibility by a real simple module}
Let $L \in \Irr\mathscr{C}$ be a real module, and let
$S, H \in \Irr\mathscr{C}$ such that
\[
\Lambda(L, S) - \Lambda(L, H) = -2,
\qquad
\mathfrak{d}(L, S) = \mathfrak{d}(L, H) = 0.
\]
Assume that $[L] \nmid [S], [H]$ in
$\C \otimes_{\Z} K(\mathscr{C})_{\Z}$.
Suppose that there exist
$\phi \in \C \otimes_{\Z} K(\mathscr{C})_{\Z}$
and $\alpha \in \C^{\times}$ such that
\[
[L]\phi = [S] - \alpha[H]
\]
in $\C \otimes_{\Z} K(\mathscr{C})_{\Z}$.
Then we have $\alpha = -1$, and there exists
$L' \in \Irr\mathscr{C}$ such that $\phi = [L']$.
Moreover, $\mathfrak{d}(L, L') = 1$, and it holds that
\[
0 \longrightarrow S \longrightarrow L \otimes L'
\longrightarrow H \longrightarrow 0
\]
is an exact sequence.
\end{Lem}

\begin{proof}
We prove it in a similar way to
\cite[Theorem 4.1]{KK19} and \cite[Theorem 4.13]{KKOP20}.
Let
\[
\phi = \sum_{s \in K} c_s[X_s]
\]
for $c_s \in \C^{\times}$ and $X_s \in \Irr\mathscr{C}$.
It follows that
\[
\sum_{s \in K} c_s[L \otimes X_s] = [S] - \alpha[H].
\]
Take $s_0 \in K$ such that
\[
\Lambda(X_{s_0}, L)
= \max\{\Lambda(X_s, L) \mid s \in K\}.
\]
Suppose that there exists $s_0' \in K \setminus\{s_0\}$
such that $L \Delta X_{s_0}$ is a simple subquotient of
$L \otimes X_{s_0'}$.
Note that $L \Delta X_{s_0'} \not\simeq L \Delta X_{s_0}$
by Lemma~\ref{Lem : hd soc injectivity}.
By Lemma~\ref{Lem : Lambda inequalities for heads and socles},
it follows that
\[
\begin{aligned}
\Lambda(X_{s_0}, L)
&= \Lambda(L \Delta X_{s_0}, L) \\
&< \Lambda(L \Delta X_{s_0'}, L) \\
&= \Lambda(X_{s_0'}, L).
\end{aligned}
\]
It is a contradiction to the choice of $s_0$.
Thus we have $L \Delta X_{s_0} \simeq S$ or $H$,
and the coefficient of $[L \Delta X_{s_0}]$
in the right-hand side is $c_{s_0}$.
Since
\[
\begin{aligned}
&\max\{\Lambda(T, L) \mid
T \text{ is a simple subquotient of } X_s \otimes L
\text{ for some } s \in K\} \\
&\qquad = \Lambda(L \Delta X_{s_0}, L)
\end{aligned}
\]
by Lemma~\ref{Lem : Lambda inequalities for heads and socles}
and $\Lambda(S, L) = \Lambda(H, L) + 2$ by the assumption,
$L \Delta X_{s_0} \simeq S$ holds.
Thus we have $c_{s_0} = 1$.
Since $[L] \nmid [S]$ by the assumption,
we have $\mathfrak{d}(L, X_{s_0}) \geq 1$.
Take $t_0 \in K$ such that
\[
\Lambda(L, X_{t_0})
= \max\{\Lambda(L, X_s) \mid s \in K\}.
\]
By the same argument, we have $L \nabla X_{t_0} \simeq H$
and $c_{t_0} = -\alpha$.
Suppose that there exists $t_1 \neq t_0$ such that
$\Lambda(L, X_{t_1}) = \Lambda(L, X_{t_0})$.
Then $L \nabla X_{t_0} \simeq L \nabla X_{t_1} \simeq H$
holds, and it is a contradiction to
Lemma~\ref{Lem : hd soc injectivity}.

Let us prove that $s_0 = t_0$.
It is enough to prove
$\Lambda(L, X_{s_0}) = \Lambda(L, X_{t_0})$.
By Lemma~\ref{Lem : Lambda inequalities for heads and socles}
and the assumption, we have
\[
\begin{aligned}
\Lambda(L, X_{s_0})
&= 2\mathfrak{d}(L, X_{s_0}) - \Lambda(X_{s_0}, L) \\
&= 2\mathfrak{d}(L, X_{s_0})
   - \Lambda(L \Delta X_{s_0}, L) \\
&= 2\mathfrak{d}(L, X_{s_0}) - \Lambda(S, L) \\
&= 2\mathfrak{d}(L, X_{s_0}) + \Lambda(L, S) \\
&\geq 2 + \Lambda(L, S) \\
&= \Lambda(L, H), \\
\Lambda(L, X_{t_0})
&= \Lambda(L, L \nabla X_{t_0}) = \Lambda(L, H).
\end{aligned}
\]
By the choice of $t_0$, we have
$\Lambda(L, X_{s_0}) = \Lambda(L, X_{t_0})$.
Therefore, we have $-\alpha = c_{t_0} = c_{s_0} = 1$ and
\[
\begin{aligned}
\mathfrak{d}(L, X_{s_0})
&= \frac{1}{2}\bigl(
\Lambda(L, L \nabla X_{s_0})
+ \Lambda(L \Delta X_{s_0}, L)\bigr) \\
&= \frac{1}{2}\bigl(\Lambda(L, H) + \Lambda(S, L)\bigr) \\
&= \frac{1}{2}\bigl(\Lambda(L, H) - \Lambda(L, S)\bigr) \\
&= 1.
\end{aligned}
\]
By Lemma~\ref{Lem : d one exact sequence},
we have an exact sequence
\[
0 \longrightarrow S \longrightarrow L \otimes X_{s_0}
\longrightarrow H \longrightarrow 0.
\]
\end{proof}

\begin{Lem}
\label{Lem : Lambda admissibility and binomial divisibility}
Let $\mathcal{M} = \{L_i \mid i \in J\}$ be a maximal real set.
Then the following conditions are equivalent:
\begin{enumerate}
    \item[(i)] $(\mathcal{M}, \widetilde{B}_{\mathcal{M}})$
    is a $\Lambda$-admissible monoidal seed.
    \item[(ii)] For each $k \in J_{\mathrm{uf}}$,
    there exists $\alpha_k \in \C^{\times}$ such that
    \[
    \frac{[L^{[\mathbf{v}_k]_+}]
    - \alpha_k[L^{[-\mathbf{v}_k]_+}]}{[L_k]}
    \in \C \otimes K(\mathscr{C})_{\Z}.
    \]
\end{enumerate}
\end{Lem}

\begin{proof}
The implication (i) $\Rightarrow$ (ii) follows from
Definition~\ref{Def : Lambda admissible monoidal seed}.
Assume (ii).
By Lemma~\ref{Lem : binomial divisibility by a real simple module},
there exists $L'_k \in \Irr\mathscr{C}$ such that
$\mathfrak{d}(L_k, L'_k) = 1$ and
\[
0 \longrightarrow L^{[\mathbf{v}_k]_+}
\longrightarrow L_k \otimes L'_k
\longrightarrow L^{[-\mathbf{v}_k]_+}
\longrightarrow 0
\]
is exact.
By Proposition~\ref{Prop : denominators and commutativity for Laurent families},
$L'_k$ commutes with all $L_i$ for $i \neq k$.
Therefore it is enough to prove that $L'_k$ is real prime.

First we prove the primeness of $L'_k$.
Suppose $L'_k$ is not prime.
Since $\mathfrak{d}(L_k, L'_k) = 1$,
there exists a prime factor $P$ such that $\mathfrak{d}(L_k, P) = 0$.
Note that we can write $L'_k \simeq P \otimes P'$.
Since $\mathfrak{d}(L_k, P') = 1$,
the length of $L_k \otimes P'$ is at least $2$.
Because $L_k \otimes P \otimes P'$ has length $2$,
we have that the length of $L_k \otimes P'$ is two and
$P$ strongly commutes with both simple subquotients.
Thus we have
\[
[P] \mid [L^{[\mathbf{v}_k]_+}],\ [L^{[-\mathbf{v}_k]_+}].
\]
It contradicts Proposition~\ref{Prop : real prime modules give prime elements}.

Next we prove the reality of $L'_k$.
Suppose $L'_k$ is imaginary.
Since $(\mathbf{v}_k)_k = 0$,
$L'_k$ commutes with $L^{[\mathbf{v}_k]_+}$ and
$L^{[-\mathbf{v}_k]_+}$.
Then we have
\[
[L_k \otimes L'_k \otimes L'_k]
= [L^{[\mathbf{v}_k]_+} \otimes L'_k]
+ [L^{[-\mathbf{v}_k]_+} \otimes L'_k],
\]
and we obtain that $L_k \otimes L'_k \otimes L'_k$ has length $2$.
Because $L'_k$ is imaginary,
the length of $L'_k \otimes L'_k$ is at least $2$.
Since $L_k \otimes L'_k \otimes L'_k$ and $L'_k \otimes L'_k$
have the same length,
$L_k$ commutes with all simple subquotients of $L'_k \otimes L'_k$.
Then we obtain that
\[
[L_k] \mid [L^{[\mathbf{v}_k]_+} \otimes L'_k]
\qquad\text{in $K(\mathscr{C})_{\Z}$}.
\]
By Proposition~\ref{Prop : real prime modules give prime elements}, we have
\[
[L_k] \mid [L'_k]
\qquad\text{in $K(\mathscr{C})_{\Z}$}.
\]
By Proposition~\ref{Prop : divisibility by a real simple module},
there exists $Y \in \Irr\mathscr{C}$ such that
$L'_k \simeq L_k \otimes Y$.
It contradicts the primeness of $L'_k$.
\end{proof}

The following theorem is a monoidal analog of
\cite[Theorem 5.1]{BMRRT}, \cite[Theorem 5.3]{IY08}.

\begin{Thm}
\label{Thm : Lambda admissibility of maximal real sets}
Let $\mathcal{M} = \{L_i \mid i \in J\}$ be a maximal real set.
Then $(\mathcal{M}, \widetilde{B}_{\mathcal{M}})$ is a
$\Lambda$-admissible monoidal seed.
\end{Thm}

\begin{proof}
Let $k \in J_{\mathrm{uf}}$.
First we prove that there is a simple module which does not commute
with $L_k$.
Let $t \in \mathbb{T}$ and $\mathcal{M}_t$ be the corresponding
monoidal cluster.
If $L_k \notin \mathcal{M}_t$, then $L_k$ does not commute with
some element of $\mathcal{M}_t$ because $\mathcal{M}_t$ is a
maximal real set.
Suppose that there exists $i \in J$ such that $L_k \simeq M_{i;t}$.
Since $k \in J_{\mathrm{uf}}$, $L_k$ does not correspond to
frozen variables and we have $i \in J_{\mathrm{uf}}$.
Then $M_{i; \mu_i t}$ does not commute with $L_k$. 

Let $X$ be a simple module which does not commute with $L_k$.
By Proposition~\ref{Prop : denominators and commutativity for Laurent families},
$\mathbf{d}(X)_k > 0$ holds.
Let us consider
\[
\left(\bigotimes_{i \neq k}
L_i^{\otimes\max\{0, \mathbf{d}(X)_i\}}\right) \otimes X.
\]
Since each coefficient of Laurent polynomials corresponding to modules
is positive, there exists some simple subquotient $S$ of
\[
\left(\bigotimes_{i \neq k}
L_i^{\otimes\max\{0, \mathbf{d}(X)_i\}}\right) \otimes X
\]
such that $\mathbf{d}(S)_k > 0$ and $\mathbf{d}(S)_i \leq 0$
for all $i \neq k$.

Let
\[
[L_k^{\otimes\mathbf{d}(S)_k} \otimes S]
= \sum_{s=1}^{\ell}c_s[L^{\mathbf{b}(s)}]
\]
for $c_s \in \Z_{>0}$ and $\mathbf{b}(s) \in \Z_{\geq 0}^N$.
Let
\[
L_k^{\otimes\mathbf{d}(S)_k} \nabla S \simeq L^{\mathbf{h}},
\qquad
L_k^{\otimes\mathbf{d}(S)_k} \Delta S \simeq L^{\mathbf{s}}.
\]
By the same argument as
Theorem~\ref{Thm : exact sequence for an exchange pair}, we have
\[
\mathbf{s} - \mathbf{h}
= \mathfrak{d}(L_k, S)\mathbf{v}_k
\]
and
\[
\mathbf{b}(s) - \mathbf{h} = a_s\mathbf{v}_k
\qquad\text{for }0 < a_s < \mathfrak{d}(L_k, S)
\]
for all $\mathbf{b}(s) \neq \mathbf{s}, \mathbf{h}$.
Thus there exists
\[
P(t) = 1 + \sum_{j=1}^{\mathfrak{d}(L_k, S)-1}a_jt^j
+ t^{\mathfrak{d}(L_k, S)} \in \Z_{\geq 0}[t]
\]
such that
\[
[L_k^{\otimes\mathbf{d}(S)_k} \otimes S]
= [L^{\mathbf{h}}]
P\left(\prod_{i \in J}[L_i]^{(\mathbf{v}_k)_i}\right).
\]
Let
\[
P(t) = \prod_{j=1}^{\mathfrak{d}(L_k, S)}(t - \alpha_j)
\qquad\text{for }\alpha_j \in \C^{\times}.
\]
Then it follows that
\[
[L_k^{\otimes\mathbf{d}(S)_k} \otimes S]
= [L^{\mathbf{h}}]
\prod_{j=1}^{\mathfrak{d}(L_k, S)}
\left(\prod_{i \in J}[L_i]^{(\mathbf{v}_k)_i} - \alpha_j\right).
\]
Multiplying by $[L^{[-\mathbf{v}_k]_+}]^{\mathfrak{d}(L_k, S)}$,
we obtain that
\[
\begin{aligned}
&[L^{\mathfrak{d}(L_k, S)[-\mathbf{v}_k]_+}
\otimes L_k^{\otimes\mathbf{d}(S)_k} \otimes S] \\
&\qquad = [L^{\mathbf{h}}]
\prod_{j=1}^{\mathfrak{d}(L_k, S)}
\left([L^{[\mathbf{v}_k]_+}]
- \alpha_j[L^{[-\mathbf{v}_k]_+}]\right)
\end{aligned}
\]
in $\C \otimes K(\mathscr{C})_{\Z}$.
If $h_k > 0$, then $\mathbf{b}(s)_k > 0$ for any
$1 \leq s \leq \ell$ since
$\mathbf{b}(s) - \mathbf{h} \in \Z\mathbf{v}_k$
and $(\mathbf{v}_k)_k = 0$.
It contradicts
\[
\min_{1 \leq s \leq \ell}\mathbf{b}(s)_k = 0.
\]
Thus $[L_k] \nmid [L^{\mathbf{h}}]$ in
$\C \otimes K(\mathscr{C})_{\Z}$ by
Lemma~\ref{Lem : real prime modules give prime elements over C}.
Therefore
\[
[L_k] \mid
\left([L^{[\mathbf{v}_k]_+}]
- \alpha_j[L^{[-\mathbf{v}_k]_+}]\right)
\]
in $\C \otimes K(\mathscr{C})_{\Z}$ for some
$1 \leq j \leq \mathfrak{d}(L_k, S)$.
Now the claim follows by
Lemma~\ref{Lem : Lambda admissibility and binomial divisibility}.
\end{proof}

\begin{Cor}
\label{Cor : complete Lambda admissibility of maximal real sets}
Let $\mathcal{M} = \{L_i \mid i \in J\}$ be a maximal real set.
Then $(\mathcal{M}, \widetilde{B}_{\mathcal{M}})$ is a completely
$\Lambda$-admissible monoidal seed.
\end{Cor}

\begin{proof}
By Theorem~\ref{Thm : Lambda admissibility of maximal real sets},
it is enough to prove that
\[
\mu_k\widetilde{B}_{\mathcal{M}}
= \widetilde{B}_{\mu_k\mathcal{M}}
\qquad\text{for any }k \in J_{\mathrm{uf}}.
\]
Here
$\mu_k\mathcal{M}
= \{L_i \mid i \in J \setminus \{k\}\} \cup \{L'_k\}$
is the maximal real set with
\[
[L'_k] =
\frac{[L^{[\mathbf{v}_k]_+}] + [L^{[-\mathbf{v}_k]_+}]}{[L_k]}.
\]
Since
\[
\frac{1}{2}
\left(\Lambda_{\mu_k\mathcal{M}}
- \Lambda^{\infty}_{\mu_k\mathcal{M}}\right)
\in \mathrm{GL}_N(\Z),
\]
$\widetilde{B}_{\mu_k\mathcal{M}}$ is characterized by
\[
\widetilde{B}_{\mu_k\mathcal{M}}^T\Lambda_{\mu_k\mathcal{M}}
= (2I_{|J_{\mathrm{uf}}|}\ 0),
\qquad
\Lambda^{\infty}_{\mu_k\mathcal{M}}
\widetilde{B}_{\mu_k\mathcal{M}} = 0.
\]
Let
$\widetilde{B}_{\mathcal{M}}
= (b_{ij})_{i \in J,\, j \in J_{\mathrm{uf}}}$.
As in the proof of \cite[Proposition 6.4]{KKOP20},
we can prove that
\[
\begin{aligned}
\Lambda(L_i, L'_k)
&= -\Lambda(L_i, L_k)
+ \sum_{\ell \in J}\Lambda(L_i, L_\ell)[-b_{\ell k}]_+, \\
\Lambda(L'_k, L_i)
&= -\Lambda(L_k, L_i)
+ \sum_{\ell \in J}\Lambda(L_\ell, L_i)[b_{\ell k}]_+
\end{aligned}
\]
for all $i \neq k$.
Thus we have
$\Lambda_{\mu_k\mathcal{M}} = \mu_k\Lambda_{\mathcal{M}}$.
Therefore
\[
(\mu_k\widetilde{B}_{\mathcal{M}})^T
\Lambda_{\mu_k\mathcal{M}} = (2I_{|J_{\mathrm{uf}}|}\ 0)
\]
holds by the fact in Subsection~\ref{ssec : cluster algebras}.

By \cite[(3.2)]{BFZ05},
\[
\mu_k\widetilde{B}_{\mathcal{M}}
= E_{k,+}^{\widetilde{B}_{\mathcal{M}}}
\widetilde{B}_{\mathcal{M}}
F_{k,+}^{\widetilde{B}_{\mathcal{M}}},
\]
where
\begin{itemize}
    \item
    $E_{k,+}^{\widetilde{B}_{\mathcal{M}}} = (e_{ij})_{i,j \in J}$,
    \[
    e_{ij} =
    \begin{cases}
        \delta_{ij}, & \text{if }j \neq k, \\
        -1, & \text{if }(i, j) = (k, k), \\
        [-b_{ik}]_+, & \text{if }i \neq k,\ j = k,
    \end{cases}
    \]
    \item
    $F_{k,+}^{\widetilde{B}_{\mathcal{M}}}
    = (f_{ij})_{i,j \in J_{\mathrm{uf}}}$,
    \[
    f_{ij} =
    \begin{cases}
        \delta_{ij}, & \text{if }i \neq k, \\
        -1, & \text{if }(i, j) = (k, k), \\
        [b_{kj}]_+, & \text{if }i = k,\ j \neq k.
    \end{cases}
    \]
\end{itemize}
Note that
$\left(E_{k,+}^{\widetilde{B}_{\mathcal{M}}}\right)^2 = I_N$.
By Lemma~\ref{Lem : Lambda infinity tensor products and subquotients},
we have
\[
\begin{aligned}
\Lambda^{\infty}(M, L'_k)
&= -\Lambda^{\infty}(M, L_k)
+ \sum_{\ell \in J}\Lambda^{\infty}(M, L_\ell)[-b_{\ell k}]_+, \\
\Lambda^{\infty}(L'_k, M)
&= -\Lambda^{\infty}(L_k, M)
+ \sum_{\ell \in J}\Lambda^{\infty}(L_\ell, M)[b_{\ell k}]_+
\end{aligned}
\]
for all $M \in \Irr\mathscr{C}$.
Thus we obtain that
\[
\Lambda^{\infty}_{\mu_k\mathcal{M}}
= \left(E_{k,+}^{\widetilde{B}_{\mathcal{M}}}\right)^T
\Lambda^{\infty}_{\mathcal{M}}
E_{k,+}^{\widetilde{B}_{\mathcal{M}}}.
\]
Then it holds that
\[
\begin{aligned}
\Lambda^{\infty}_{\mu_k\mathcal{M}}
(\mu_k\widetilde{B}_{\mathcal{M}})
&= \left(E_{k,+}^{\widetilde{B}_{\mathcal{M}}}\right)^T
\Lambda^{\infty}_{\mathcal{M}}
\left(E_{k,+}^{\widetilde{B}_{\mathcal{M}}}\right)^2
\widetilde{B}_{\mathcal{M}}
F_{k,+}^{\widetilde{B}_{\mathcal{M}}} \\
&= \left(E_{k,+}^{\widetilde{B}_{\mathcal{M}}}\right)^T
\Lambda^{\infty}_{\mathcal{M}}\widetilde{B}_{\mathcal{M}}
F_{k,+}^{\widetilde{B}_{\mathcal{M}}} \\
&= 0.
\end{aligned}
\]
\end{proof}

\section{Tropical $F$-invariants and monoidal categorifications}
\label{Sec : tropical F invariants and monoidal categorifications}

\subsection{Tropical $F$-invariants}
\label{ssec : tropical F invariants}

Let $J = J_{\mathrm{uf}} \sqcup J_{\mathrm{fr}}$ be a finite set.
Let
$\mathcal{S}_{\mathbf{x}}^{\Lambda}
= \{(\mathbf{x}_t, \widetilde{B}_t, \Lambda_t)
\mid t \in \mathbb{T}\}$
be a $\Lambda$-cluster pattern such that
$(\widetilde{B}_t, \Lambda_t)$ is a compatible pair of type
$D = \operatorname{diag}(d_i)_{i \in J_{\mathrm{uf}}}$.

\begin{Def}[{\cite[$\S$4.2]{Cao},
\cite[Definition 3.1]{Fei23}}]
\label{Def : tropical evaluation of Laurent polynomials}
Let $K$ be a finite set.
Let
\[
F(u_i \mid i \in K)
= \sum_{\mathbf{v} \in \Z^K}c_{\mathbf{v}}\mathbf{u}^{\mathbf{v}}
\in \Z[u_i^{\pm 1} \mid i \in K] \setminus \{0\}.
\]
We define the map
\[
F[-] \colon \mathbb{R}^{|K|} \longrightarrow \mathbb{R}
\]
by
\[
F[\mathbf{r}]
\coloneq \max\{\mathbf{v}^T\mathbf{r} \mid c_{\mathbf{v}} \neq 0\}.
\]
\end{Def}

\begin{Def}[{\cite[Definition 3.6]{Cao}}]
\label{Def : good elements}
We call $u \in \mathcal{U}(\mathcal{S}_{\mathbf{x}}^{\Lambda})$
\emph{good} if $u$ is compatibly pointed and
\[
u \in \bigcap_{t \in \mathbb{T}}
\Z_{\geq 0}[x_{i;t}^{\pm 1} \mid i \in J].
\]
\end{Def}

\begin{Prop}[{\cite[Proposition 4.4, Theorem 4.10]{Cao}}]
\label{Prop : tropical invariants of good elements}
Let $u_1, u_2 \in \mathcal{U}(\mathcal{S}_{\mathbf{x}}^{\Lambda})$
be good elements.
Then
\[
\begin{gathered}
(\deg^t u_1)^T\Lambda_t(\deg^t u_2)
+ F_{u_1}^t[(D\ 0)\deg^t u_2], \\
F_{u_1}^t[(D\ 0)\deg^t u_2]
+ F_{u_2}^t[(D\ 0)\deg^t u_1]
\end{gathered}
\]
are independent of the choice of $t \in \mathbb{T}$.
\end{Prop}

In \cite{Cao}, the above two values are defined using
a subtraction-free semifield, and the above result is shown
for good elements.
However, the above two values can be defined for compatibly pointed
elements $u, u'$ without positivity.
In fact, the following holds.

\begin{Prop}
\label{Prop : tropical invariants of compatibly pointed elements}
Let $u_1, u_2 \in \mathcal{U}(\mathcal{S}_{\mathbf{x}}^{\Lambda})$
be compatibly pointed elements.
Then
\[
\begin{gathered}
(\deg^t u_1)^T\Lambda_t(\deg^t u_2)
+ F_{u_1}^t[(D\ 0)\deg^t u_2], \\
F_{u_1}^t[(D\ 0)\deg^t u_2]
+ F_{u_2}^t[(D\ 0)\deg^t u_1]
\end{gathered}
\]
are independent of the choice of $t \in \mathbb{T}$.
\end{Prop}

\begin{proof}
Since $\Lambda_t$ is skew-symmetric,
\[
(\deg^t u_1)^T\Lambda_t(\deg^t u_2)
+ (\deg^t u_2)^T\Lambda_t(\deg^t u_1) = 0
\]
holds for any $t \in \mathbb{T}$.
Therefore, it is enough to show only the first claim.

Let $t \in \mathbb{T}$ and $k \in J_{\mathrm{uf}}$.
Let $t' \coloneq \mu_k t$.
Let
\[
u_1^t(z_i \mid i \in J),\quad
u_1^{t'}(z_i \mid i \in J)
\in \Z[z_i^{\pm 1} \mid i \in J]
\]
be Laurent polynomials such that
\[
u_1 = u_1^t(x_{i;t} \mid i \in J)
= u_1^{t'}(x_{i;t'} \mid i \in J).
\]
Let
\[
u_1^{t'}(z_i \mid i \in J)
= \sum_{\mathbf{v} \in \Z^J}c'_{\mathbf{v}}\mathbf{z}^{\mathbf{v}}
\qquad\text{for }c'_{\mathbf{v}} \in \Z.
\]
For
$[\mathbf{a}] = \{\mathbf{a}_s \mid s \in \mathbb{T}\}
\in S_X(\Z^{\max})$,
let us prove that $u_1^t[\mathbf{a}_t]
= u_1^{t'}[\mathbf{a}_{t'}]$.
We have
\[
\begin{aligned}
u_1
&= \sum_{\mathbf{v} \in \Z^J}c'_{\mathbf{v}}\mathbf{x}_{t'}^{\mathbf{v}} \\
&= \sum_{\mathbf{v} \in \Z^J}c'_{\mathbf{v}}
\prod_{i \neq k}x_{i;t}^{v_i}
\left(
\frac{\mathbf{x}_t^{[\widetilde{B}_t e_k]_+}
+ \mathbf{x}_t^{[-\widetilde{B}_t e_k]_+}}{x_{k;t}}
\right)^{v_k}.
\end{aligned}
\]
Let
\[
m \coloneq
\max\bigl(\{-v_k \mid c'_{\mathbf{v}} \neq 0\} \cup \{0\}\bigr)
\]
and consider
\[
\begin{aligned}
f^t(z_i \mid i \in J)
&\coloneq
\left(
\frac{\mathbf{z}^{[\widetilde{B}_t e_k]_+}
+ \mathbf{z}^{[-\widetilde{B}_t e_k]_+}}{z_k}
\right)^m u_1^t(z_i \mid i \in J) \\
&\in \Z[z_i^{\pm 1} \mid i \in J].
\end{aligned}
\]
Note that, for $\mathbf{a}, \mathbf{b}, \mathbf{c} \in \Z^N$ and $\ell \in \Z_{\geq 0}$,
\[
\max\{\mathbf{a}^T(k\mathbf{b} + (\ell-k)\mathbf{c}) \mid 0 \leq k \leq \ell\}
= \max\{\ell \mathbf{a}^T\mathbf{b},\, \ell \mathbf{a}^T\mathbf{c}\}.
\]
We have
\[
\begin{aligned}
f^t[\mathbf{a}_t]
&\leq \max\left\{
\begin{aligned}
&\mathbf{a}_t^T\left(
\sum_{i \neq k}\bigl(v_i+(m+v_k)[b_{ik}^t]_+\bigr)e_i
-(m+v_k)e_k\right), \\
&\mathbf{a}_t^T\left(
\sum_{i \neq k}\bigl(v_i+(m+v_k)[-b_{ik}^t]_+\bigr)e_i
-(m+v_k)e_k\right)
\end{aligned}
\;\middle|\; c'_{\mathbf{v}} \neq 0
\right\} \\
&= \max\left\{
\begin{aligned}
\sum_{i \neq k}a_i^tv_i
+(m+v_k)\Bigl(-a_k^t
+\max\Bigl\{&\sum_{i \neq k}[b_{ik}^t]_+a_i^t, \\
&\sum_{i \neq k}[-b_{ik}^t]_+a_i^t\Bigr\}\Bigr)
\end{aligned}
\;\middle|\; c'_{\mathbf{v}} \neq 0
\right\} \\
&= m a_k^{t'}
+\max\{\mathbf{a}_{t'}^T\mathbf{v} \mid c'_{\mathbf{v}} \neq 0\} \\
&= m a_k^{t'}+u_1^{t'}[\mathbf{a}_{t'}].
\end{aligned}
\]
Since
\[
\begin{aligned}
f^t[\mathbf{a}_t]
&= \left(
\frac{\mathbf{z}^{[\widetilde{B}_t e_k]_+}
+ \mathbf{z}^{[-\widetilde{B}_t e_k]_+}}{z_k}
\right)^m[\mathbf{a}_t]
+u_1^t[\mathbf{a}_t] \\
&= m a_k^{t'}+u_1^t[\mathbf{a}_t],
\end{aligned}
\]
we have
\[
u_1^t[\mathbf{a}_t]\leq u_1^{t'}[\mathbf{a}_{t'}].
\]
Interchanging the role of $t$ and $t'$, we obtain
\[
u_1^t[\mathbf{a}_t]\geq u_1^{t'}[\mathbf{a}_{t'}].
\]
Thus we have $u_1^t[\mathbf{a}_t]=u_1^{t'}[\mathbf{a}_{t'}]$.

By Lemma~\ref{Lem : tropical points and Lambda}, we have
\[
u_1^t[\Lambda_t\deg^t u_2]
=u_1^{t'}[\Lambda_{t'}\deg^{t'}u_2].
\]
Because
\[
u_1^t(z_i \mid i \in J)
=\mathbf{z}^{\deg^t u_1}
F_{u_1}^t(\mathbf{z}^{\widetilde{B}_t e_s}
\mid s \in J_{\mathrm{uf}}),
\]
we have
\[
\begin{aligned}
u_1^t[\Lambda_t\deg^t u_2]
&=(\deg^t u_1)^T\Lambda_t(\deg^t u_2)
+F_{u_1}^t[\widetilde{B}_t^T\Lambda_t\deg^t u_2] \\
&=(\deg^t u_1)^T\Lambda_t(\deg^t u_2)
+F_{u_1}^t[(D\ 0)\deg^t u_2].
\end{aligned}
\]
Similarly, we have
\[
\begin{aligned}
u_1^{t'}[\Lambda_{t'}\deg^{t'}u_2]
&=(\deg^{t'}u_1)^T\Lambda_{t'}(\deg^{t'}u_2) \\
&\qquad+F_{u_1}^{t'}[(D\ 0)\deg^{t'}u_2].
\end{aligned}
\]
Then the first claim follows.
\end{proof}

\subsection{Tropical $F$-invariants and $\mathfrak{d}$-invariants}
\label{ssec : tropical F invariants and d invariants}

Let
\[
\iota \colon K(\mathscr{C})_{\Z}
\xrightarrow{\sim}
\overline{\mathcal{A}}(\mathcal{S}_{\mathbf{x}}^{\Lambda}),
\qquad J = \{1, \ldots, N\},
\]
$J_{\mathrm{uf}}$, and $J_{\mathrm{fr}}$ be as in
Subsection~\ref{ssec : Lambda cluster algebra structure}.
Note that $(\widetilde{B}_t, \Lambda_t)$ is a compatible pair
of type $2I_{|J_{\mathrm{uf}}|}$.
By Proposition~\ref{Prop : compatible pointedness of simple modules},
$\iota([X]) \in \mathcal{U}(\mathcal{S}_{\mathbf{x}}^{\Lambda})$
is compatibly pointed for $X \in \Irr\mathscr{C}$.

\begin{Thm}[{\cite[Theorem 5.16]{Cao}}]
\label{Thm : tropical F invariants and reachable modules}
Let $L_1, L_2 \in \Irr\mathscr{C}$ and assume that one of them
is reachable.
For $t \in \mathbb{T}$, let
\[
\widehat{\mathbf{g}}_i^t \coloneq \deg^t\iota([L_i]),
\qquad
\mathbf{g}_i^t \coloneq
\operatorname{pr}_{\mathrm{uf}}\widehat{\mathbf{g}}_i^t
\qquad\text{for }i = 1, 2.
\]
Then we have
\[
\begin{aligned}
\Lambda(L_1, L_2)
&= (\widehat{\mathbf{g}}_1^t)^{\mathsf{T}}
\Lambda_t\widehat{\mathbf{g}}_2^t
+ F_{\iota([L_1])}^t[2\mathbf{g}_2^t], \\
\mathfrak{d}(L_1, L_2)
&= F_{\iota([L_1])}^t[\mathbf{g}_2^t]
+ F_{\iota([L_2])}^t[\mathbf{g}_1^t].
\end{aligned}
\]
\end{Thm}

The following inequalities hold without any additional assumption.

\begin{Thm}
\label{Thm : tropical F invariants and Lambda inequalities}
Let $L_1, L_2 \in \Irr\mathscr{C}$ and $t \in \mathbb{T}$.
Let
\[
\widehat{\mathbf{g}}_i^t \coloneq \deg^t\iota([L_i]),
\qquad
\mathbf{g}_i^t \coloneq
\operatorname{pr}_{\mathrm{uf}}\widehat{\mathbf{g}}_i^t
\qquad\text{for }i = 1, 2.
\]
Then we have
\[
\begin{aligned}
\Lambda(L_1, L_2)
&\geq (\widehat{\mathbf{g}}_1^t)^{\mathsf{T}}
\Lambda_t\widehat{\mathbf{g}}_2^t
+ F_{\iota([L_1])}^t[2\mathbf{g}_2^t], \\
\mathfrak{d}(L_1, L_2)
&\geq F_{\iota([L_1])}^t[\mathbf{g}_2^t]
+ F_{\iota([L_2])}^t[\mathbf{g}_1^t].
\end{aligned}
\]
\end{Thm}

\begin{proof}
It is enough to show the first claim.
Let $\{M_{i;t} \mid i \in J\}$ be a monoidal cluster.
Then we can write
\[
[L_1]
= [M_t]^{\widehat{\mathbf{g}}_1^t}
F_{\iota([L_1])}^t
\bigl([M_t]^{\widetilde{B}_t e_k} \mid k \in J_{\mathrm{uf}}\bigr).
\]
Let
\[
F_{\iota([L_1])}^t(u_k \mid k \in J_{\mathrm{uf}})
= \sum_{\mathbf{v} \in \Z_{\geq 0}^{J_{\mathrm{uf}}}}
c_{\mathbf{v}}\mathbf{u}^{\mathbf{v}}
\]
and let $\mathbf{a} \in \Z_{>0}^J$ be such that
\[
\widehat{\mathbf{g}}_1^t + \widetilde{B}_t\mathbf{v} + \mathbf{a}
\in \Z_{\geq 0}^J
\qquad\text{for all $\mathbf{v}$ with $c_{\mathbf{v}} \neq 0$}.
\]
By Lemma~\ref{Lem : Lambda tensor products and subquotients},
it holds that
\[
\begin{aligned}
\Lambda(L_1, L_2) + \Lambda(M_t^{\mathbf{a}}, L_2)
&= \Lambda(L_1 \otimes M_t^{\mathbf{a}}, L_2) \\
&\geq \max\left\{
\Lambda\left(
M_t^{\widehat{\mathbf{g}}_1^t + \widetilde{B}_t\mathbf{v} + \mathbf{a}},
L_2\right)
\;\middle|\; c_{\mathbf{v}} \neq 0
\right\}.
\end{aligned}
\]
Since $M_t^{\mathbf{a}}$ and
$M_t^{\widehat{\mathbf{g}}_1^t + \widetilde{B}_t\mathbf{v} + \mathbf{a}}$
are reachable, we have
\[
\begin{aligned}
\Lambda(M_t^{\mathbf{a}}, L_2)
&= \mathbf{a}^{\mathsf{T}}\Lambda_t\widehat{\mathbf{g}}_2^t, \\
\Lambda\left(
M_t^{\widehat{\mathbf{g}}_1^t + \widetilde{B}_t\mathbf{v} + \mathbf{a}},
L_2\right)
&= (\widehat{\mathbf{g}}_1^t + \widetilde{B}_t\mathbf{v}
+ \mathbf{a})^{\mathsf{T}}\Lambda_t\widehat{\mathbf{g}}_2^t.
\end{aligned}
\]
Here, note that the $F$-polynomial at $t$ of a cluster monomial
at $t$ is $1$.
Then it holds that
\[
\begin{aligned}
&\max\left\{
\Lambda\left(
M_t^{\widehat{\mathbf{g}}_1^t + \widetilde{B}_t\mathbf{v} + \mathbf{a}},
L_2\right)
\;\middle|\; c_{\mathbf{v}} \neq 0
\right\} \\
&\quad = (\widehat{\mathbf{g}}_1^t + \mathbf{a})^{\mathsf{T}}
\Lambda_t\widehat{\mathbf{g}}_2^t
+ \max\left\{
\mathbf{v}^{\mathsf{T}}\widetilde{B}_t^{\mathsf{T}}
\Lambda_t\widehat{\mathbf{g}}_2^t
\;\middle|\; c_{\mathbf{v}} \neq 0
\right\} \\
&\quad = (\widehat{\mathbf{g}}_1^t + \mathbf{a})^{\mathsf{T}}
\Lambda_t\widehat{\mathbf{g}}_2^t
+ F_{\iota([L_1])}^t
[\widetilde{B}_t^{\mathsf{T}}\Lambda_t\widehat{\mathbf{g}}_2^t] \\
&\quad = (\widehat{\mathbf{g}}_1^t + \mathbf{a})^{\mathsf{T}}
\Lambda_t\widehat{\mathbf{g}}_2^t
+ F_{\iota([L_1])}^t[2\mathbf{g}_2^t].
\end{aligned}
\]
\end{proof}

\section{Monoidal reachability for acyclic cluster algebras}
\label{Sec : monoidal reachability for acyclic cluster algebras}

\subsection{Cluster categories of acyclic quivers}
\label{ssec : cluster categories of acyclic quivers}

Let $Q$ be a finite acyclic quiver and $\C Q$ be the path algebra.
Let $\operatorname{mod}\C Q$ be the category of finitely generated
modules over $\C Q$ and $D^b(\operatorname{mod}\C Q)$ be the bounded
derived category of $\operatorname{mod}\C Q$.
Let $\tau$ be the Auslander--Reiten translation in
$D^b(\operatorname{mod}\C Q)$ and $\Sigma$ be the shift functor.

\begin{Def}[{\cite[\S1]{BMRRT}}]
\label{Def : cluster category}
We define the \emph{cluster category}
$\mathcal{C}_Q \coloneq
D^b(\operatorname{mod}\C Q)/(\tau\Sigma^{-1})$ as follows:
\[
\operatorname{obj}(\mathcal{C}_Q)
= \operatorname{obj}(D^b(\operatorname{mod}\C Q)),
\]
\[
\operatorname{Hom}_{\mathcal{C}_Q}(X, Y)
= \bigoplus_{i \in \Z}
\operatorname{Hom}_{D^b(\operatorname{mod}\C Q)}
\bigl((\tau\Sigma^{-1})^iX, Y\bigr).
\]
\end{Def}

\begin{Prop}[{\cite{Kel05, BMRRT}}]
\label{Prop : properties of cluster categories}
Let $Q$ be an acyclic quiver.
\begin{enumerate}
    \item[(i)] $\mathcal{C}_Q$ is a triangulated category.
    \item[(ii)] The canonical functor
    $\pi \colon D^b(\operatorname{mod}\C Q) \longrightarrow
    \mathcal{C}_Q$ is a triangle functor.
    \item[(iii)] $\mathcal{C}_Q$ is Krull--Schmidt.
    \item[(iv)]
    \[
    \operatorname{ind}\mathcal{C}_Q
    = \operatorname{ind}(\operatorname{mod}\C Q)
    \sqcup \{\Sigma P_i \mid i \in Q_0\},
    \]
    where $\operatorname{ind}(-)$ is the set of isomorphism classes
    of indecomposable objects and $P_i$ is the indecomposable
    projective module associated to $i \in Q_0$.
\end{enumerate}
\end{Prop}

\begin{Def}[{\cite{BMRRT}}]
\label{Def : rigid and cluster tilting objects}
We say that $X \in \mathcal{C}_Q$ is \emph{rigid} if
$\operatorname{Ext}^1_{\mathcal{C}_Q}(X, X) = 0$.
We say that $T \in \mathcal{C}_Q$ is a \emph{cluster-tilting object}
if $T$ is rigid and $T$ has a maximal number of non-isomorphic
indecomposable direct summands.
An object $X \in \mathcal{C}_Q$ is called \emph{basic}
if all of its indecomposable direct summands are non-isomorphic.
\end{Def}

\begin{Thm}[{\cite[Theorem 3.3, Theorem 5.1]{BMRRT}}]
\label{Thm : cluster tilting summands and mutation}
\begin{enumerate}
    \item[(i)] Any cluster-tilting object has $|Q_0|$
    non-isomorphic indecomposable direct summands.
    \item[(ii)] For a basic cluster-tilting object
    $T = \bigoplus_{i \in Q_0}T_i$ and $k \in Q_0$,
    there exists a unique rigid indecomposable object
    $T_k^* \not\simeq T_k$ such that
    \[
    T_k^* \oplus \bigoplus_{i \in Q_0 \setminus \{k\}}T_i
    \]
    is a basic cluster-tilting object.
\end{enumerate}
\end{Thm}

For a basic cluster-tilting object
$T = \bigoplus_{i \in Q_0}T_i$ and $k \in Q_0$, we call
$T_k^* \oplus \bigoplus_{i \in Q_0 \setminus \{k\}}T_i$
the \emph{mutation} of $T$ at $k$.

\begin{Thm}[{\cite[Proposition 3.5]{BMRRT}}]
\label{Thm : additive reachability}
Let $Q$ be an acyclic quiver.
Then any basic cluster-tilting object $T$ is reachable,
i.e., $T$ is obtained from the basic cluster-tilting object
$\bigoplus_{i \in Q_0}P_i$ by finitely many mutations.
\end{Thm}

For $\mathbf{a} \in \Z_{\geq 0}^{Q_0}$, we set
\[
P(\mathbf{a}) \coloneq
\bigoplus_{i \in Q_0}P_i^{\oplus a_i}.
\]

\begin{Prop}[{\cite{Palu08}}]
\label{Prop : presentations and index}
For $X \in \mathcal{C}_Q$, there exists a triangle
\[
P(\mathbf{g}_1) \longrightarrow P(\mathbf{g}_0)
\longrightarrow X \xrightarrow{+1}
\]
for $\mathbf{g}_0, \mathbf{g}_1 \in \Z_{\geq 0}^{Q_0}$,
and $\mathbf{g}_0 - \mathbf{g}_1 \in \Z^{Q_0}$ is independent
of the choice of the presentation.
\end{Prop}

For $X \in \mathcal{C}_Q$, we define
$\operatorname{ind}X \coloneq
\mathbf{g}_0 - \mathbf{g}_1 \in \Z^{Q_0}$.

\begin{Def}[{\cite[\S2.2]{CK06}}]
\label{Def : quiver Grassmannians}
Let $M \in \operatorname{mod}\C Q$
and $\mathbf{d} \in \Z_{\geq 0}^{Q_0}$.
We define the \emph{quiver Grassmannian}
$\operatorname{Gr}_{\mathbf{d}}(M)$ by
\[
\{N \subset M \mid N \text{ is a submodule of } M
\text{ such that } \underline{\dim}N = \mathbf{d}\}.
\]
\end{Def}

\begin{Def}[{\cite{CK06}}]
\label{Def : cluster character}
Let $X \in \mathcal{C}_Q$.
We write $\Sigma X = M \oplus \Sigma P$
for $M \in \operatorname{mod}\C Q$ and a projective module $P$.
We define the \emph{cluster character}
$\mathrm{CC}(X) \in \Z[x_i^{\pm 1} \mid i \in Q_0]$ by
\[
\mathrm{CC}(X) \coloneq
\mathbf{x}^{\operatorname{ind}X}
\sum_{\mathbf{d} \in \Z_{\geq 0}^{Q_0}}
\chi\bigl(\operatorname{Gr}_{\mathbf{d}}(M)\bigr)
\mathbf{x}^{B_{Q^{\mathrm{op}}}\mathbf{d}},
\]
where $\chi$ is the Euler characteristic and $B_{Q^{\mathrm{op}}}$ is
the skew-symmetric matrix corresponding to $Q^{\mathrm{op}}$.
\end{Def}

\begin{Thm}[{\cite[Theorem 4]{CK06}}]
\label{Thm : cluster tilting objects and clusters}
Let $\overline{\mathcal{A}}_{Q^{\mathrm{op}}}$ be the cluster algebra with no frozen variables 
with the initial seed
\[
(\mathbf{x}_{t_0}, B_{t_0})
= \bigl((x_i)_{i \in Q_0}, B_{Q^{\mathrm{op}}}\bigr).
\]
Then
\[
\mathrm{CC} \colon \mathcal{C}_Q \longrightarrow
\Z[x_i^{\pm 1} \mid i \in Q_0]
\]
induces the bijection
\[
\{\text{indecomposable rigid objects in } \mathcal{C}_Q\}
\longrightarrow
\{\text{cluster variables of }
\overline{\mathcal{A}}_{Q^{\mathrm{op}}}\}.
\]
Moreover, it induces the bijection
\[
\{\text{basic cluster-tilting objects in } \mathcal{C}_Q\}
\longrightarrow
\{\text{clusters of }\overline{\mathcal{A}}_{Q^{\mathrm{op}}}\},
\]
and it is compatible with the mutations.
\end{Thm}

Note that for $X \in \mathcal{C}_Q$ such that
$\operatorname{ind}X = \mathbf{g}$, we have
\[
\deg^{t_0}\mathrm{CC}(X) = \mathbf{g}
\quad\text{in }\overline{\mathcal{A}}_{Q^{\mathrm{op}}}.
\]

\subsection{$F$-polynomials and canonical decompositions}
\label{ssec : F polynomials and canonical decompositions}

Let $Q$ be an acyclic quiver.
For $\mathbf{g} \in \Z^{Q_0}$, we set
\[
\operatorname{PHom}(\mathbf{g})
\coloneq \operatorname{Hom}_{\C Q}
\bigl(P([-\mathbf{g}]_+), P([\mathbf{g}]_+)\bigr).
\]
For $d \in \operatorname{PHom}(\mathbf{g})$, set
$X_d \coloneq \operatorname{Cone}(d)$.
Note that $\operatorname{ind}X_d = \mathbf{g}$
for $d \in \operatorname{PHom}(\mathbf{g})$.

\begin{Def}[{\cite[Definition 3.3]{DF15}}]
\label{Def : E invariant of presentations}
Let $\mathbf{g}_i \in \Z^{Q_0}$ and
$d_i \in \operatorname{PHom}(\mathbf{g}_i)$ for $i = 1, 2$.
We define
\[
E(d_1, d_2)
\coloneq
\operatorname{Hom}_{K^b(\operatorname{proj}\text{-}\C Q)}
(d_1, \Sigma d_2),
\]
where $K^b(\operatorname{proj}\text{-}\C Q)$ is the bounded homotopy
category of $\operatorname{proj}\text{-}\C Q$.
We set
\[
e(\mathbf{g}_1, \mathbf{g}_2)
\coloneq
\min\left\{
\dim E(d_1, d_2)
\;\middle|\;
 d_i \in \operatorname{PHom}(\mathbf{g}_i) \text{ for } i = 1, 2
\right\}.
\]
We define
\[
\operatorname{ext}(\mathbf{g}_1, \mathbf{g}_2)
\coloneq
\min\left\{
\dim\operatorname{Ext}^1_{\mathcal{C}_Q}(X_{d_1}, X_{d_2})
\;\middle|\;
 d_i \in \operatorname{PHom}(\mathbf{g}_i) \text{ for } i = 1, 2
\right\}.
\]
Since
\[
\dim E,\ \dim\operatorname{Ext}^1_{\mathcal{C}_Q}
\colon
\operatorname{PHom}(\mathbf{g}_1) \times \operatorname{PHom}(\mathbf{g}_2)
\longrightarrow \Z
\]
are upper semicontinuous \cite[Corollary 3.7]{DF15},
$e$ and $\operatorname{ext}$ are the generic values of
$\dim E$ and $\dim\operatorname{Ext}^1_{\mathcal{C}_Q}$.
\end{Def}

\begin{Prop}[{\cite[Lemma 3.3]{Palu08}}]
\label{Prop : generic cluster extensions and E invariants}
Let $\mathbf{g}_1, \mathbf{g}_2 \in \Z^{Q_0}$.
Then
\[
\operatorname{ext}(\mathbf{g}_1, \mathbf{g}_2)
= e(\mathbf{g}_1, \mathbf{g}_2) + e(\mathbf{g}_2, \mathbf{g}_1).
\]
\end{Prop}

\begin{Def}[{\cite[Definition 4.3]{DF15}}]
\label{Def : canonical decomposition of indices}
Let $\mathbf{g} \in \Z^{Q_0}$ and $d \in \operatorname{PHom}(\mathbf{g})$.
We say that $d$ is \emph{indecomposable} if $d$ is nonzero and
is not decomposable into a direct sum of nontrivial presentations.
We say that $\mathbf{g}$ is \emph{indecomposable} if a generic
$d \in \operatorname{PHom}(\mathbf{g})$ is indecomposable.
For $\mathbf{g} \in \Z^{Q_0}$,
$\mathbf{g} = \mathbf{g}_1 \oplus \cdots \oplus \mathbf{g}_s$ is a
\emph{canonical decomposition} if each $\mathbf{g}_i$ is indecomposable and
a generic element of $\operatorname{PHom}(\mathbf{g})$ has the form
$f_1 \oplus \cdots \oplus f_s$ with
$f_i \in \operatorname{PHom}(\mathbf{g}_i)$.
\end{Def}

\begin{Thm}[{\cite[Theorem 4.4]{DF15}}]
\label{Thm : characterization of canonical decompositions}
Let $\mathbf{g} \in \Z^{Q_0}$.
Then $\mathbf{g} = \mathbf{g}_1 \oplus \cdots \oplus \mathbf{g}_s$
is the canonical decomposition
if and only if each $\mathbf{g}_i$ is indecomposable, and
$e(\mathbf{g}_i, \mathbf{g}_j) = 0$ whenever $i \neq j$.
\end{Thm}

\begin{Prop}[{\cite[Corollary 21]{DW02}}]
\label{Prop : number of real and imaginary canonical summands}
Let $\{\mathbf{g}_1, \ldots, \mathbf{g}_r,
\mathbf{g}_{r+1}, \ldots, \mathbf{g}_{r+s}\}$ be a subset
of $\Z^{Q_0}$ satisfying the following conditions:
\begin{itemize}
    \item
    $\mathbf{g}_1 \oplus \cdots \oplus \mathbf{g}_{r+s}$
    is the canonical decomposition
    of $\sum_{i=1}^{r+s}\mathbf{g}_i$,
    \item
    each $\mathbf{g}_i$ $(1 \leq i \leq r)$ is an index of some
    indecomposable rigid object,
    \item
    each $\mathbf{g}_i$ $(r+1 \leq i \leq r+s)$ is not an index of any
    indecomposable rigid object.
\end{itemize}
Then we have $r + 2s \leq |Q_0|$.
\end{Prop}

\begin{Cor}
\label{Cor : full canonical decomposition and cluster tilting}
Let $\{\mathbf{g}_i \mid 1 \leq i \leq |Q_0|\}$ be a subset of $\Z^{Q_0}$
such that $\mathbf{g}_1 \oplus \cdots \oplus \mathbf{g}_{|Q_0|}$
is the canonical
decomposition of $\sum_{i=1}^{|Q_0|}\mathbf{g}_i$.
Then there exists a cluster-tilting object
$T = \bigoplus_{i=1}^{|Q_0|}T_i$ such that
$\operatorname{ind}T_i = \mathbf{g}_i$ for any $1 \leq i \leq |Q_0|$.
\end{Cor}

\begin{proof}
By Proposition~\ref{Prop : number of real and imaginary canonical summands},
each $\mathbf{g}_i$ is the index of some indecomposable rigid object $T_i$.
Since $T_i$ is rigid, $T_i$ is a generic object in
$\{X_{d_i} \mid d_i \in \operatorname{PHom}(\mathbf{g}_i)\}$.
Therefore,
\[
0 = \operatorname{ext}(\mathbf{g}_i, \mathbf{g}_j)
= \dim\operatorname{Ext}^1_{\mathcal{C}_Q}(T_i, T_j).
\]
Thus $T = \bigoplus_{i=1}^{|Q_0|}T_i$ is rigid and has $|Q_0|$
nonisomorphic direct summands, so $T$ is a basic cluster-tilting object.
\end{proof}

For $\Sigma X = M \oplus \Sigma P \in \mathcal{C}_Q$, we put
\[
F_X(z_i \mid i \in Q_0)
\coloneq
\sum_{d \in \Z_{\geq 0}^{Q_0}}
\chi\bigl(\operatorname{Gr}_d(M)\bigr)\mathbf{z}^d
\in \Z[z_i \mid i \in Q_0]
\]
and call it the \emph{$F$-polynomial} of $X$.

\begin{Prop}[{\cite[\S 3.4]{Plam13}}]
\label{Prop : generic cluster characters and F polynomials}
Let $\mathbf{g} \in \Z^{Q_0}$.
Then the map
\[
\operatorname{PHom}(\mathbf{g})
\longrightarrow \Z[x_i^{\pm 1} \mid i \in Q_0],
\qquad d \longmapsto \mathrm{CC}(X_d)
\]
has a generic value.
Moreover, the map
\[
\operatorname{PHom}(\mathbf{g})
\longrightarrow \Z[z_i \mid i \in Q_0],
\qquad d \longmapsto F_{X_d}(z_i \mid i \in Q_0)
\]
has a generic value.
We write it
$F_{\mathbf{g}}^{Q,\mathrm{gen}}(z_i \mid i \in Q_0)
\coloneq F_{\mathbf{g}}^{\mathrm{gen}}(z_i \mid i \in Q_0)$.
\end{Prop}

\begin{Prop}[{\cite{Schofield92}, \cite[Corollary 3.26]{Fei23}}]
\label{Prop : generic tropical F polynomials and E invariants}
For $\mathbf{g}_1, \mathbf{g}_2 \in \Z^{Q_0}$, we have
\[
F_{\mathbf{g}_1}^{\mathrm{gen}}[\mathbf{g}_2] = e(\mathbf{g}_1, \mathbf{g}_2).
\]
\end{Prop}

\begin{Rem}
\label{Rem : comparison with Fei notation}
Note that, in the notation of \cite{Fei23},
if $M = \operatorname{Ker}(\check{\delta})$, then
$f_M = F_{-\check{\delta}}^{\mathrm{gen}}[-]$.
\end{Rem}

\begin{Def}
\label{Def : Newton polytope}
For
\[
f(z_i \mid 1 \leq i \leq m)
= \sum_{\mathbf{v} \in \Z^m}c_{\mathbf{v}}\mathbf{z}^{\mathbf{v}}
\in \Z[z_i^{\pm 1} \mid 1 \leq i \leq m],
\]
the \emph{Newton polytope} of $f$ is the convex hull of
$\{\mathbf{v} \in \Z^m \mid c_{\mathbf{v}} \neq 0\}$ in $\mathbb{R}^m$.
We write the Newton polytope of $f$ as $\operatorname{Newt}(f)$.
\end{Def}

\begin{Lem}
\label{Lem : vertices of Newton polytopes}
Let $0 \neq F(z_i \mid i \in Q_0)
\in \Z[z_i^{\pm 1} \mid i \in Q_0]$
and $p \in \mathbb{R}^{Q_0}$.
The following are equivalent:
\begin{enumerate}
    \item[(i)]
    $p$ is a vertex of $\operatorname{Newt}(F)$.
    \item[(ii)]
    There exist $\{\mathbf{g}_i \in \Z^{Q_0} \mid i \in Q_0\}$ such that
    they are linearly independent, and
    \[
    F\left[\sum_{i \in Q_0}\mathbf{g}_i\right]
    = \sum_{i \in Q_0}F[\mathbf{g}_i],
    \qquad
    \{p\}
    = \bigcap_{i \in Q_0}
    \left\{\gamma \in \mathbb{R}^{Q_0}
    \;\middle|\; \mathbf{g}_i^{\mathsf{T}}\gamma = F[\mathbf{g}_i]\right\}.
    \]
\end{enumerate}
\end{Lem}

\begin{proof}
Assume (i).
Then the normal cone
\[
N_F(p)
= \left\{\gamma \in \mathbb{R}^{Q_0}
\;\middle|\;
\gamma\cdot(p-\mathbf{v}) \geq 0
\text{ for all } \mathbf{v} \in \operatorname{Newt}(F)\right\}
\]
is a $|Q_0|$-dimensional cone.
Therefore, there exist
$\{\mathbf{g}_i \in \Z^{Q_0} \mid i \in Q_0\} \subset N_F(p)$
such that they are linearly independent.
Then we have
\[
\begin{aligned}
F\left[\sum_{i \in Q_0}\mathbf{g}_i\right]
&= \max\left\{
\left(\sum_i \mathbf{g}_i\right)^{\mathsf{T}}\mathbf{v}
\;\middle|\; \mathbf{v} \in \operatorname{Newt}(F)\right\} \\
&= \left(\sum_i \mathbf{g}_i\right)^{\mathsf{T}}p.
\end{aligned}
\]
Therefore, we have
\[
F\left[\sum_i \mathbf{g}_i\right] = \sum_i F[\mathbf{g}_i],
\qquad
\{p\} = \bigcap_{i \in Q_0}
\left\{\gamma \in \mathbb{R}^{Q_0}
\;\middle|\; \mathbf{g}_i^{\mathsf{T}}\gamma = F[\mathbf{g}_i]\right\}.
\]
Assume (ii).
Let $p_0 \in \operatorname{Newt}(F)$ be a point which satisfies
\[
F\left[\sum_i \mathbf{g}_i\right]
= \left(\sum_i \mathbf{g}_i\right)^{\mathsf{T}}p_0.
\]
Then we have
\[
\begin{aligned}
0
&= \sum_i F[\mathbf{g}_i] - F\left[\sum_i \mathbf{g}_i\right] \\
&= \sum_i F[\mathbf{g}_i]
- \left(\sum_i \mathbf{g}_i\right)^{\mathsf{T}}p_0.
\end{aligned}
\]
Since $F[\mathbf{g}_i] \geq \mathbf{g}_i^{\mathsf{T}}p_0$ for all $i \in Q_0$,
we have $F[\mathbf{g}_i] = \mathbf{g}_i^{\mathsf{T}}p_0$ for all $i$.
Therefore,
\[
p_0 \in \bigcap_{i \in Q_0}
\left\{\gamma \in \mathbb{R}^{Q_0}
\;\middle|\; \mathbf{g}_i^{\mathsf{T}}\gamma = F[\mathbf{g}_i]\right\}.
\]
Because
\[
\left|
\bigcap_{i \in Q_0}
\left\{\gamma \in \mathbb{R}^{Q_0}
\;\middle|\; \mathbf{g}_i^{\mathsf{T}}\gamma = F[\mathbf{g}_i]\right\}
\right| = 1
\]
by the assumption, we have $p = p_0$.
Since
\[
\operatorname{Newt}(F)
\subset
\bigcap_{i \in Q_0}
\left\{\gamma \in \mathbb{R}^{Q_0}
\;\middle|\; \mathbf{g}_i^{\mathsf{T}}\gamma \leq F[\mathbf{g}_i]\right\}
\]
and $p \in \operatorname{Newt}(F)$ is a vertex of
\[
\bigcap_{i \in Q_0}
\left\{\gamma \in \mathbb{R}^{Q_0}
\;\middle|\; \mathbf{g}_i^{\mathsf{T}}\gamma \leq F[\mathbf{g}_i]\right\},
\]
$p$ is a vertex of $\operatorname{Newt}(F)$.
\end{proof}

\begin{Prop}[{\cite[Lemma 5.2, Theorem 6.4]{Fei23}}]
\label{Prop : cluster cones and generic Newton polytopes}
Let $\mathbf{g} \in \Z^{Q_0}$.
For a cluster-tilting object $T = \bigoplus_{i \in Q_0}T_i$,
\[
F_{\mathbf{g}}^{\mathrm{gen}}\left[\sum_i c_i\operatorname{ind}T_i\right]
= \sum_i c_i F_{\mathbf{g}}^{\mathrm{gen}}[\operatorname{ind}T_i]
\]
holds for all $c_i \in \Z_{\geq 0}$.
Moreover, for any vertex $p$ of
$\operatorname{Newt}(F_{\mathbf{g}}^{Q,\mathrm{gen}})$,
there exists some cluster-tilting object
$T = \bigoplus_i T_i$ such that
\[
\{p\} = \bigcap_{i \in Q_0}
\left\{\gamma \in \mathbb{R}^{Q_0}
\;\middle|\;
\gamma^{\mathsf{T}}(\operatorname{ind}T_i)
= F_{\mathbf{g}}^{\mathrm{gen}}[\operatorname{ind}T_i]\right\}.
\]
\end{Prop}

\subsection{Monoidal reachability for acyclic cluster algebras}
\label{ssec : monoidal reachability for acyclic cluster algebras}

Let
$\iota \colon K(\mathscr{C})_{\Z}
\xrightarrow{\sim}
\overline{\mathcal{A}}(\mathcal{S}_{\mathbf{x}}^{\Lambda}),J = 
\{1, \ldots, N\},J_{\mathrm{uf}},J_{\mathrm{fr}}$
be as in Subsection~\ref{ssec : Lambda cluster algebra structure}.
In this subsection, assume that
$\overline{\mathcal{A}}(\mathcal{S}_{\mathbf{x}}^{\Lambda})$
is of acyclic type.
We fix $t_1 \in \mathbb{T}$ such that the principal part $B_{t_1}$
of the exchange matrix $\widetilde{B}_{t_1}$ corresponds to
an acyclic quiver $Q^{\mathrm{op}}$.
Note that $J_{\mathrm{uf}} = Q_0$.
Since $t_1$ is injective-reachable, the following facts hold by
\cite[Theorem 1.3]{Plam13} and
\cite[Theorem 5.2.2]{Qin24}.

\begin{Lem}
\label{Lem : generic basis and compatible pointedness}
The set
\[
\left\{
L_{\widehat{\mathbf{g}}}
\coloneq \mathbf{x}_{t_1}^{\widehat{\mathbf{g}}}
F^{\mathrm{gen}}_{\mathbf{g}}
(\widehat{y}_{i;t_1} \mid i \in J_{\mathrm{uf}})
\;\middle|\; \widehat{\mathbf{g}} \in \Z^J
\right\}
\]
forms a basis of $\mathcal{U}(\mathcal{S}_{\mathbf{x}}^{\Lambda})$,
where $\mathbf{g} = \operatorname{pr}_{\mathrm{uf}}
\widehat{\mathbf{g}} \in \Z^{J_{\mathrm{uf}}} = \Z^{Q_0}$.
Moreover, each $L_{\widehat{\mathbf{g}}}$ is compatibly pointed and
$\deg^{t_1}L_{\widehat{\mathbf{g}}} = \widehat{\mathbf{g}}$.
\end{Lem}

This basis is called the \emph{generic basis} of
$\mathcal{U}(\mathcal{S}_{\mathbf{x}}^{\Lambda})$.
Note that the generic basis is not necessarily positive.

\begin{Prop}
\label{Prop : tropical F dominates generic tropical F}
Let $L_1, L_2 \in \Irr\mathscr{C}$ and
\[
\widehat{\mathbf{g}}_i
\coloneq \deg^{t_1}\iota([L_i]) \in \Z^J
\qquad\text{for }i = 1, 2.
\]
Let $\mathbf{g}_i
\coloneq \operatorname{pr}_{\mathrm{uf}}\widehat{\mathbf{g}}_i$
for $i = 1, 2$.
Then we have
\[
F^{t_1}_{\iota([L_1])}[\mathbf{g}_2]
\geq F^{Q,\mathrm{gen}}_{\mathbf{g}_1}[\mathbf{g}_2].
\]
\end{Prop}

\begin{proof}
It is enough to prove that
\[
\operatorname{Newt}(F^{t_1}_{\iota([L_1])})
\supset \operatorname{Newt}(F^{Q,\mathrm{gen}}_{\mathbf{g}_1}).
\]
Let us prove that each vertex $p$ of
$\operatorname{Newt}(F^{Q,\mathrm{gen}}_{\mathbf{g}_1})$
is a vertex of $\operatorname{Newt}(F^{t_1}_{\iota([L_1])})$.
By Proposition~\ref{Prop : cluster cones and generic Newton polytopes},
there exists a basic cluster-tilting object
$T = \bigoplus_{i \in Q_0}T_i$ such that
\[
F^{\mathrm{gen}}_{\mathbf{g}_1}
\left[\sum_i c_i\operatorname{ind}T_i\right]
= \sum_i c_i F^{\mathrm{gen}}_{\mathbf{g}_1}
[\operatorname{ind}T_i]
\qquad\text{for all }c_i \in \Z_{\geq 0},
\]
and
\[
\{p\} = \bigcap_{i \in Q_0}
\left\{\gamma \in \mathbb{R}^{Q_0}
\;\middle|\;
\gamma^{\mathsf{T}}\operatorname{ind}T_i
= F^{\mathrm{gen}}_{\mathbf{g}_1}[\operatorname{ind}T_i]
\right\}.
\]
By Lemma~\ref{Lem : vertices of Newton polytopes},
it is enough to prove that
\[
F^{t_1}_{\iota([L_1])}
\left[\sum_i c_i\operatorname{ind}T_i\right]
= F^{\mathrm{gen}}_{\mathbf{g}_1}
\left[\sum_i c_i\operatorname{ind}T_i\right]
\qquad\text{for all }c_i \in \Z_{\geq 0}.
\]
By Theorems~\ref{Thm : additive reachability} and
\ref{Thm : cluster tilting objects and clusters},
there exists $w \in \mathbb{T}$ such that
\[
\operatorname{ind}T_i
= \operatorname{pr}_{\mathrm{uf}}\deg^{t_1}x_{i;w}
\qquad\text{for }i \in J_{\mathrm{uf}}.
\]
Let
$\{\widehat{\mathbf{g}}_i^t \mid t \in \mathbb{T}\}
\in S_Y(\Z^{\max})$
be such that
$\widehat{\mathbf{g}}_i^{t_1} = \widehat{\mathbf{g}}_i$
for $i = 1, 2$.
By Proposition~\ref{Prop : tropical invariants of compatibly pointed elements},
Proposition~\ref{Prop : compatible pointedness of simple modules}, and
Lemma~\ref{Lem : generic basis and compatible pointedness}, we have
\[
\begin{aligned}
&\widehat{\mathbf{g}}_1^{\mathsf{T}}\Lambda_{t_1}
\bigl(\deg^{t_1}\mathbf{x}_{w}^{\mathbf{c}}\bigr)
+ 2F^{t_1}_{\iota([L_1])}
\left[\sum_i c_i\operatorname{ind}T_i\right] \\
&\qquad = (\widehat{\mathbf{g}}_1^{w})^{\mathsf{T}}
\Lambda_{w}
\bigl(\deg^{w}\mathbf{x}_{w}^{\mathbf{c}}\bigr)
+ 2F^{w}_{\iota([L_1])}
\left[\sum_i c_i\operatorname{pr}_{\mathrm{uf}}
\deg^{w}x_{i;w}\right]
\end{aligned}
\]
and
\[
\begin{aligned}
&\widehat{\mathbf{g}}_1^{\mathsf{T}}\Lambda_{t_1}
\bigl(\deg^{t_1}\mathbf{x}_{w}^{\mathbf{c}}\bigr)
+ 2F^{\mathrm{gen}}_{\mathbf{g}_1}
\left[\sum_i c_i\operatorname{ind}T_i\right] \\
&\qquad = (\widehat{\mathbf{g}}_1^{w})^{\mathsf{T}}
\Lambda_{w}
\bigl(\deg^{w}\mathbf{x}_{w}^{\mathbf{c}}\bigr)
+ 2F^{w}_{L_{\widehat{\mathbf{g}}_1}}
\left[\sum_i c_i\operatorname{pr}_{\mathrm{uf}}
\deg^{w}x_{i;w}\right].
\end{aligned}
\]
Since
\[
\sum_{i \in J_{\mathrm{uf}}}c_i\operatorname{pr}_{\mathrm{uf}}
\deg^{w}x_{i;w}
= \sum_{i \in J_{\mathrm{uf}}}c_i e_i
\in \Z_{\geq 0}^{J_{\mathrm{uf}}},
\]
we have
\[
F^{w}_{\iota([L_1])}
\left[\sum_i c_i\operatorname{pr}_{\mathrm{uf}}
\deg^{w}x_{i;w}\right]
= F^{w}_{L_{\widehat{\mathbf{g}}_1}}
\left[\sum_i c_i\operatorname{pr}_{\mathrm{uf}}
\deg^{w}x_{i;w}\right]
\]
by Proposition~\ref{Prop : compatible pointedness and bipointedness}.
Then we have
\[
F^{t_1}_{\iota([L_1])}
\left[\sum_i c_i\operatorname{ind}T_i\right]
= F^{\mathrm{gen}}_{\mathbf{g}_1}
\left[\sum_i c_i\operatorname{ind}T_i\right]
\qquad\text{for all }c_i \in \Z_{\geq 0}.
\]
\end{proof}

\begin{Cor}
\label{Cor : d dominates generic extensions}
Let $L_1, L_2 \in \Irr\mathscr{C}$.
Let $\widehat{\mathbf{g}}_i
\coloneq \deg^{t_1}\iota([L_i])$ and
$\mathbf{g}_i
\coloneq \operatorname{pr}_{\mathrm{uf}}\widehat{\mathbf{g}}_i$.
Then we have
\[
\mathfrak{d}(L_1, L_2)
\geq e(\mathbf{g}_1, \mathbf{g}_2)
+ e(\mathbf{g}_2, \mathbf{g}_1).
\]
\end{Cor}

\begin{proof}
By Theorem~\ref{Thm : tropical F invariants and Lambda inequalities},
Proposition~\ref{Prop : tropical F dominates generic tropical F}, and
Proposition~\ref{Prop : generic tropical F polynomials and E invariants}, we
have
\[
\begin{aligned}
\mathfrak{d}(L_1, L_2)
&\geq F^{t_1}_{\iota([L_1])}[\mathbf{g}_2]
+ F^{t_1}_{\iota([L_2])}[\mathbf{g}_1] \\
&\geq F^{\mathrm{gen}}_{\mathbf{g}_1}[\mathbf{g}_2]
+ F^{\mathrm{gen}}_{\mathbf{g}_2}[\mathbf{g}_1] \\
&= e(\mathbf{g}_1, \mathbf{g}_2)
+ e(\mathbf{g}_2, \mathbf{g}_1).
\end{aligned}
\]
\end{proof}

\begin{Ex}
\label{Ex : positive self d invariant}
Let
\[
A = \begin{pmatrix}2 & -n \\ -n & 2\end{pmatrix}
\qquad\text{for }n \in \Z_{\geq 3}.
\]
Consider the corresponding quiver Hecke algebra $R$ and let
$\widetilde{w} \coloneq s_1s_2s_1s_2$.
Let $\mathscr{C} = \mathscr{C}_w$ and
\[
\iota \colon K(\mathscr{C})_{\Z}
\xrightarrow{\sim}
\overline{\mathcal{A}}(\mathcal{S}_{\mathbf{x}}^{\Lambda})
\]
be the monoidal categorification as in
Subsection~\ref{ssec : Lambda cluster algebra structure}.
Here, $J = \{1, 2, 3, 4\}$, $J_{\mathrm{uf}} = \{1, 2\}$,
and $J_{\mathrm{fr}} = \{3, 4\}$.
Let $t_1 = t_0$.
Then $Q^{\mathrm{op}} = (1 \xrightarrow{n} 2)$
and $Q = (1 \xleftarrow{n} 2)$.
Note that $\dim P_1 = (1, 0)$ and $\dim P_2 = (n, 1)$
for $P_1, P_2 \in \operatorname{mod}\C Q$.
Let $\mathbf{g} = (-1, n-1)^{\mathsf{T}}$.
We can identify
\[
\operatorname{PHom}(\mathbf{g})
= \operatorname{Hom}_{\C Q}(P_1, P_2^{\oplus(n-1)})
\simeq M_{(n-1)\times n}(\C).
\]
For $d_1, d_2 \in \operatorname{PHom}(\mathbf{g})$
corresponding to $A, B \in M_{(n-1)\times n}(\C)$,
\[
\dim E(d_1, d_2)
= \dim\operatorname{Hom}_{\C Q}(P_1, P_2^{\oplus(n-1)})
- \dim\operatorname{Im}\delta_{A,B},
\]
where
\[
\delta_{A,B} \colon
\C \times M_{(n-1)\times(n-1)}(\C)
\longrightarrow M_{(n-1)\times n}(\C);
\qquad (a, C) \longmapsto Ba - CA.
\]
Therefore
\[
\dim E(d_1, d_2)
\geq n(n-1)-\bigl(1+(n-1)^2\bigr)
= n-2.
\]
Since $\delta_{A,B}$ is injective for generic $A, B$, we have
$e(\mathbf{g}, \mathbf{g}) = n-2$.
Let $L \in \Irr\mathscr{C}$ be such that
$\operatorname{pr}_{\mathrm{uf}}\deg^{t_1}\iota([L]) = \mathbf{g}$.
For example, such an $L$ appears as a composition factor of
$M_{1;\mu_1t_1} \otimes M_{2;\mu_2t_1}$.
Then $\mathfrak{d}(L, L) \geq 2(n-2) > 0$
by Corollary~\ref{Cor : d dominates generic extensions}.
In particular, there is a simple module whose self $\mathfrak{d}$-invariant
is non-zero.
\end{Ex}

\begin{Thm}
\label{Thm : monoidal reachability for acyclic cluster algebras}
For $L \in \Irr\mathscr{C}$,
$L$ is real if and only if $L$ is reachable.
Moreover, $L$ is real prime if and only if
$L$ corresponds to a cluster variable.
\end{Thm}

\begin{proof}
Let $\mathcal{M} = \{L_i \mid i \in J\}$ be a maximal real set
such that $L$ is a monomial in $\mathcal{M}$.
Let $\widehat{\mathbf{g}}_i \coloneq \deg^{t_1}\iota([L_i])$
and $\mathbf{g}_i
\coloneq \operatorname{pr}_{\mathrm{uf}}\widehat{\mathbf{g}}_i$
for $i \in J$.
Note that $\widehat{\mathbf{g}}_i = e_i$ for $i \in J_{\mathrm{fr}}$.
Let $\{\delta_1, \ldots, \delta_s\} \subset \Z^{Q_0}$
be the set of indecomposable indices which appear in
the canonical decomposition of some $\mathbf{g}_i$
for $i \in J_{\mathrm{uf}}$.
By Corollary~\ref{Cor : d dominates generic extensions}, we have
$e(\delta_i, \delta_j) = 0$ for $1 \leq i, j \leq s$.
Since
$\{\mathbf{g}_i \mid i \in J_{\mathrm{uf}}\}
\subset \Z^{J_{\mathrm{uf}}}$
is linearly independent by
Theorem~\ref{Thm : Laurent families and maximal real sets} and
Proposition~\ref{Prop : degrees of Laurent families form an integral basis},
we obtain that $s \geq |Q_0|$.
By Corollary~\ref{Cor : full canonical decomposition and cluster tilting},
$s = |Q_0|$ and there exists a basic cluster-tilting object
$T = \bigoplus_{i \in Q_0}T_i$ such that
$\delta_i = \operatorname{ind}T_i$ for $i \in J_{\mathrm{uf}}$.
By Theorem~\ref{Thm : additive reachability},
there exist $w \in \mathbb{T}$ and a permutation $\sigma$
of $J_{\mathrm{uf}}$ such that
\[
\operatorname{pr}_{\mathrm{uf}}\deg^{t_1}x_{\sigma(i);w}
= \delta_i.
\]
Let $\mathbf{g}_i = \sum_{j \in J_{\mathrm{uf}}}a_{ij}\delta_j$
for $i \in J_{\mathrm{uf}}$ and $a_{ij} \in \Z_{\geq 0}$.
Since
\[
\operatorname{pr}_{\mathrm{uf}}\deg^{t_1}\iota([L_i])
= \operatorname{pr}_{\mathrm{uf}}\deg^{t_1}
\iota\left(\left[
\bigotimes_{k \in J_{\mathrm{uf}}}
M_{\sigma(k);w}^{\otimes a_{ik}}
\right]\right),
\]
there exist $\mathbf{a}, \mathbf{b}
\in \Z_{\geq 0}^{J_{\mathrm{fr}}}$ such that
\[
\begin{aligned}
&\deg^{t_1}\iota\left(\left[
L_i \otimes \bigotimes_{j \in J_{\mathrm{fr}}}
M_{j;t_1}^{\otimes a_j}\right]\right) \\
&\qquad = \deg^{t_1}\iota\left(\left[
\bigotimes_{k \in J_{\mathrm{uf}}}
M_{\sigma(k);w}^{\otimes a_{ik}}
\otimes \bigotimes_{j \in J_{\mathrm{fr}}}
M_{j;t_1}^{\otimes b_j}\right]\right).
\end{aligned}
\]
By Lemma~\ref{Lem : injectivity of degrees for Laurent families},
we have
\[
L_i \otimes \bigotimes_{j \in J_{\mathrm{fr}}}
M_{j;t_1}^{\otimes a_j}
\simeq \bigotimes_{k \in J_{\mathrm{uf}}}
M_{\sigma(k);w}^{\otimes a_{ik}}
\otimes \bigotimes_{j \in J_{\mathrm{fr}}}
M_{j;t_1}^{\otimes b_j}.
\]
By Corollary~\ref{Cor : unique prime tensor factorization in C},
we have
\[
L_i \simeq \bigotimes_{k \in J_{\mathrm{uf}}}
M_{\sigma(k);w}^{\otimes a_{ik}}.
\]
Therefore, $L$ is a monomial in a monoidal seed
$\{M_{i;w} \mid i \in J\}$.
In particular, if $L$ is real prime,
$L$ is isomorphic to some $M_{i;w}$ by primeness.
\end{proof}

\bibliographystyle{amsplain}
\bibliography{refs}

\end{document}